\documentclass[11pt]{article} 

\usepackage{amsmath}
\usepackage{amssymb}
\usepackage[T1]{fontenc}
\usepackage[latin1]{inputenc}
\def\r{\rightarrow}

\newcommand{\fdem}{\hspace*{\fill}~$\Box$\par\endtrivlist\unskip}

\newcommand{\E}{\mathbb{E}}     
\renewcommand{\P}{\mathbb{P}}     
\renewcommand{\L}{\mathbb{L}}

\newcommand{\N}{\mathbb{N}}     
\newcommand{\Z}{\mathbb{Z}}
\newcommand{\R}{\mathbb{R}}     
     
\newcommand{\C}{\mathbb{C}}

\renewcommand{\r}{\mathop{\rightarrow}}

\newcommand{\cB}{\mbox{$\cal B$}}
\newcommand{\cC}{\mbox{$\cal C$}}
\newcommand{\cD}{\mbox{$\cal D$}}
\newcommand{\cE}{\mbox{$\cal E$}}
\newcommand{\cF}{\mbox{$\cal F$}}
\newcommand{\cG}{\mbox{$\cal G$}}

\newcommand{\cL}{\mbox{$\cal L$}}
\newcommand{\cM}{\mbox{$\cal M$}}
\newcommand{\cN}{\mbox{$\cal N$}}

\newcommand{\cR}{\mbox{$\cal R$}}

\newcommand{\cU}{\mbox{$\cal U$}}
\newcommand{\cV}{\mbox{$\cal V$}}

\begin{document}

\begin{center} {\large\bf THE NAGAEV-GUIVARC'H METHOD 
VIA THE KELLER-LIVERANI THEOREM}  
\end{center}
\vskip 3mm
\begin{center}
Lo\"{\i}c HERV\'E \footnote{
Universit\'e
Europ\'eenne de Bretagne, I.R.M.A.R. (UMR-CNRS 6625), Institut National des Sciences 
Appliqu\'ees de Rennes. Loic.Herve@insa-rennes.fr} 
and Fran\c{c}oise P\`ENE \footnote{Universit\'e
Europ\'eenne de Bretagne, Universit\'e de Brest,
laboratoire de Math\'ematiques (UMR CNRS 6205). 
francoise.pene@univ-brest.fr}
\end{center}



\vskip 3mm

\noindent{\bf Abstract.} {\scriptsize \it The Nagaev-Guivarc'h method, via the perturbation 
operator theorem of Keller and Liverani, has been exploited in recent papers to establish limit theorems for unbounded functionals of strongly ergodic Markov chains. 
The main difficulty of this approach is to prove Taylor expansions for the dominating eigenvalue of the Fourier kernels. 
The paper outlines this method and extends it by stating a multidimensional local limit theorem, 
a one-dimensional Berry-Esseen theorem, a first-order Edgeworth expansion, and a multidimensional Berry-Esseen type theorem in the sense of the Prohorov metric. When applied to the exponentially $\L^2$-convergent Markov chains, to the $v$-geometrically ergodic Markov chains and to the iterative  Lipschitz models, the three first above cited limit theorems hold  under moment conditions similar, or close (up to $\varepsilon>0$), to those of the i.i.d.~case.}
\vskip 10mm

{\bf Titre (en français).} {La méthode de Nagaev-Guivarc'h via le théorème de Keller-Liverani}

Mots-clés : chaîne de Markov, théorème limite central, développement d'Edgeworth, méthode spectrale

\noindent{\bf Résumé.} {\scriptsize \it La méthode de Nagaev-Guivarc'h, via le théorème de perturbation de Keller et Liverani, a été appliquée récemment en vu d'établir des théorèmes limites pour des fonctionnelles non bornées de chaînes de Markov fortement ergodiques. La difficulté principale dans cette approche est de démontrer des développements de Taylor pour la valeur propre perturbée de l'opérateur de Fourier. Dans ce travail, nous donnons une présentation générale de cette méthode, et nous l'étendons en démontrant un théorème limite local multidimensionnel, un théorème de Berry-Esseen unidimensionnel, un développement d'Edgeworth d'ordre 1, et enfin un théorème de Berry-Esseen multidimensionnel au sens de la distance de Prohorov. Nos applications concernent les chaînes de Markov $\L^2$-fortement ergodiques, $v$-géométriquement ergodiques, et les modèles itératifs. Pour ces exemples, les trois premiers théorèmes limites cités précédemment sont satisfaits sous des conditions de moment dont l'ordre est le même (parfois à $\varepsilon>0$ près) que dans le cas indépendant. }
\vskip 10mm

\indent $\ \ \ \ \ \ \ \ \ \ \ \ \ \ \ \ \ \ \ \ \ \ \ \ \ \ \ \ \ \ \ \ \ \ \ \ \ \ \ \ \  $ {\bf CONTENTS} \\

\noindent 1. Introduction, setting and notations  \\[0.15cm]
\noindent 2. A central limit theorem in the stationary case \\[0.15cm]
\noindent 3. The usual Nagaev-Guivarc'h method \\[0.15cm]
\noindent 4. The Nagaev-Guivarc'h method via the Keller-Liverani theorem (Conditions ($\widetilde{K}$) (K)) \\[0.15cm]
\noindent 5. A multidimensional local limit theorem (Conditions (S) ($\widehat K$)) \\[0.15cm]
\noindent 6. A one-dimensional uniform Berry-Esseen theorem \\[0.15cm]
\noindent 7. Regularity of the eigen-elements of the Fourier kernels (condition ${\cC}(m)$) \\[0.15cm]
\noindent 8. A one-dimensional first-order Edgeworth expansion \\[0.15cm]
\noindent 9. A multidimensional Berry-Esseen theorem (for the Prohorov metric) \\[0.15cm] 
\noindent 10. Application to $v$-geometrically ergodic Markov chains \\[0.15cm]
\noindent 11. Application to iterative Lipschitz models \\[0.15cm]
\noindent 12. More on non-arithmeticity and nonlattice conditions \\[0.15cm]
\noindent Appendix A: proof of Proposition 7.1 \\[0.15cm]
\noindent Appendix B: proof of Propositions 11.4-8.

\newpage 

\section{Introduction, setting and notations} 
\noindent Let $(X_n)_n$ be a Markov chain with values in 
$(E,{\mathcal E})$, with transition probability $Q$ and
with stationary distribution $\pi$. 
Let $\xi$ be a $\pi$-centered random 
variable with values in ${\R}^d$ (with $d\ge 1$).
We are interested in probabilistic limit theorems for
$(\xi(X_n))_n$ namely: \\[0.15cm]
\indent {\scriptsize $\bullet$} central limit theorem (c.l.t.), \\[0.15cm]
\indent {\scriptsize $\bullet$} rate of convergence in the central limit theorem: 
Berry Esseen type theorem, \\[0.15cm]
\indent {\scriptsize $\bullet$} multidimensional local limit theorem, \\[0.15cm]
\indent {\scriptsize $\bullet$} First-order
Edgeworth expansion (when $d=1$). \\[0.15cm]
We want to establish these results under moment
conditions on $\xi$ as close as possible to those of the
i.i.d.~case (as usual i.i.d.~is the short-hand for 
``independent and identically distributed''). 
Let us recall some facts about the case when $(Y_n)_n$ is a sequence of i.i.d.~$\R^d$-valued random variables (r.v.) with null expectation.
If $Y_1\in \L^2$, we have the central limit 
theorem and, under
some additional nonlattice type assumption, we have the local limit theorem.
If $Y_1\in \L^3$ and $d=1$, we have the uniform Berry-Esseen theorem, and the
first-order Edgeworth expansion (under the nonlattice assumption).  
All these results can be proved thanks to Fourier
techniques.
If $Y_1\in\L^3$, 
$(Y_n)_n$ satisfies a multidimensional Berry-Esseen 
type theorem 
(in the sense of the Prohorov metric). 
The proof of this last result uses Fourier techniques
and a truncation argument. \\[0.1cm]
\indent To get analogous results for Markov chains, 
we shall use and adapt the Nagaev-Guivarc'h method, introduced 
in \cite{nag1} \cite{nag2} in the case $d=1$.
This method is based on Fourier techniques and on the usual perturbation 
operator theory applied to the Fourier kernels 
$Q(t)(x,dy) = e^{it\xi(y)}Q(x,dy)$ ($t\in\R$). 
The idea is that ${\E}\big[e^{it\sum_{k=1}^n\xi(X_k)}\big]$ is close enough to
an expression of the form $\lambda(t)^n$, and the calculations are then similar to those of the
i.i.d.~case. Indeed, let us recall that, if $(Y_n)_n$ is a sequence
of i.i.d.~random variables, then we have ${\E}\big[e^{it\sum_{k=1}^nY_k}\big] = \big(\E[e^{itY_1}]\big)^n$.  \\[0.1cm]
\indent The Nagaev-Guivarc'h method, also called the spectral method, has been widely strengthened and extended, especially since the 80's with the contribution of Le Page \cite{lep82}, Rousseau-Egele \cite{rou}, Guivarc'h \cite{gui}, Guivarc'h and Hardy \cite{guihar}, Milhaud and Raugi \cite{mira}. This is fully described by Hennion and the first author in \cite{hulo}, where other references are given. Roughly speaking, to operate the Nagaev-Guivarc'h method, one needs the following strong ergodicity assumption (specified below) w.r.t.~some Banach space $\cB$, namely: $Q^n\r\pi$ in the operator norm topology of $\cB$. Under this assumption, the sequence $(\xi(X_n))_n$ then satisfies the usual distributional limit theorems provided that $(Q,\xi)$ verifies some operator-moment conditions on $\cB$.  
This method is especially efficient when $\cB$ is a Banach algebra and $\xi$ is in $\cB$. 
Unfortunately, on the one hand, since Banach algebras are often composed of bounded functions, the condition 
$\xi\in \cB$ implies that $\xi$ must be bounded. On the other hand, usual models 
as $v$-geometrically ergodic Markov chains or iterative Lipschitz models (typically $E=\R^p$) are strongly ergodic w.r.t.
some 
weighted supremum normed space or weighted Lipschitz-type space which are not Banach algebras, and the above mentioned 
operator-moment conditions then hold under very restrictive assumptions involving both $Q$ and $\xi$. For instance, in these models, 
the usual  spectral method cannot be efficiently applied to the sequence $(X_n)_n$ (i.e.~$\xi(x)=x$); an explicit and typical counter-example will be presented in Section 3. \\[0.1cm]
\indent In recent works  \cite{aap,ihp1,chazottes-gouezel,gui-lepage,ihp2,guibourg,seb-08}, a new procedure, based on 
the perturbation theorem of Keller-Liverani \cite{keli} (see also \cite{bal} p.~177), 
allows to get round the previous difficulty and to greatly improve the 
Nagaev-Guivarc'h method when applied to unbounded functionals $\xi$. 
Our work outlines this new approach, and presents the applications, namely~:
a multidimensional local limit theorem, a one-dimensional Berry Esseen theorem, 
a first-order Edgeworth expansion.
We establish these results under hypotheses close
to the i.i.d. case. 
We also establish a multidimensional Berry-Esseen
type theorem in the sense of the Prohorov metric
under hypotheses analogous to $Y_1\in{\mathbb L}^m$
with $m=\max\left(3,\lfloor{d/2}\rfloor+1\right)$
instead of $Y_1\in{\mathbb L}^3$. The reason is that,
when adapting \cite{Yurinskii},
we can use Yurinskii's smoothing inequality
(valid for r.v.~in ${\mathbb L}^m$) but
we cannot adapt Yurinskii's truncation argument.\\[0.1cm]
\indent When the usual perturbation theorem is replaced with that 
of Keller-Liverani, the main difficulty consists in proving Taylor expansions for the dominating eigenvalue $\lambda(t)$ 
of the Fourier kernel $Q(t)$. 
This point is crucial here.
Such expansions may be obtained as follows: \\[0.1cm]
(A) To get Taylor expansion at $t=0$, 
one can combine the spectral method with more probabilistic arguments such as martingale techniques \cite{ihp2}. 
In this paper, this method is just  
outlined: the local limit theorem obtained in \cite{ihp1} is extended 
to the multidimensional case, and the one-dimensional uniform Berry-Esseen theorem of \cite{ihp2} 
is here just recalled for completeness. \\[0.1cm]
(B) To establish the others limit theorems, we shall use a stronger property: the 
regularity of the eigen-elements of $Q(\cdot)$ on a neighbourhood of $t=0$.
We shall see that this can be done by considering the action of $Q(t)$ on 
a ``chain'' of suitable Banach spaces instead of a single 
one as in the classical approach. This method,  already used for other purposes in \cite{lep89,hen1,gouliv}, has been introduced in the spectral method \cite{aap} to investigate the c.l.t.~for iterative Lipschitz models. It is here specified and extended to general strongly ergodic Markov chains, and it will provide the one-dimensional Edgeworth expansion and the multidimensional Berry-Esseen type theorem.  

\indent Next, we  introduce our probabilistic setting, and the functional notations and definitions, helpful in defining the operator-type procedures of the next sections. 

\noindent {\bf Probabilistic setting.} $(X_n)_{n\geq0}$ is a Markov chain with 
general state space $(E,\cE)$, transition probability $Q$, stationary distribution $\pi$, initial distribution $\mu$, and 
$\xi = (\xi_1,\ldots,\xi_d)$ is a $\R^d$-valued $\pi$-integrable function on $E$ such that $\pi(\xi) = 0$ (i.e.~the $\xi_i$'s are 
$\pi$-integrable and  $\pi(\xi_i) = 0$). The associated random walk in $\R^d$ is denoted by \\
\indent $\displaystyle \ \ \ \ \ \ \ \ \ \ \ \ \ \ \ \ \ \ \ \ \ \ \ \ \ \ \ \ \ \ \ \ \ \ \ \ \ \ \ \ 
S_n = \sum_{k=1}^{n} \xi(X_k)$. \\
\noindent We denote by $|\cdot|_2$ and $\langle\cdot,\cdot\rangle$ the euclidean norm and the canonical scalar product on 
$\R^d$. For any $t\in\R^d$ and $x\in E$, we define the Fourier kernels of ($Q,\xi$) as 
$$Q(t)(x,dy) = e^{i\langle t,\, \xi(y)\rangle}Q(x,dy).$$
$\cN(0,\Gamma)$ denotes the centered normal distribution associated to 
a covariance matrix $\Gamma$, and ``$^{_{\underline{\ \, \tiny{\cD}\, \ }_>}}$'' means ``convergence in distribution''. 
Although $(X_n)_{n\geq0}$ is not a priori the canonical version, we shall slightly abuse notation and write $\P_\mu$, $\E_\mu$ to refer 
to the initial distribution. For any $\mu$-integrable function $f$, we shall often write $\mu(f)$ for $\int fd\mu$. For $x\in E$, 
$\delta_x$ will stand for the Dirac mass: $\delta_x(f) = f(x)$.  Finally, a set $A\in\cE$ is said to be $\pi$-full if $\pi(A)=1$, and $Q$-absorbing if $Q(a,A)=1$ for all $a\in A$.  

\noindent {\bf Functional setting.}  Let $\cB,X$ be complex Banach spaces. We denote by $\cL(\cB,X)$ the space of the bounded linear operators from $\cB$ to $X$, and by $\|\cdot\|_{\scriptsize{\cB,X}}$ the associated operator 
norm, with the usual simplified notations $\cL(\cB) = \cL(\cB,\cB)$, $\cB' = \cL(\cB,\C)$, for which the 
associated norms are simply denoted by $\|\cdot\|_{\scriptsize \cB}$. If $T\in\cL(\cB)$, $r(T)$ denotes its spectral radius, and $r_{ess}(T)$ its essential spectral radius. For the next use of the notion of essential spectral radius, we refer for instance to \cite{hen2,rab-wolf,wu} and \cite[Chap.~XIV]{hulo}. The notation ``$\cB\hookrightarrow X$'' means that $\cB\subset X$ and that the identity map is continuous from $\cB$ into $X$. \\[0.1cm]
We denote by $\cL^1(\pi)$ the vector space of the complex-valued $\pi$-integrable functions on $E$, and by $Cl(f)$ the class of $f$ modulo $\pi$. We call $\cB^{^\infty}$ the space of all bounded measurable functions on $E$ equipped with the supremum norm, and $\L^p(\pi)$, $1\leq p\leq+\infty$, the usual Lebesgue space. If $\cB\subset\cL^1(\pi)$ and $X\subset\L^1(\pi)$, we shall also use the notation ``$\cB\hookrightarrow X$'' to express that we have $Cl(f)\in X$ for all $f\in\cB$ and that the map $f\mapsto Cl(f)$ is continuous from $\cB$ to $X$. \\[0.1cm]
\noindent If $f\in\cL^1(\pi)$, it can be easily seen that the following function \\[0.2cm]
\indent $\displaystyle ({\cal Q})\ \ \ \ \ \ \ \ \ \ \ \ \ \ \ \ \ \ \ \ \ \ \ \ \ \ (Qf)(x) = \int_E f(y)\, Q(x,dy)$ \\[0.2cm] 
is defined $\pi$-a.s.~and is $\pi$-integrable with: $\pi(|Qf|) \leq \pi(|f|)$. If $\cB\subset\cL^1(\pi)$, $Q(\cB)\subset \cB$ and $Q\in\cL(\cB)$, we say that $Q$ continuously acts on $\cB$. If $\cB\subset\L^1(\pi)$, we shall use the same definition with $Q$ given by $Q\big(Cl(f)\big) = Cl(Qf)$ (which is possible since $Cl(f) = Cl(g)$ implies $Cl(Qf) = Cl(Qg)$). Clearly, $Q$ is a contraction on $\cB^{^\infty}$ and $\L^p(\pi)$. 

\noindent {\bf Strong ergodicity assumption.} {\it Unless otherwise indicated, all the normed spaces $(\cB,\|\cdot\|_{\scriptsize \cB})$ considered in this paper satisfy the following assumptions: $(\cB,\|\cdot\|_{\scriptsize \cB})$ is a Banach space such that, either $\cB\subset\cL^1(\pi)$ and $1_E\in\cB$, or $\cB\subset\L^1(\pi)$ and $Cl(1_E)\in\cB$, and we have in both cases $\cB\hookrightarrow \L^1(\pi)$. We then have $\pi\in\cB'$, so we can define the rank-one projection $\Pi$ on $\cB$:  
$$\Pi f = \pi(f)1_E\ \ (f\in\cB),$$
and we shall say that $Q$ (or merely $(X_n)_n$) is {\it strongly ergodic w.r.t.~$\cB$} if the following holds: 

\noindent (K1) $\ \ Q\in\cL(\cB)$ and $\ \lim_n\|Q^n-\Pi\|_{\scriptsize \cB} = 0$.}

\noindent One could also say `` geometrically ergodic w.r.t.~$\cB$ ``. Indeed, one can easily see that the last property in (K1) is equivalent to: \\[0.2cm]
\noindent (K'1) $\displaystyle \ \ \exists \kappa_0<1,\, \exists C>0,\, \forall n\geq1,\ \|Q^n-\Pi\|_{\scriptsize \cB} \leq C\, \kappa_0^{n}$.  \\[0.2cm]
\noindent We shall repeatedly use the following obvious fact. If $Q$ is  strongly ergodic w.r.t.~$\cB$, and if $f\in\cB$ 
is such that $\pi(f)=0$, then the series $\sum_{k\geq0}Q^kf$ is absolutely convergent in $\cB$. 

\indent Now, let us return to more probabilistic facts. When $(X_n)_n$ is Harris recurrent and strongly mixing, the so-called regenerative (or splitting) method provides limit theorems,  including the uniform Berry-Esseen theorem \cite{bolt} and Edgeworth expansions \cite{malin}. We want to point out that here the Harris recurrence is not assumed a priori. Moreover, the Markov chains in Examples~1-2 below are strongly mixing, but for these two examples, our results will be as efficient as all the others hitherto known ones, even better in many cases. The random iterative models of Example 3 are not automatically, either strongly mixing, or even Harris recurrent (see \cite{als}). 

\noindent {\bf Example 1: The strongly ergodic Markov chains on $\L^2(\pi)$} (see e.g. \cite{rosen}).  
We assume here that the $\sigma$-algebra $\cE$ is countably generated. 
%, and that $(X_n)_n$ is 
Let us recall that the strong ergodicity property on $\L^2(\pi)$ (namely, (K1) on $\cB=\L^2(\pi)$) implies that (K1) holds on $\L^p(\pi)$ for any $p\in(1,+\infty)$, see \cite{rosen}. This assumption, introduced in \cite{rosen} and called the exponential $\L^2(\pi)$-convergence in the literature, corresponds to ergodic and aperiodic Markov chains with spectral gap on $\L^2(\pi)$, see for instance the recent works \cite{wu,gong-wu} (and the references therein). \\[0.1cm]
The previous assumption is for instance satisfied if we have (K1) on $\cB^{^\infty}$ (see \cite{rosen}): in this case, according to the terminology of \cite{mey}, we will say that $(X_n)_n$ is {\it uniformly ergodic}. Equivalently, $(X_n)_n$ is aperiodic, ergodic, and satisfies the so-called Doeblin condition, see \cite{rosen}. This simple example was used in Nagaev's works \cite{nag1,nag2} (see Section 3). \\[0.1cm]
The strong ergodicity on $\L^2(\pi)$ provides a first motivation and a good understanding of the present improvements. Indeed, (except for the multidimensional
Berry-Esseen theorem) for results requiring $Y_1\in\L^m$ in the i.i.d.~case, whereas the usual Nagaev-Guivarc'h method needs the assumption  $\sup_{x\in E}\int |\xi(y)|^m Q(x,dy) < +\infty$ \cite{nag2,gharib, datta}, the present method appeals to the moment conditions $\xi\in\L^m(\pi)$ or $\xi\in\L^{m+\varepsilon}(\pi)$. \\[0.1cm]
\noindent In more concrete terms, let $(X_n)_n$ be a strongly ergodic Markov chain on $\L^2(\pi)$, and for convenience let us assume that $(X_n)_n$ is stationary (i.e.~$\mu=\pi$). From Gordin's theorem (Section~2), if $\pi(|\xi|_2^2)<+\infty$, then $(S_n/\sqrt n)_n$ converges in distribution to a normal law $\cN(0,\Gamma)$ (see also \cite{chen,jones}). It is understood below that the covariance matrix $\Gamma$ is invertible. The nonlattice condition will mean that the following property is fulfilled: there is no $a\in\R^d$, no closed subgroup $H$ in $\R^d$, $H\neq \R^d$, 
no $\pi$-full $Q$-absorbing set $A\in\cE$, and finally no bounded measurable function $\theta\, :\, E\r\R^d$ such that:  
$\forall x\in A, \ \ \xi(y)+\theta(y)-\theta(x)\in a+H\ \ Q(x,dy)-a.s.$.  \\[0.1cm]
\noindent The next statements, that will be specified  and established as corollaries of the abstract results of 
Sections 5-9, are new to our knowledge. Some further details and comparisons with prior results will be presented together with the corollaries cited below: 

\noindent{\it 
%{\scriptsize $\bullet$} 
(a) If $\pi(|\xi|_2^2)<+\infty$  and $\xi$ is nonlattice, then 
$(\xi(X_n))_n$ satisfies a multidimensional local limit theorem} (Corollary 5.5). \\[0.1cm]
{\it 
%{\scriptsize $\bullet$} 
(b) ($d=1$) If $\pi(|\xi|^3)<+\infty$, then  $(\xi(X_n))_n$ satisfies a one-dimensional 
uniform Berry-Esseen theorem} (Corollary 6.3). \\[0.1cm]
{\it 
%{\scriptsize $\bullet$} 
(c) ($d=1$) If $\pi(|\xi|^\alpha)<+\infty$ with some $\alpha>3$ and $\xi$ is nonlattice, 
then $(\xi(X_n))_n$ satisfies a  one-dimensional first-order Edgeworth expansion} (Corollary 8.2). \\[0.1cm]
{\it 
%{\scriptsize $\bullet$} 
(d) If $\pi(|\xi|_2^\alpha)<+\infty$ for some $\alpha>\max\left(3,{\lfloor d/2\rfloor+1}\right)$, then 
$(\xi(X_n))_n$ satisfies  a multidimensional Berry-Esseen theorem in the sense of the Prohorov metric} (Corollary 9.2). 

\noindent{\it Application to the Knudsen gas model.} 
Corollary 9.2 just above summarized enables us to specify the slightly incorrect Theorem 2.2.4 of \cite{FPAAP} concerning the Knudsen gas model studied by Boatto and Golse in \cite{BoattoGolse}. Let us briefly recall the link with the uniform ergodicity hypothesis, see \cite{FPAAP} for details. Let $(E,{\cal E},\pi)$ be a probability space, let $T$ be a $\pi$-preserving transformation. The Knudsen gas model can be investigated with the help of the Markov chain $(X_n)_n$ on $(E,{\cal E},\pi)$, whose transition operator $Q$ is defined as follows, for some  $\delta\in(0,1)$: \\[0.1cm]
\indent $\ \ \ \ \ \ \ \ \ \ \ \ \ \ \ \ \ \ \ \ \ \ \ \ \ \ \ \ \ \ \ \ \ \ \ \ Qf=\delta\, \pi(f) + (1-\delta)\, f\circ T.$\\[0.1cm]
Then $(X_n)_n$ is clearly uniformly ergodic. Theorem 2.2.4 of \cite{FPAAP} gave a rate of convergence in $n^{-{1/2}}$ (in the sense of 
the Prohorov metric) in the multidimensional c.l.t.~for $(\xi(X_n))_n$ under the hypothesis $\xi\in \L^3(\pi)\cap \L^{\lfloor d/2\rfloor +1}(\pi)$. However, the proof of this statement is not correct as it is written in \cite{FPAAP} \footnote{Proposition 2.4.2 of \cite{FPAAP} stated that, if $\xi\in \L^3(\pi)\cap \L^{\lfloor d/2\rfloor +1}(\pi)$, then $Q(\cdot)$ defines a regular family of operators when acting on the single space $\cB^{^\infty}$: this result is not true. As already mentioned, it holds under some more restrictive condition of the type $\sup_{x\in E}\int |\xi(y)|^m Q(x,dy) < +\infty$.}. By Corollary 9.2 of the present paper, the  above mentioned rate of convergence is valid if we have 
$\xi\in \L^{3+\varepsilon}(\pi)\cap \L^{\lfloor d/2\rfloor+1+\varepsilon}(\pi)$ for some $\varepsilon>0$.

\indent Of course Example 1 is quite restrictive, and another motivation of this work 
is to present applications to the two next Markov models  of more practical interest. 

\noindent {\bf Example 2: the $v$-geometrically ergodic Markov chains} (see e.g \cite{mey,kontmey}). This example constitutes a natural extension of the previous one. Let $v : E\r[1,+\infty)$ be an unbounded function. Then $(X_n)_n$ is said to be {\it $v$-geometrically 
ergodic} if its transition operator $Q$ satisfies (K1) on the weighted supremum normed space $(\cB_v,\|\cdot\|_v)$ composed of the  measurable complex-valued functions $f$ on $E$ such that $\|f\|_v  = \sup_{x\in E} |f(x)|/v(x) < +\infty$. \\[0.1cm]
Applications of our abstract results to this example are given in Section 10. 
For all our limit theorems (except for the multidimensional Berry Esseen theorem),
when $Y_1\in\L^m$ is needed in the i.i.d.~case, the usual spectral method requires for these models the condition $\sup_{x\in E}v(x)^{-1}\, \int |\xi(y)|^m v(y) Q(x,dy) < +\infty$ (see e.g \cite{fuhlai}) which, in practice, often amounts to assuming that $\xi$ is bounded \cite{kontmey}. Our method only requires that $|\xi|^m\leq C\, v$ or $|\xi|^{m+\varepsilon}\leq C\, v$, which extends the well-known condition $|\xi|^2\leq C\, v$ used for proving the c.l.t.~\cite{mey}.  

\noindent {\bf Example 3: the iterated random Lipschitz models} (see e.g \cite{duf,diaco}).  Except when Harris recurrence and strong mixing hypotheses are assumed, not many works have been devoted to the refinements of the c.l.t.~for the iterative models. As in \cite{mira,hulo,aap}, the important fact here is that these models are Markov chains satisfying (K1) on the weighted Lipschitz-type spaces, first introduced in \cite{lep83}, and slightly modified here according to a  definition due to Guibourg.  
Applications of our results to this example are detailed in Section~11: by considering the general weighted-Lipschitz functionals $\xi$ of \cite{duf},  the limit theorems are stated under some usual moment and mean contraction conditions, which extend those of \cite{duf} \cite{benda} used to prove the c.l.t.. When applied to some classical random iterative models, these assumptions again reduce to the moment conditions of the i.i.d~case (possibly up to $\varepsilon>0$). \\[0.1cm] 
For instance, let us consider in $\R^d$ the affine iterative model $X_n = AX_{n-1}+\theta_n$ where $A$ is a strictly contractive $d\times d$-matrix and $X_0,\theta_1,\theta_2,\ldots$ are $\R^d$-valued independent r.v.. Then, in the case $\xi(x)=x$,
our limit theorems (except the multidimensional
Berry-Esseen theorem) hold if $\theta_1\in\L^m$, 
where $m$ is the corresponding optimal order of the 
i.i.d.~case (up to $\varepsilon>0$ as above for 
the Edgeworth expansion), whereas the usual spectral method requires exponential moment conditions for these statements \cite{mira}.  

\noindent{\bf Extensions.} The operator-type derivation procedure (B) may be also used to investigate renewal theorems \cite{guibourg} \cite{guiher}, and to study the rate of convergence of statistical estimators for strongly ergodic Markov chains (thanks to the control of the constants in (B)), see \cite{jlv}. \\[0.1cm]
Anyway, our method may be employed in other contexts where Fourier operators occur. First, by an easy adaptation of the hypotheses, 
the present limit theorems may be extended to the general setting of Markov random walks (extending the present results to sequence $(X_n,S_n)_n$). Second, these theorems may be stated for the Birkhoff sums stemming from dynamical systems, by adapting the hypotheses to the so-called Perron-Frobenius operator (to pass from Markov chains to dynamical systems, see e.g \cite{hulo} Chap.~XI). \\[0.1cm] 
\noindent The Nagaev-Guivarc'h method can be also used to prove the convergence to stable laws. For this study, the standard perturbation theorem sometimes operates, see \cite{aar2,aar3,bapei,gou,gui-lejan,hh-lh}. But, since the r.v.~which are in the domain of attraction of a stable law are unbounded, the Keller-Liverani theorem is of great interest for these questions. This new 
approach has been introduced in \cite{gouezel} in the context of the stadium billiard, and it has been recently developed in \cite{gui-lepage} for affine random walks and in \cite{seb-08} for Gibbs-Markov maps. \\[0.1cm] 
An important question to get further applications will be to find some others "good" families of spaces to apply the operator-type derivation procedure (B). To that effect, an efficient direction is to use interpolation spaces as in \cite{seb-08}. 

\noindent{\bf Plan of the present paper.} Section 2 presents a well-known central limit theorem based on Gordin's method, with further statements concerning the associated covariance matrix. In Section 3, we summarize the usual spectral method, and we give an 
explicit example (belonging to example 2) to 
which this method cannot be applied. Section 4 
presents the Keller-Liverani perturbation theorem and some first applications concerning the link between the characteristic function of $S_n$ and the eigen-elements of the Fourier kernels $Q(t)$. These preliminary results are then directly applied to prove a  multidimensional local limit theorem in Section 5, and to recall in Section 6 the Berry-Esseen theorem of \cite{ihp2}. Some useful additional results on the non-arithmeticity condition are presented in Section 5.2: these results are detailed in Section 12. Section 7 states the derivation statement mentioned in the above 
procedure (B), and this statement is then applied to prove a first-order Edgeworth expansion (Section 8) and a multidimensional Berry-Esseen type theorem for the Prohorov metric (Section 9). Let us mention that all the operator-type assumptions introduced in the sections 4 and 7, as well as all our limit theorems, will be directly afterward investigated and illustrated through the example 
of the strongly ergodic Markov chains on $\L^2(\pi)$ (Example 1). The applications to Examples 2-3 are deferred to Sections 10-11.  Finally, mention that the proof of the main result of Section 7, and the technical computations involving the weighted Lipschitz-type spaces of Section 11, are relegated to Appendices A-B. 

\noindent{\bf Acknowledgments.} The authors are grateful to the referee for many very helpful comments which allowed to greatly enhance the content and the presentation of this paper. The weighted Lipschitz-type spaces used in Section 11.2 have been introduced by Denis Guibourg in a work (in preparation) concerning the multidimensional Markov renewal theorems. We thank him for accepting that we use in our work this new definition, which allowed us to divide by 2 the order of the moment conditions for the iterative models.

%========================================================
\section{A central limit theorem in the stationary case} 
As a preliminary to the next limit theorems, we state here a well-known 
c.l.t.~for $(\xi(X_n))_n$, which is a standard consequence of a theorem due to Gordin \cite{gor}. We shall then deduce a 
corollary based on Condition (K1). 
In this section, we only consider the stationary case. Let us observe that, concerning distributional questions on $(\xi(X_n))_n$, 
one may without loss of generality assume that $(X_n)_{n\geq0}$ is the canonical Markov chain associated to $Q$. \\[0.1cm]
\noindent So we consider here the usual probability 
space $(E^\N,{\mathcal E}^{\otimes\N},\P_\pi)$ for the 
canonical Markov chain, still denoted by $(X_n)_{n\geq0}$, with transition probability $Q$ and initial stationary distribution $\pi$. 
Let $\theta$ be the shift operator on $E^\N$. As usual we shall say that $(X_n)_{n\geq0}$ is ergodic if the dynamical system 
$(E^\N,{\mathcal E}^{\otimes\N},\P_\pi,\theta)$ is ergodic.  

\noindent{\bf Theorem} {\it (Gordin).  \noindent Assume that $(X_n)_{n\geq0}$ is ergodic, and  \\[0.2cm] 
\indent $\ \ \ \ \forall i=1,\ldots,d,\ \ \ \ \xi_i\in\L^2(\pi)\ $ and 
$\ \breve\xi_i := \sum_{n\geq0} Q^n\xi_i\ \ \mbox{converges in}\ \L^2(\pi).$\\[0.2cm] 
Then $\frac{S_n}{\sqrt n}\, ^{_{\underline{\ \, \tiny{\cD}\, \ }_>}}\, \cN(0,\Gamma)$, where $\Gamma$ is the covariance matrix 
defined by $\langle \Gamma t,t \rangle = \pi(\breve\xi_t^2) - \pi((Q\breve\xi_t)^2)$, where we set 
$\breve\xi_t = \sum_{i=1}^d t_i\, \breve\xi_i$. }

\noindent{\bf Corollary 2.1.} {\it Let us suppose 
that $(X_n)_n$ is ergodic, that (K1) holds on $\cB\hookrightarrow\L^2(\pi)$, and $\xi_i\in\cB$ 
($i=1,\ldots d$). Then the c.l.t.~of the previous theorem holds. }

\noindent{\it Proof of Corollary.} Since we have (K1) on $\cB$, $\xi_i\in\cB$ and $\pi(\xi_i)=0$, 
the series $\breve\xi_i = \sum_{n=0}^{+\infty} Q^n\xi_i$ converges in $\cB$, thus in $\L^2(\pi)$. \fdem

\noindent For instance, if $(X_n)_n$ is strongly ergodic on $\L^2(\pi)$ (see Example 1), then $(X_n)_n$ is ergodic \cite{rosen}, and we find again the well-known fact that the central limit theorem holds in the stationary case when  $\pi(|\xi|_2^2)<+\infty$. In order to make easier the use of Corollary~2.1 in other models, let us recall the following sufficient condition for $(X_n)_n$ to be ergodic. This statement, again in relation with Condition (K1), is 
established in \cite{hulo} (Th.~IX.2) with the help of standard arguments based on the monotone class theorem. 

\noindent{\bf Proposition 2.2.} {\it  Let us suppose 
that (K1) holds with $\cB$ satisfying the additional following conditions: 
$\cB$ generates the $\sigma$-algebra $\cE$, $\delta_x\in\cB'$ for all $x\in E$, and $\cB\cap\cB^{^\infty}$is stable under 
product. Then $(X_n)_n$ is ergodic. }

\noindent Of course, other methods exist to investigate the c.l.t.~for Markov chains, but Corollary 2.1 is sufficient 
for our purposes: indeed, it is easily applicable to our examples, and it enables us to define the asymptotic covariance matrix $\Gamma$ which will occur in all the others limit theorems.

\noindent The above definition of $\Gamma$ provides the following classical characterisation of the case when $\Gamma$ is degenerate. 

\noindent{\bf Proposition 2.3.} {\it Under the hypotheses of Corollary 2.1, $\Gamma$ is 
non invertible if and only if }
$$\exists t\in\R^d,\ t\neq 0,\ \exists g\in \cB,\ \ \ \langle t,\xi(X_1)\rangle = g(X_0) - g(X_1)\ \ \P_\pi-a.s..$$
\noindent Let us notice that this equivalence is still true 
for $\cB=\L^2(\pi)$ if we know that: \\[0.15cm] 
\indent $\ \ \ \ \ \ \ \ \ \ \ \ \ \ \ \ \ \ \ \ \ \ \ \ \ 
\forall t\in{\R}^d,\ \ \sup_{n\ge 1}\left\vert n\langle \Gamma t,t\rangle
        -{\E}_\pi[\langle t,S_n\rangle^2]\right\vert<+\infty$. 
 
\noindent {\it Proof of Proposition 2.3.} If $\langle t,\xi(X_1)\rangle 
= g(X_0) - g(X_1)\ \ \P_\pi$-a.s., then 
$\left(\frac{\langle t,S_n\rangle}{\sqrt n}\right)_n$ 
converges in distribution to the Dirac mass at 0 (which 
proves that $\Gamma$ is non invertible). Indeed, 
by stationarity, we have $\langle t,\xi(X_n)\rangle 
= g(X_{n-1}) - g(X_n)\ \ \P_\pi$-a.s.~for 
all $n\geq1$, so $\langle t,S_n\rangle = g(X_0) - g(X_n)$. Since we have $g\in\cB\hookrightarrow\L^2(\pi)$, 
this implies that $\lim_n\E_\pi[(\frac{\langle 
t,S_n\rangle}{\sqrt n})^2] = 0$ and hence the desired statement. 
Conversely, let us suppose that $\Gamma$ is not invertible. Then there exists $t\in\R^d$, $t\neq 0$, such that $\langle \Gamma t,t \rangle = 0$.
From the definition of $\Gamma$ given in the above theorem and from the obvious 
equality $\E_\pi[(\breve\xi_t(X_1) - Q\breve\xi_t(X_0))^2] = 
\pi(\breve\xi_t^2) - \pi((Q\breve\xi_t)^2)$, it follows that   \\[0.12cm] 
\indent $\ \ \ \ \ \ \ \ \ \ \ \ \ \ \ \ \ \ \ \ \ \ \ \ \ \ \ \ \ \ \ \ \ \ 
\E_\pi[(\breve\xi_t(X_1) - Q\breve\xi_t(X_0))^2] = 0$.  \\[0.12cm]
Thus $\breve\xi_t(X_1) - Q\breve\xi_t(X_0) = 0\ \ \P_\pi-a.s.$. 
Set $\xi_t(\cdot) =  \langle t,\xi(\cdot)\rangle$. By definition of $\breve\xi_t$, we have 
$\breve\xi_t = \xi_t + Q\breve\xi_t$, so 
$$\xi_t(X_1) + Q\breve\xi_t(X_1)-Q\breve\xi_t(X_0) = 0 \ \ \P_\pi-a.s..$$
This yields 
$\xi_t(X_1) = g(X_0) - g(X_1)\ \ \P_\pi-a.s.$ with $g=Q\breve\xi_t$. \fdem

\noindent The previous proposition can be specified as follows. 

\noindent{\bf Proposition 2.4.} {\it
Let $t\in\R^d$, $t\neq 0$, and let $g$ be a measurable function on $E$ such that: \\[0.15cm] 
\indent $\ \ \ \ \ \ \ \ \ \ \ \ \ \ \ \ \ \ \ \ \ \ \ \ \ \ \ \ 
\langle t,\xi(X_1)\rangle = g(X_0) - g(X_1)\ \ 
\P_\pi-a.s.$. \\[0.15cm] 
Then there exists a $\pi$-full $Q$-absorbing set $A\in\cE$ such that we have: \\[0.15cm] 
\indent $\ \ \ \ \ \ \ \ \ \ \ \ \ \ \ \ \ \ \ \ \ \ \ \ \ 
\forall x\in A,\ \langle t,\xi(y)\rangle = g(x) - g(y)\ \, Q(x,dy)$-a.s..} 

\noindent {\it Proof.} 
For $x\in E$, set $B_x = \{y\in E\, :\, \langle t,
\xi(y)\rangle = g(x) - g(y)\}$. By hypothesis we have 
$\int Q(x,B_x)d\pi(x) = 1$, and since 
$Q(x,B_x)\leq 1$, this gives $Q(x,B_x) = 1$ $\pi$-a.s.. Thus 
there exists a $\pi$-full set $A_0\in\cE$ such that $Q(x,B_x) = 1$ for $x\in A_0$. 
From $\pi(A_0) = 1$ and the invariance of $\pi$, we also have $\pi(Q1_{A_0}) = 1$, and since 
$Q1_{A_0}\leq Q1_E = 1_E$, this implies that 
$Q(\cdot,A_0) = 1\ \pi$-a.s.. Again there exists 
a $\pi$-full set $A_1\in\cE$ such that $Q(x,A_0) = 1$ for $x\in A_1$. Repeating this procedure, 
one then obtains a family $\{A_n, n\geq1\}$ of 
$\pi$-full sets satisfying by construction the condition: $\forall n\geq1,\, 
\forall x\in A_n,\ Q(x,A_{n-1})=1$. 
Now the set $A := \cap_{n\geq0} A_n$ is $\pi$-full and, for any 
$a\in A$, we have $Q(a,A_{n-1})=1$ for all $n\geq1$, thus $Q(a,A)=1$. This 
proves that $A$ is $Q$-absorbing, and 
the desired equality follows from the inclusion $A\subset A_0$. \fdem 
%==========================================================================
\section{The usual Nagaev-Guivarc'h method} 
\noindent The characteristic function of $S_n$ is linked to the Fourier kernels 
$Q(t)(x,dy) = e^{i\langle t,\, \xi(y)\rangle}Q(x,dy)$ of ($Q,\xi$) by the following formula (see e.g \cite{hulo} p.~23)\\[0.23cm]
\indent  {\bf (CF)} $\displaystyle\ \ \ \ \ \ \ \ \ \ \ \ \ \ \ \ \ \ 
\forall n\geq 1,\ \forall t\in\R^d,\ \E_\mu[e^{i\langle t,S_n\rangle}] = \mu(Q(t)^n1_E)$, \\[0.23cm]
and the Nagaev-Guivarc'h method consists in applying to $Q(t)$ the standard perturbation theory~\cite{ds}. For this to make sense, 
one must assume that $Q$ satisfies Condition (K1) (of Section~1) on $\cB$, that $Q(t)\in\cL(\cB)$, and 
that $Q(\cdot)$ is $m$ times continuously differentiable from $\R^d$ to $\cL(\cB)\, $ ($m\in\N^*$). In this case,  $Q(t)^n$, hence $\E_\mu[e^{itS_n}]$, can be expressed in function of $\lambda(t)^n$, where $\lambda(t)$, the dominating eigenvalue of $Q(t)$, is also  $m$ times continuously differentiable. Then, the classical limit theorems (based on Fourier techniques), requiring $Y_1\in\L^m$ for a i.i.d.~sequence $(Y_n)_n$, extend to $(\xi(X_n))_n$, see for example \cite{lep82,rou,guihar,broi,hulo}. Unfortunately, the previous regularity assumption on $Q(\cdot)$ (in case $d=1$ for simplicity) requires that the kernel $\xi(y)^mQ(x,dy)$ continuously acts on $\cB$: this is what we called an operator-moment condition in Section~1, and we already mentioned that, if $\xi$ is unbounded, this assumption is in general very restrictive. \\[0.1cm] 
\indent Actually Nagaev established in \cite{nag1} a c.l.t., and a local limit theorem in the countable case, for the uniformly ergodic Markov chains (see Ex.~1 of Section~1), and he did not appeal to operator-moment conditions: indeed, Nagaev first applied the standard perturbation theorem for bounded functionals $\xi$, and by using some intricate truncation techniques, he extended his results under the condition $\pi(|\xi|^2)<+\infty$. However afterward, this truncation method has not been used any more. In particular, the Berry-Esseen theorem in \cite{nag2} was stated under the operator-moment assumption $\sup_{x\in E} \int_E|\xi(y)|^3\, Q(x,dy) < +\infty$, which is clearly necessary and sufficient for $Q(\cdot)$ to be three times continuously differentiable from 
$\R$ to $\cL(\cB^{^\infty})$.  \\[0.1cm]
\indent The use of the standard perturbation theory is even more difficult in Examples 2-3 of Section~1: the typical example below shows that, neither the operator-moment conditions, nor even the simple assumption $\|Q(t)-Q\|_{\scriptsize \cB}\r0$, hold in general when $\xi$ is unbounded.  

\noindent{\it Counter-example.} Let $(X_n)_{n\geq0}$ be the real-valued autoregressive chain defined by \\[0.1cm]
\indent $\displaystyle\ \ \ \ \ \ \ \ \ \ \ \ \ \ \ \ \ \ \ \ \ \ \ \ \ \ \ \ \ \ \ \ \ \  
X_n = a X_{n-1} +  \theta_n\ \ (n\in\N^*),$ \\[0.1cm]
where $a\in(-1,1)$, $a\neq 0$, $X_0$ is a real r.v.
and $( \theta_n)_{n\geq1}$ is a sequence of i.i.d.r.v., independent of $X_0$. 
Assume that $ \theta_1$ has a positive density $p$ with finite variance. 
It is well-known that $(X_n)_{n\geq0}$ is a Markov chain whose 
transition probability is: $(Qf)(x) = \int_\R f(ax+y)\, p(y)\, dy$.  \\[0.1cm]
\noindent Set $v(x)=1+x^2$ ($x\in\R$). Using the so-called drift condition (see \cite{mey}, Section 15.5.2), one can prove that 
$(X_n)_{n\geq0}$ is $v$-geometrically ergodic (see Example 2 in Section 1). 
Now let us consider the functional $\xi(x) = x$. 
We have for any $x\in\R$ \\
\indent $\displaystyle\ \ \ \ \ \ \ \ \ \ \ \ \ \ \ \ \ \ \ \ \ \ \ \  \ \ \ \ \ \ \ \ \ 
Q(\xi^2\, v)(x) \geq \int_\R (ax+y)^4\, p(y)\, dy$. \\
If $\int_\R y^4\, p(y)\, dy = +\infty$, then $Q(\xi^2\, v)$ is not defined. If $\int_\R y^4\, p(y)\, dy < +\infty$, then 
$Q(\xi^2\, v)$ is a polynomial function of degree 4, so that \\
\indent $\displaystyle\ \ \ \ \ \ \ \ \ \ \ \ \ \ \ \ \ \ \ \ \ \ \ \  \ \ \ \ \ \ \ \ \ \ \ \ \ \ \ 
\sup_{x\in E} \frac{|Q(\xi^2\, v)(x)|}{1+x^2} = +\infty$, \\
that is, $Q(\xi^2\, v)\notin\cB_v$. Similarly we have $Q(|\xi|\, v)\notin\cB_v$.  Thus neither $\xi(y)Q(x,dy)$, nor $\xi(y)^2Q(x,dy)$, 
continuously act on $\cB_v$. Actually, even the continuity condition $\|Q(t)-Q\|_{\scriptsize \cB_v}\r0$ is not valid. 
To see that, it suffices to establish that, if $g(x) = x^2$, 
then $\|Q(t)g-Qg\|_v = \sup_{x\in\R}(1+x^2)^{-1}\, |Q(t)g(x) - Qg(x)|$ does not converge to 0 when $t\r0$. 
Set $p_1(y) = yp(y)$ and $p_2(y) = y^2p(y)$, and denote by 
$\hat\phi(t) = \int_\R\phi(y)e^{ity}dy$ the Fourier transform of any integrable function $\phi$ on $\R$. Then \\[0.1cm]
\indent $\displaystyle\ \ \ \ \ \ \ \ \ 
Q(t)g(x) = \int_\R e^{it(ax+y)}\, (y+ax)^2\, p(y)\, dy = e^{iatx}\, [ \hat p_2(t) + 2ax\hat p_1(t) + a^2x^2\, \hat p(t)]$. \\[0.1cm]
So $\displaystyle Q(t)g(x) - Qg(x) =   \bigg(e^{iatx}\hat p_2(t) - \hat p_2(0) +
2ax\, [e^{iatx}\hat p_1(t) - \hat p_1(0)]\bigg) +  a^2x^2\, [e^{iatx}\hat p(t) - 1]$. Using the inequality 
$|e^{iu} - 1| \leq |u|$, we easily see that there exists a constant $C>0$ such that 
$$\sup_{x\in\R} (1+x^2)^{-1}\, 
\bigg| e^{iatx}\hat p_2(t) - \hat p_2(0) + 2ax\, [e^{iatx}\hat p_1(t) - \hat p_1(0)]\bigg| \leq 
C\, \bigg(|t| +  |\hat p_2(t)-\hat p_2(0)| + |\hat p_1(t)-\hat p_1(0)|\bigg).$$
By continuity of $\hat p_1$ and $\hat p_2$, the last term converges to 0 as $t\r0$. Now set \\[0.1cm]
\indent $\displaystyle\ \ \ \ \ \ \ \ \ \ \ \ \ \ \ \ \ \ \ \ \ \ \ \ \ \ \ 
\psi(x,t) = (1+x^2)^{-1}\,  a^2x^2\, |e^{iatx}\hat p(t) - 1|$. \\[0.1cm]
We have 
$\sup_{x\in\R} \psi(x,t) \geq \psi(\frac{\pi}{at},t) = \frac{a^2\pi^2}{\pi^2+a^2t^2}\, |\hat p(t) + 1|$. 
Since this last term converges to $2a^2\neq0$ as $t\r0$, this clearly implies the desired statement. 
%=================================================================================
\section{The Nagaev-Guivarc'h method via the Keller-Liverani theorem} 
The next statement is the perturbation theorem of Keller-Liverani, when applied to the Fourier Kernels $Q(t)$ under Condition (K1) of Section 1. The present assumptions will be discussed, and illustrated in the case of the strongly ergodic Markov chains on $\L^2(\pi)$.  Finally we shall present a first probabilistic application to the characteristic function of $S_n$. 

\noindent {\bf The perturbation operator theorem of Keller-Liverani.} 

\noindent {\it Condition ($\widetilde{K}$):  $Q$ satisfies Condition (K1) (of Section 1) on $\cB$, and 
there exists a neighbourhood ${\cal O}$ of $0$ in $\R^d$ and a Banach space $\widetilde{\cB}$ satisfying  $\cB\hookrightarrow\widetilde{\cB}\hookrightarrow\L^1(\pi)$, such that we have $Q(t)\in\cL(\cB)\cap\cL(\widetilde{\cB})$ for each 
$t\in {\cal O}$, and: 

\noindent {($\widetilde{K2}$)} $\displaystyle 
\forall t\in{\cal O},\ \lim_{h\r0}\|Q(t+h)-Q(t)\|_
{\scriptsize{\cB,\widetilde{\cB}}}=0$  

\noindent {($\widetilde{K3})$} $\ \exists \kappa_1<1,\ \exists C> 0,\ \forall n\geq1,\ \forall f\in\cB,\ \forall t\in {\cal O}, \ \
\|Q(t)^nf\|_{\scriptsize \cB} \leq C\, \kappa_1^n\, \|f\|_{\scriptsize \cB} + C\, \|f\|_{\scriptsize{\widetilde{\cB}}}$.  

}
\noindent {\it Condition (${K}$):  
Condition ($\widetilde{K}$) with 
$\widetilde{\cB} = \L^1(\pi)$.}

\noindent Under Condition ($\widetilde{K}$), we 
denote by $\kappa$ any real number such that 
$\max\{\kappa_0,\kappa_1\} < \kappa < 1$, where 
$\kappa_0$ is given in Condition~(K'1) of Section~1, and we define the following set 
$$\cD_\kappa = \bigg\{z : z\in\C,|z|\geq \kappa,\ |z-1|\geq \frac{1-\kappa}{2}\bigg\}.$$
\noindent{\bf Theorem} (K-L) {\it \cite{keli,liverani} (see also \cite{bal}). 
Let us assume that Condition ($\widetilde{K}$) holds. Then, for all $t\in {\cal O}$ (with possibly ${\cal O}$ reduced), $Q(t)$ admits a dominating eigenvalue $\lambda(t)\in\C$, 
with a corresponding rank-one eigenprojection 
$\Pi(t)$ satisfying
$\Pi(t)Q(t)=Q(t)\Pi(t)=\lambda(t)\Pi(t)$, such that we have the following properties: 
$$\lim_{t\r0}\lambda(t) = 1, \ \ \ \ \sup_{t\in {\cal O}}\|Q(t)^n -  \lambda(t)^n\, \Pi(t)\|_{\scriptsize \cB} = O(\kappa^n),\ \ \ \  \lim_{t\r0}\|\Pi(t) - \Pi\|_{\scriptsize{\cB,\widetilde{\cB}}}=0,$$
and finally $\, \cM := \sup\big\{\|(z-Q(t))^{-1}\|_{\scriptsize \cB},\ t\in {\cal O},\ z\in\cD_\kappa \big\}< +\infty$. }

\noindent Let us moreover mention that $\lambda(t)$ and $\Pi(t)$ can be expressed in terms of $(z-Q(t))^{-1}$ (see the proof of
Corollary 7.2 where the explicit formulas are given and used). 

\noindent {\bf Remark.} {\it The conclusions of Theorem (K-L) still hold when Condition {($\widetilde{K2}$)} is replaced with: $\lim_{h\r0}\|Q(h)-Q\|_{\scriptsize{\cB,\widetilde{\cB}}}=0$. In fact, Condition {($\widetilde{K2}$)} provides the following additional property, that will be used in Section~5.1: $\lambda(\cdot)$ is continuous on ${\cal O}$ (see \cite{ihp1}). Anyway, in most of cases, the previous continuity condition at $t=0$ implies {($\widetilde{K2}$)} (see for instance Rk.~(a) below). Let us also recall that the neighbourhood ${\cal O}$ and the bound $\cM$ of Theorem (K-L) depend on $\kappa$ (with $\kappa$ fixed as above) and on the following quantities (see \cite{keli} p.~145): \\[0.12cm]
\indent - the constant $H := \sup\{\|(z-Q)^{-1}\|_{\scriptsize \cB},\ z\in\cD_\kappa\}$, which is finite by (K1), \\[0.1cm]
\indent - the rate of convergence of $\|Q(t)-Q\|_{\scriptsize{\cB,\widetilde{\cB}}}$ to $0$ when $t\r0$, \\
\indent - the operator norms $\|Q\|_{\scriptsize \cB}$, $\|Q\|_{\scriptsize \widetilde\cB}$, and the constants $C$, $\kappa_1$ of Condition {($\widetilde{K3})$}. \\[0.12cm]
This remark is relevant since the asymptotic properties of Theorem (K-L) depend on $\cM$.}   

\noindent {\bf Some comments on  Condition ($\widetilde{K}$).} 

The hypotheses in \cite{keli} are stated with the help of an auxiliary norm on $\cB$ (which can be easily replaced by a semi-norm). In practice, this auxiliary norm is the restriction of the norm of a usual  Banach space $\widetilde{\cB}\hookrightarrow\L^1(\pi)$. It is the reason why Condition ($\widetilde{K}$) has been presented with an auxiliary space. The dominated hypothesis between the norms stated in \cite{keli} is here replaced with our assumption $\cB\hookrightarrow\widetilde{\cB}$. The fact that $\widetilde{\cB}$ is complete and $\widetilde{\cB}\hookrightarrow\L^1(\pi)$ is not necessary for the validity of Theorem (K-L), but these two hypotheses are satisfied in practice. Moreover the assumption $\widetilde{\cB}\hookrightarrow\L^1(\pi)$ ensures that $\pi\in\widetilde{\cB}'$, which is important 
for our next probabilistic applications. 
Let us also mention that \cite{keli} appeals to the following additional condition on the essential spectral radius of $Q(t)$: $\forall t\in {\cal O},\ r_{ess}(Q(t)) \leq \kappa_1$. As explained in \cite{liverani}, this assumption is not necessary for Theorem~(K-L), thanks to Condition~(K1). It will be assumed in Section  5.1 for applying \cite{keli} to $Q(t)$ for $t$ close to $t_0\neq0$. \\[0.1cm]
It is worth noticing that the continuity Condition ($\widetilde{K2}$) is less 
restrictive than the condition $\|Q(t+h)-Q(t)\|_{\scriptsize \cB}\r0$ required in the usual perturbation theorem. In fact, despite their not very probabilistic appearance, the conditions ($\widetilde{K2}$) ($\widetilde{K3}$) are suited to many examples of strongly ergodic Markov chains: for instance, they hold for any measurable functional $\xi$ in the case of the strongly ergodic Markov chains on $\L^2(\pi)$ 
and of the $v$-geometrically ergodic Markov chains (see Prop.~4.1, Lem.~10.1), and they are valid under simple mean contraction and moment conditions for  iterative Lipschitz models 
(see section 11).  

\noindent {\bf Some comments on Condition (${K}$).} 

In the special case $\widetilde{\cB} = \L^1(\pi)$, we shall use repeatedly the next simple remarks. 

\noindent (a) First observe that we have $\sup_{t\in\R^d}\|Q(t)\|_{\scriptsize{\L^1(\pi)}} < +\infty$ (use $|Q(t)^nf|\leq Q^n|f|$ and the $Q$-invariance of $\pi$). Besides the following condition 
$$\sup\big\{\pi(|e^{i\langle t,\, \xi\rangle}-1|\, |f|),\ f\in\cB,\ \|f\|_{\scriptsize \cB}\leq 1\big\}\ \ \mbox{converges to}\ 0\ \mbox{when}\ t\r0,$$  
which is for instance satisfied if $\cB\hookrightarrow\L^p(\pi)$ for some $p>1$ (by H\"older's inequality and Lebesgue's theorem) is a sufficient condition for the continuity assumption of Condition (K). More precisely, the above property implies that 
$$\forall t\in\R^d,\ \ \lim_{h\r0}\|Q(t+h)-Q(t)\|_{\scriptsize{\cB,\L^1(\pi)}}=0.$$
Indeed we have for any $f\in\cB$
$$\pi\big(|Q(t+h)f-Q(t)f|\big) \leq  \pi\big(Q|e^{i\langle h,\, \xi\rangle}-1|\, |f|\big) = \pi\big(|e^{i\langle h,\, \xi\rangle}-1|\, |f|\big).$$
\noindent  (b) Recall that $\cB$ is a Banach lattice if we have: $|f|\leq |g|\ \Rightarrow\ \|f\|_{\scriptsize \cB}\leq \|g\|_{\scriptsize \cB}$ for any $f,g\in\cB$. The examples of Banach lattices in our work are: $\cB=\cB^{^\infty}$, $\cB=\L^p(\pi)$ (used in Ex.~1 of Section~1), and $\cB=\cB_v$ (used in Ex.~2). Another classical example is the space of the bounded continuous functions on $E$. \\
Let us assume that $\cB$ is a Banach lattice such that: $\forall t\in\R^d,\ \forall f\in\cB,\ \ e^{i\langle t,\, \xi\rangle}\cdot f\in\cB$. Then Condition~(K1) implies ($\widetilde{K3})$ with $\widetilde{\cB} = \L^1(\pi)$ and ${\cal O} = \R^d$. \\[0.1cm]
Indeed, 
we have $|Q(t)^nf| \leq Q^n|f|$,  so $\|Q(t)^nf\|_{\scriptsize \cB} \leq \|\, Q^n|f|\, \|_{\scriptsize \cB}$, and (K1) then gives for 
all $n\geq1$, $f\in\cB$ and $t\in \R^d$: $\|Q(t)^nf\|_{\scriptsize \cB} \leq C\, \kappa_0^n\, \|f\|_{\scriptsize \cB} + \pi(|f|)\, \|1_E\|_{\scriptsize \cB}$. 

\noindent (c) If ($\widetilde{K3})$ is fulfilled with $\widetilde{\cB} = \L^1(\pi)$, then it holds for any $\widetilde{\cB}\hookrightarrow\L^1(\pi)$.  

\noindent {\bf Example} (the strongly ergodic Markov chains on $\L^2(\pi)$, see Ex.~1 of Section~1): 

\noindent {\bf Proposition 4.1.} {\it Assume that $(X_n)_{n\geq0}$ is a strongly ergodic Markov chain on $\L^2(\pi)$, that $\xi$ is any $\R^d$-valued measurable function, and let $1\leq p' <  p<+\infty$. Then 
we have ($\widetilde{K}$) with ${\cal O} = \R^d$, $\cB = \L^p(\pi)$, and $\widetilde{\cB} = \L^{p'}(\pi)$. }

\noindent {\it Proof.} We know that $Q$ satisfies Condition (K1) of Section 1 on $\L^p(\pi)$ (see \cite{rosen}). From the above remarks  (b) (c), we then have ($\widetilde{K3}$) with ${\cal O} = \R^d$, $\cB = \L^p(\pi)$, and $\widetilde{\cB} = \L^{p'}(\pi)$. Condition ($\widetilde{K2}$) follows from the next lemma. \fdem 

\noindent {\bf Lemma 4.2.} {\it Let $1\leq p' <  p$, and $t\in \R^d$. Then 
$\lim_{h\r0}\|Q(t+h)-Q(t)\|_{\scriptsize{\L^p(\pi),\L^{p'}(\pi)}}=0$. }

\noindent{\it Proof.} Let us denote $\|\cdot\|_p$ for $\|\cdot\|_{\L^p(\pi)}$.  Using the inequality $|e^{ia}-1|\leq 2\min\{1,|a|\}\ $ 
($a\in \R$) and the H\"older inequality, 
one gets for $t,h\in \R^d$ and $f\in\L^p(\pi)$,  
\begin{eqnarray*}
\left\Vert Q(t+h)f-Q(t)f \right\Vert_{ p'} &\leq&
\left\Vert Q(|e^{i\langle h,\xi\rangle}-1|\, |f|)\right\Vert_{ p'} \\
&\leq& 2\left\Vert \min\{1,|\langle h,\xi\rangle|\}|f|\right\Vert_{ p'}\\
&\leq & 2\left\Vert \min\{1,|\langle h,\xi\rangle|\}\right\Vert_{ p p'\over p- p'}
     \Vert f\Vert_ p,
\end{eqnarray*}
with $\left \Vert \min\{1,|\langle h,\xi\rangle|\}\right\Vert_{ p p'\over p- p'}
\r0$ when $h\r 0$ by Lebesgue's theorem. \fdem 

\indent To end this section, let us return to our general setting and present a first probabilistic application of Theorem~(K-L). 

\noindent {\bf Link between $\lambda(t)$ and the characteristic function of $S_n$.}  

For convenience, let us repeat the basic formula (CF), already formulated in Section 3, which links the characteristic function of $S_n$ with the Fourier kernels of ($Q,\xi$): 
$\forall n\geq 1,\ \forall t\in\R^d,\ \ \E_\mu[e^{i\langle t,S_n\rangle}] = \mu(Q(t)^n1_E)$,
where $\mu$ is the initial distribution of the chain. We appeal here to Theorem (K-L), in particular to the dominating eigenvalue $\lambda(t)$ of $Q(t)$, $t\in{\cal O}$, to the associated rank-one eigenprojection $\Pi(t)$, and finally to the real number $\kappa$ for which we just 
recall that $\kappa<1$.  

\noindent{\bf Lemma 4.3.} {\it Assume ($\widetilde{K}$) and $\mu\in\widetilde{\cB}'$, and set $\ell(t) = \mu(\Pi(t)1_E)$. Then we have: }
$$\lim_{t\r0} \ell(t) = 1\ \ \ \ \ \mbox{and}\ \ \ \ \sup_{t\in {\cal O}}\, \big|\E_\mu[e^{i\langle t,S_n\rangle}] -  \lambda(t)^n\, \ell(t)\big| =  O(\kappa^{n}).$$
\noindent {\it Proof.} Lemma 4.3 directly follows from Theorem (K-L) and Formula (CF). \fdem 
%===============================================================================
\section{A multidimensional local limit theorem} 
The previous lemma constitutes the necessary preliminary to employ Fourier techniques.  
However, it is worth noticing that, except $\lim_{t\r0} \lambda(t) = 1$, 
the perturbation theorem of Keller-Liverani cannot yield anyway the Taylor expansions needed 
for $\lambda(t)$ in Fourier techniques. 
An abstract operator-type hypothesis will be presented in Section 7 
in order to ensure the existence of $m$ continuous derivatives for $\lambda(\cdot)$. 
But we want before to recall another method, based on weaker and more simple probabilistic 
c.l.t.-type assumptions, 
which provides second or third-order Taylor expansions of $\lambda(t)$ near $t=0$. 
As in the i.i.d.~case, these expansions are sufficient to establish a multidimensional local limit theorem, 
this is the goal of the present section, and a one-dimensional uniform Berry-Esseen theorem which will be presented in Section 6. 
\footnote{These two limit theorems could also be deduced from respectively Conditions $\cC(2)$ and  $\cC(3)$ 
of Section~7, but in practice, these two conditions are slightly  more 
restrictive than those of Sections 5-6. For instance, compare $\cC(2)$ and $\cC(3)$ for the strongly ergodic Markov chains on $\L^2(\pi)$ (see Prop.~7.3) with the conditions of Coro.~5.5 and 6.3.} \\[0.1cm] 
Theorem 5.1 below has been established for real-valued functionals in \cite{ihp1} under slightly different hypotheses. Here we present an easy extension to the multidimensional case. Section~5.2 states some expected statements on the Markov non-arithmetic and nonlattice conditions. The application to the strongly ergodic Markov chains on $\L^2(\pi)$ in Section 5.3 is new.    
%===========================
\subsection{A general statement} 
To state the local limit theorem, one needs to introduce the two following conditions. 
The first one is the central limit assumption stated under $\P_\pi$
for which one may appeal to Corollary~2.1 for instance.  
The second one is a spectral non-arithmeticity condition. Recall that, by hypothesis, we have $\pi(\xi)=0$, so that $\E_\pi[S_n] = 0$. 

\noindent {\it Condition (CLT): Under $\P_\pi$, $\ \frac{S_n}{\sqrt n}\, ^{_{\underline{\ \, \tiny{\cD}\, \ }_>}}\cN(0,\Gamma)$, 
with a non-singular matrix $\Gamma$. }

\noindent{\it Condition (S): For all $t\in\R^d$, $Q(t)\in\cL(\cB)$, and for each compact set $K_0$ in $\R^d\setminus\{0\}$, 
there exist $\rho<1$ and $c\geq 0$ such that we have, for all $n\geq1$ and $t\in K_0$, $\|Q(t)^n\|_{_{\scriptsize \cB}} \leq c\, \rho^n$.}

\noindent Condition (S) constitutes the tailor-made hypothesis to operate in the spectral method the proofs of the 
i.i.d.~limit theorems involving the so-called nonlattice assumption. Condition (S) will be reduced to more practical hypotheses in Section 5.2. 

\noindent We want to prove that, given some fixed 
positive function $f$ on $E$ and some fixed real-valued measurable function $h$ on $E$, we have 
$${(LLT)} \ \ \lim_n 
  \sup_{a\in\R^d}\left\vert
   \sqrt{\det\Gamma}\, (2\pi n)^{\frac{d}{2}}\, 
\E_\mu[\, f(X_n)\, g(S_n-a)\, h(X_0)\, ] - 
   e^{-\frac{1}{2n}\langle 
\Gamma^{-1}a,a\rangle}\,  \mu(h)\, \pi(f)\, \int_{\R^d} g(x)dx
   \right\vert = 0,$$
{\it for all compactly supported continuous 
function $g : \R^d\r\R$}. \\
The conditions on $f$, $h$ and $\mu$ are specified below. 

\noindent {\bf Theorem 5.1.} {\it Assume that Condition (CLT) holds, that Condition ($\widetilde{K}$) (of Section 4) holds 
w.r.t.~some spaces $\cB$, $\widetilde{\cB}$, and 
that Condition (S) holds on $\cB$. Finally assume 
$(h\mu)\in\widetilde{\cB}'$ and $f\in\cB$, $f\geq 0$. Then we have (LLT). }

Before going into the proof, let us notice that this result
can be easily extended to any real-valued function $f\in\cB$ such that
$\max(f,0)$ and $\min(f,0)$ belong to $\cB$.\medskip

\noindent {\it Proof of Theorem 5.1.}  In order to use Lemma 4.3 and to write out the Fourier techniques of the i.i.d.~case \cite{bre}, 
one needs to establish a second-order Taylor expansion for $\lambda(t)$. 

\noindent{\bf Lemma 5.2.}  {\it Under the conditions ($\widetilde{K}$) and (CLT), we have for $u\in\R^d$ 
close to $0$:  }
$$\lambda(u) = 1 - \frac{1}{2}\langle \Gamma u,u\rangle + o(\|u\|^2).$$
\noindent {\it Proof (sketch)}. For $d=1$, the proof of Lemma 5.2 is presented in \cite{ihp1}, let us just recall the main ideas. 
By hypothesis, we have $\frac{S_n}{\sqrt n}\, ^{_{\underline{\ \, \tiny{\cD}\, \ }_>}}\cN(0,\sigma^2)$ under $\P_\pi$, with $\sigma^2>0$. Besides, $\widetilde{\cB}\hookrightarrow\L^1(\pi)$ implies $\pi\in\widetilde{\cB}'$. So, from L\'evy's theorem and Lemma 4.3 (applied here with $\mu=\pi$), it follows that $\lim_n\lambda(\frac{t}{\sqrt n})^n = e^{-\frac{\sigma^2}{2}t^2}$, with uniform convergence on any compact set in $\R$. 
Then the fact that $\log \lambda(\frac{t}{\sqrt n})^n = n\log\lambda(\frac{t}{\sqrt n})$ and 
$\log \lambda(\frac{t}{\sqrt n}) \sim  \lambda(\frac{t}{\sqrt n}) -1$ 
gives for $t\neq0$~: \\[0.1cm]
\indent $\displaystyle \ \ \ \ \ \ \ \ \ \ \ \ \ \ \ \ \ \ \ \ \ \ \ \ 
$$(\frac{\sqrt n}{t})^2 \big(\lambda(\frac{t}{\sqrt n}) - 1\big) + \frac{\sigma^2}{2} = o(1)\ \ \mbox{when}\ n\r+\infty.$ \\[0.1cm]  
Setting $u=\frac{t}{\sqrt n}$, it is then not hard 
to deduce the stated Taylor expansion
(see \cite{ihp1} Lem.~4.2). 
  
\noindent These arguments can be readily repeated for $d\geq2$. (To get 
$\log \lambda(\frac{t}{\sqrt n})^n = n\log\lambda(\frac{t}{\sqrt n})$ in $d\geq2$, proceed as in \cite{ihp1} with 
$\psi(x) = \lambda(x \frac{t}{\sqrt n})$, $x\in[0,1]$; the continuity of $\lambda(\cdot)$ on some neighbourhood of 0, helpful for this part \footnote{This continuity property is proved in \cite{ihp1} by applying \cite{keli} to the family $\{Q(t),\ t\in{\cal O}\}$ when $t$ goes to any fixed $t_0\in\cal O$. To that effect, notice that, according to theorem (K-L), we have $r_{ess}(Q(t))\le\kappa$
for all $t\in{\cal O}$.}, obviously extends to $d\geq2$). \fdem 

\noindent If $f=h=1_E$, then (LLT) follows from Lemma 4.3, by writing out the i.i.d.~Fourier techniques of \cite{bre}. In particular, Condition (S) plays the same 
role as the nonlattice condition of 
\cite{bre}. If $f\in\cB$, $f\geq 0$, and $h : E\r\R$ is measurable, one can proceed in the same way by using 
the following equality, of which (CF) is a special case (see e.g \cite{hulo} p.~23), 
\\[0.15cm]
\indent  {\bf (CF')} $\displaystyle\ \ \ \ \ \ \ \ \ \ \ 
\forall n\geq 1,\ \forall t\in\R^d,\ \ \ 
\E_\mu[f(X_n)\, e^{itS_n}\, h(X_0)] = (h\mu)(Q(t)^nf),$\\[0.15cm]
and by using an obvious extension of Lemma 4.3. \fdem 
%=====================================
\subsection{Study of Condition (S)} 
When the spectral method is applied with the standard perturbation theory, 
it is well-known that Condition (S) can be reduced to more practical non-arithmetic or nonlattice assumptions, see e.g \cite{gui} \cite{guihar} \cite{hulo}. These reductions are based on some spectral arguments, and on simple properties of strict convexity. In this section, we generalize these results under the next Condition ($\widehat K$), close to ($\widetilde K$) of Section 4, but involving the whole family $\{Q(t), t\in\R^d\}$ and an additional condition on the essential spectral radius of $Q(t)$. Condition ($\widehat K$) will be satisfied in all our examples. 

\noindent {\it Condition ($\widehat K$):  $Q$ satisfies Condition (K1) (of Section 1) 
on $\cB$, and there exists  a Banach space $\widehat{\cB}$ such that $\cB\hookrightarrow\widehat{\cB}$, $Q(t)\in\cL(\cB)\cap\cL(\widehat{\cB})$  for each $t\in\R^d$, and: 
\footnote{As in ($\widetilde K$), the fact that $\widehat{\cB}$ is complete is not necessary, but always satisfied in practice. Contrary to ($\widetilde K$), it is not convenient for the next statements to assume $\widehat{\cB}\hookrightarrow\L^1(\pi)$ (except for Proposition~12.4). } 

\noindent ($\widehat{K2}$)  $\ \forall t\in{\mathbb R}^d,
\ \displaystyle\lim_{h\r0}\|Q(t+h)-Q(t)\|_{\scriptsize{\cB,\widehat{\cB}}}=\, 0$  

\noindent and, for all compact set $K_0$ in $\R^d$, there exists $\kappa\in (0,1)$ such that:  

\noindent ($\widehat{K3}$) $\exists C> 0,\ \forall n\geq1,\ \forall f\in\cB,\ \forall t\in K_0, \ \
\|Q(t)^nf\|_{\scriptsize \cB} \leq C\, \kappa^n\, \|f\|_{\scriptsize \cB} + C\, \|f\|_{\scriptsize{\widehat{\cB}}}$ 

\noindent ($\widehat{K4}$) $\ \forall t\in K_0$, $\, r_{ess}(Q(t)) \leq \kappa$. }

\noindent Clearly, if $\widehat{\cB}\hookrightarrow\L^1(\pi)$, then ($\widehat K$) implies ($\widetilde K$) of Section 4. Besides, when $\widehat{\cB} = \L^1(\pi)$, the condition introduced  in Remark (a) of Section 4 implies ($\widehat{K2}$). 

\noindent We also need the next assumption (fulfilled in practice under Condition (K1), see Rk. below):  {\it 

\noindent (P) We have, for any $\lambda\in\C$ such that $|\lambda| \geq 1$, and for any nonzero element $f\in \cB$: 
$$\big[\, \exists n_0,\ \forall n\ge n_0,\ \ |\lambda|^n|f| \leq Q^n|f|\, \big]\ \ \Rightarrow\ \ 
\big[\, |\lambda| =1\ 
\mbox{and}\ |f|\leq\pi(|f|)\, \big].$$ 
\noindent The previous inequalities hold, everywhere on $E$  if we have $\cB\subset\cL^1(\pi)$, and $\pi$-almost surely on $E$ if we have $\cB\subset\L^1(\pi)$. }

\noindent If $\cB\subset\L^1(\pi)$, we shall say 
that $w$ is a bounded element in $\cB$ if 
$w\in\cB\cap\L^\infty(\pi)$. 

\noindent {\bf A non-arithmetic condition on $\xi$.} {\it We shall say that $(Q,\xi)$, or merely $\xi$, is arithmetic w.r.t.~$\cB$
(and non-arithmetic  w.r.t.~$\cB$ in the opposite case) if there exist $t\in\R^d$, $t\neq0$, $\lambda\in\C$, $|\lambda|=1$, 
a $\pi$-full $Q$-absorbing set $A\in\cE$, and a bounded element $w$ in $\cB$ such that 
$|w|$ is nonzero constant on $A$, satisfying: }
$$(*)\ \ \ \ \displaystyle \forall x\in A,\ \ 
e^{i\langle t,\xi(y)\rangle} w(y) = \lambda w(x)\ \ 
Q(x,dy)-\mbox{\it a.s.}.$$ 
\noindent{\bf Proposition 5.3.} {\it Under the assumptions ($\widehat{K}$) and (P), Condition (S) holds on $\cB$ if and only if $\xi$ is non-arithmetic w.r.t.~$\cB$. } 

\noindent In the usual spectral method, this statement is for instance established in \cite{hulo} (Prop.~V.2) (under some additional conditions on $\cB$). The proof of Proposition 5.3, which is an easy extension of that in \cite{hulo}, is presented in Section 12.1. 

\noindent We now state a lattice-type criterion for (S) which is a natural extension of the i.i.d.~case and a well-known condition in the general context of Markov random walks. 

\noindent {\bf A nonlattice condition on $\xi$.} 
\noindent {\it We say that $(Q,\xi)$, or merely $\xi$, is lattice (and nonlattice in the opposite case) if  
there exist $a\in\R^d$, a closed subgroup $H$ in $\R^d$, $H\neq \R^d$, a $\pi$-full $Q$-absorbing set
$A\in\cE$, and a bounded measurable function $\theta\, :\, E\r\R^d$ such that }
$$(**)\ \  \ \ \ \forall x\in A, \ \ 
\xi(y)+\theta(y)-\theta(x)\in a+H\ \ Q(x,dy)-a.s..$$
\noindent{\bf Proposition 5.4.} {\it Assume that the assumptions ($\widehat{K}$) and (P) hold. If $\xi$ is nonlattice, then (S) holds on $\cB$. The converse is true when, for any real-valued 
measurable function $\psi$ on $E$, we have $e^{i\psi}\in\cB$ (or $Cl(e^{i\psi})\in\cB$). }

\noindent {\it Proof.} If (S) is not fulfilled, then $\xi$ is 
arithmetic w.r.t.~$\cB$, and one may assume that $w\in\cB$ in ($*$) is such that we have $|w|=1\ \pi$-a.s., so that we can write $w(x) = e^{ig(x)}$ for some measurable function $g : E\r [0,2\pi]$. Therefore, setting $\lambda = e^{ib}$, the property ($*$) is then equivalent to: 
$$ \forall x\in A,\ \ \langle t,\xi(y)\rangle + g(y) - g(x) -b  \in2\pi\Z\ \ \ Q(x,dy)-\mbox{a.s.}.$$
Now set $\theta(x) = g(x)\, \frac{t}{|t|_2^2}$, and $a = b\, \, \frac{t}{|t|_2^2}$.  Then we have ($**$) with $H = (2\pi\Z)\, \frac{t}{|t|_2^2}\oplus (\R\, t)^{\perp}$, so $\xi$ is lattice. Conversely, if $\xi$ is lattice, then, by considering ($**$) and $t\in H^{\perp}$, one can easily prove that ($*$) holds with $\lambda = e^{i\langle t,a\rangle}$ and $w(x) = e^{i\langle t,\theta(x)\rangle}$. Since $w\in\cB$, (S) is not fulfilled on $\cB$ (by Proposition 5.3.). \fdem 

\indent Proposition 5.4 will be specified in Section 12.2, where we shall investigate the following set: $G = \{t\in\R^d : r(Q(t))=1\}$. We conclude Section 5.2 by some further remarks. 

\noindent {\it On Conditions  ($\widehat{K3}$) ($\widehat{K4}$).} \\ 
If $(X_n)_{n\geq0}$ is strongly ergodic on a Banach lattice $\cB$ and if $\cB$ is such that  $e^{i\langle t,\, \xi\rangle}\cdot f\in\cB$ for all $t\in\R^d$ and $f\in\cB$, then we have ($\widehat{K3}$) with $\widehat{\cB} = \L^1(\pi)$ (see Rk.~(b) of Section~4). Moreover we have ($\widehat{K4}$) on $\cB$ according to \cite[Cor. 1.6]{rab-wolf}. \\[0.1cm]
The weighted Lipschitz-type spaces used in Section~11 for the iterative models are not Banach lattices, and in these models, the next remark will be helpful to prove ($\widehat{K4}$). Let us assume that ($\widehat{K3}$) is fulfilled with $\cB$ and $\widehat{\cB}$ satisfying the following property: for each $t\in\R^d$, $Q(t)({\mathcal S})$ is relatively compact in $(\widehat{\cB},\|\cdot\|_{\scriptsize{\widehat{\cB}}}) $, where $\mathcal S$ is the unit ball of $(\cB,\|\cdot\|_{\scriptsize \cB})$. Then it follows from  \cite{itm} \cite{hen2} that Condition ($\widehat{K4}$) automatically holds on $\cB$.  

\noindent{\it On Condition (P).} \\
Under Assumption (K1), the property (P) is for instance fulfilled in the following cases: \\[0.12cm]
- $\cB\subset\L^1(\pi)$, $\cB\hookrightarrow\L^1(\pi)$, and $\cB$ is dense in $\L^1(\pi)$. Indeed, since $Q$ is a contraction on $\L^1(\pi)$, one then obtains from (K1) that $\lim_nQ^n|f|= \pi(|f|)$ in $\L^1(\pi)$ for all $f\in\cB$, hence (P) (here, $f\neq0$  means that $\pi(|f|)\neq 0$). This case contains $\L^p(\pi)$, the Sobolev spaces,... \\[0.12cm]
- $\cB\subset\cL^1(\pi)$, $\cB\hookrightarrow\L^1(\pi)$, $\cB$ is stable under complex modulus
(i.e. $f\in{\mathcal B}\ \Rightarrow\ \vert f\vert
\in{\mathcal B}$), and $\delta_x\in\cB'$ for all $x\in E$. Indeed,  we then have by (K1): $\forall x\in E,\ \lim_n(Q^n|f|)(x) = \pi(|f|)$, hence (P) (here, $f\neq0$ means that $f(x)\neq 0$ for some $x\in E$). This case contains $\cB^{^\infty}$, the weighted (either supremum or Lipschitz-type) spaces, the space of bounded continuous functions, the space of functions of bounded variation (on an interval),.... \\[0.12cm]
- $\cB$ is the space of $\cC^k$ functions (on some nice $E$) equipped with its usual norm. Observe that, if $f\in\cC^k$, then $|f|$ is continuous on $E$. By using a density argument (with the supremum norm) and the property (K1) on $\cC^k$, one can easily see that $\lim_nQ^n|f| = \pi(|f|)$ uniformly on $E$, hence (P). 

\noindent{\it A case when $A=Supp(\pi)$ in ($*$) and ($**$).} \\
If ($\widehat{K}$) and (P) hold, if $\delta_x\in\cB'$ for all $x\in E$, and finally if all the functions of $\cB$ are continuous on the state space $E$ (assumed to be locally compact here), then Propositions 5.3-4 (and Proposition 12.4) apply with $A=Supp(\pi)$ in ($*$) and ($**$), where $Supp(\pi)$ is the support of $\pi$. This can be seen by an easy examination of the proof in Section 12.1.  

\noindent{\it Condition (S) and invertibility of $\Gamma$.}\\
\noindent  Let us just assume in this remark that $Q(t)\in\cL(\cB)$ for all $t\in\R^d$. If the conclusion of Proposition~2.4 
holds for some real-valued measurable function $g$ on $E$, then we clearly have 
($*$) with $w(\cdot) = e^{i\, g(\cdot)}$ and $\lambda=1$. Moreover, $(**)$
is satisfied with $a=0$, $\theta(x)=2\pi\{\frac{g(x)}{2\pi}\}\, {t\over\vert t\vert_2^2}$, and $H= (2\pi\Z) {t\over\vert t\vert_2^2} \oplus (\R\cdot t)^{\perp}  $, 
where $\{\cdot\}$ stands for the fractionary part. Condition (S) on any space $\cB$ containing $w$ is then false because, in this case, the above mentioned equality ($*$) easily implies that $r(Q(t))\geq 1$, see Lemma 12.2. One can deduce the following facts from the previous remarks and the results of Section 2. \\[0.1cm] 
If the hypotheses of Corollary 2.1 hold on some space $\cB_2$ and if $e^{i\, \psi(\cdot)}\in\cB$ for all $\psi\in\cB_2$, then we have the following implications, in which $\Gamma$ denotes the covariance matrix of Section~2 (the above condition on $\cB$ is unnecessary for the last implication):  \\[0.1cm]
\indent {\it $\ \ \ \ \ \ \ \ \ \ \ \ \ \ \ \ \ \ \ \ \ \ \ $ 
Condition (S) on $\cB\ \Rightarrow\ \Gamma$ is invertible}  \\[0.1cm] 
\indent {\it $\ \ \ \ \ \ \ \ \ \ \ \ \ \ \ \ \ \ \ \ \ \ \ $ 
Non-arithmeticity w.r.t.~$\cB\ \Rightarrow\ \Gamma$ is invertible}  \\[0.1cm] 
\indent {\it $\ \ \ \ \ \ \ \ \ \ \ \ \ \ \ \ \ \ \ \ \ \ \ $ 
$(Q,\xi)$ is nonlattice $\Rightarrow\ \Gamma$ is invertible.}   
%==================================
\subsection{(LLT) for the strongly ergodic Markov chains on $\L^2(\pi)$} 
Let us suppose that $(X_n)_{n\geq0}$ is a strongly ergodic Markov chain on $\L^2(\pi)$ (Ex.~1 of Section~1). If $\pi(|\xi|_2^2)<+\infty$, then $(n^{-\frac{1}{2}}S_n)_n$ converges in distribution to a normal distribution $\cN(0,\Gamma)$ (see Section 2). 

\noindent {\bf Corollary 5.5.} {\it Let us assume that $\pi(|\xi|_2^2)<+\infty$, that  $\xi$ is nonlattice, that 
$\mu = \pi$, and that $h\in\L^r(\pi)$ for some $r>1$. Then we have (LLT) for each function $f$ in $\L^p(\pi)$ provided that $p > \frac{r}{r-1}$. }

\noindent {\it Proof.} Let $r'=\frac{r}{r-1}$, and $p > r'$. From Proposition 4.1 and Lemma 4.2, we have ($\widehat{K2}$) and ($\widehat{K3}$) (thus ($\widetilde{K}$)) with  $\cB=\L^p(\pi)$ and $\widehat{\cB} = \L^{r'}(\pi)$. Note that $\widehat{\cB}' = \L^r(\pi)$. Since $\cB=\L^p(\pi)$ is a Banach lattice, we have ($\widehat{K4}$) on $\cB=\L^p(\pi)$ by \cite[Cor. 1.6]{rab-wolf}. Finally, from Proposition~5.4, Condition (S) on $\cB=\L^p(\pi)$ is fulfilled under the nonlattice assumption. Corollary~5.5 can be then deduced from Theorem~5.1. \fdem 

\noindent The property ($\widehat{K4}$) on $\cB=\L^p(\pi)$ has been above derived from the general statement \cite[Cor. 1.6]{rab-wolf} which is based on some sophisticated arguments of the theory of positive operators acting on a Banach lattice. Below, we present a simpler proof of this fact in the special case of the uniformly ergodic Markov chains. By repeating some arguments of \cite{hen4}, we are going to see that ($\widehat{K4}$) on $\cB=\L^p(\pi)$ then follows from Doeblin's condition. 

\noindent Let us assume that $(X_n)_{n\geq0}$ is uniformly ergodic (i.e.~we have (K1) on $\cB^{^\infty}$). Then the so-called Doeblin condition holds (use (K1) on $\cB^{^\infty}$): there exist $\ell\geq1$, $\eta>0$, and $\rho<1$ such that 
$$(\ \pi(A)\leq \eta\ )\ \Rightarrow\ (\ \forall x\in E,\ Q^\ell(x,A)\leq \rho^\ell\ ).$$
\noindent{\bf Proposition 5.6.} {\it Let $p\in(1,+\infty)$. If $\xi$ is any $\R^d$-valued measurable function on $E$, then we have: $\forall t\in\R^d,\ r_{ess}(Q(t)) \leq \rho^{\frac{p-1}{p}}$. } 

\noindent {\it Proof of Proposition 5.6.} Let $\|\cdot\|_p = \pi(|\cdot|^p)^{\frac{1}{p}}$ denote the norm on  $\L^p(\pi)$. We also use the notation $\|\cdot\|_p$ for the operator norm on $\L^p(\pi)$. Let $q$ be such that $\frac{1}{p}+\frac{1}{q} = 1$. 

\noindent {\bf Lemma 5.7.} {\it There exist 
a nonnegative bounded measurable function $\alpha$ on $E\times E$ 
and a positive kernel 
$S(x,dy),\ x\in E$, such that $Q^\ell(x,dy) = \alpha(x,y)d\pi(y) + S(x,dy)$ and $\|S\|_p \leq \rho^{\frac{\ell}{q}}$. } 

\noindent {\it Proof.} Let us summarize the beginning of the proof in \cite{hen4} (Lemma III.4): using the differentiation of measures, 
there exist a nonnegative measurable function $\alpha'$ on $E\times E$ and a positive kernel $S'(x,dy)$ such that, for all $x\in E$, 
we have $Q^\ell(x,dy) = \alpha'(x,y)d\pi(y) + S'(x,dy)$, with $\pi(C_x)=0$ and $S'(x,E\setminus C_x)=0$ for some $C_x\in\cE$. 
Set $\alpha = \alpha'\, 1_{\{\, \alpha'\leq \eta^{-1}\, \}}$, and for $x\in E$, let 
$L_x = \{\, y\in E :  \alpha'(x,y)> \eta^{-1}\, \} \setminus C_x$. Then $Q^\ell(x,dy) = \alpha(x,y)d\pi(y) + S(x,dy)$ with 
$S(x,A) = Q^\ell(x,A\cap(C_x\cup L_x))$. We have 
$$\forall x\in E,\ \ \ 1 \geq Q^\ell(x,L_x) \geq \int_{L_x} \alpha'(x,y)\, d\pi(y) \geq \eta^{-1}\,\pi(L_x),$$
thus $\pi(L_x\cup C_x) = \pi(L_x) \leq \eta$, so that $Q^\ell(x,L_x\cup C_x)\leq \rho^\ell$. \\
Now let $f\in \L^p(\pi)$. We have $Sf(x) = \int_{C_x\cup L_x}f(y)\, Q^\ell(x,dy)$, and  from H\"older's inequality 
w.r.t.~the probability 
measure $Q^\ell(x,dy)$, we have 
$$\|Sf\|_p^p = \int_E\bigg|\int_E f(y)\, 1_{C_x\cup L_x}(y)\, Q^\ell(x,dy)\bigg|^pd\pi(x) \leq 
\int_E  Q^\ell|f|^p(x)\, Q^\ell(x,C_x\cup L_x)^{\frac{p}{q}}\, d\pi(x),$$
hence $\|Sf\|_p^p \leq (\rho^\ell)^{\frac{p}{q}}\, \pi(Q^\ell|f|^p) = 
(\rho^\ell)^{\frac{p}{q}}\, \pi(|f|^p)$ which is the stated estimate on $\|S\|_p$. \fdem

\noindent Now let us  prove $r_{ess}(Q(t)) \leq \rho^{\frac{1}{q}}$ for all $t\in\R^d$.  
Since $|Q(t)^\ell f| \leq Q^\ell|f|$, there exists a complex-valued measurable function $\chi_{t}$ on 
$E\times E$ such that $Q(t)^\ell(x,dy) = \chi_{t}(x,y)\, Q^\ell(x,dy)$ with $|\chi_{t}| \leq 1$. So, by Lemma 5.7, 
$$Q(t)^\ell(x,dy) = \chi_{t}(x,y)\, \alpha(x,y)d\pi(y) + \chi_{t}(x,y)\, S(x,dy) := \alpha_t(x,y)d\pi(y) + S_t(x,dy),$$ 
and, since $\alpha_t(\cdot,\cdot)$ is bounded, the associated kernel operator is compact on $\L^p(\pi)$ \cite{ds}.  
Recall that, if $T$ is a bounded operator on a Banach space $\cB$, then 
$r_{ess}(T) = \lim_n(\inf \|T^n-V\|_{\scriptsize \cB})^{\frac{1}{n}}$ where the infimum is considered over the ideal of 
compact operators $V$ on  $\cB$. This yields $r_{ess}(Q(t)^\ell) = r_{ess}(S_t) \leq r(S_t) \leq 
\|S_t\|_p \leq \|S\|_p \leq \rho^{\frac{\ell}{q}}$ (Lem.~5.7). Hence  $r_{ess}(Q(t)) \leq \rho^{\frac{1}{q}}$. \fdem 
%=====================================================================
\section{A one-dimensional uniform Berry-Esseen theorem} 
Here we assume $d=1$ (i.e.~$\xi$ is real-valued), we denote by $\cN$ the distribution function 
of $\cN(0,1)$, we suppose that Hypothesis (CLT)
of Section 5.1 holds with $\Gamma = \sigma^2 >0$, and we set:
$$\forall u\in\R,\ \ \Delta_n(u) =  \left|\, \P_\mu(\frac{S_n}
{\sigma\sqrt n}\leq u) - \cN(u)\, \right|, \ \ \ \ \mbox{and}\ \ \ \ \ 
\Delta_n = \sup_{u\in\R}\, \Delta_n(u).$$
Theorem 6.1 and Proposition 6.2 below have been already presented in \cite{ihp2}, we state them again for completeness. 
The next application to the Markov chains with spectral gap on $\L^2(\pi)$ is new. Comparisons with prior works are presented in \cite{ihp2}, 
they will be partially recalled below and in Sections 10-11. 
%=============================

\noindent {\bf A general statement.} 

Let us reinforce Condition (CLT) by the following one: 

\noindent {\it Condition (CLT'): $\displaystyle \ \exists C>0,\ \forall t\in[-\sqrt n,\sqrt n], \ \ 
\big|\E_\pi[e^{it\frac{S_n}{\sigma\sqrt n}}]\, - \, e^{-\frac{t^2}{2}}\big| \leq C\frac{|t|}{\sqrt n}$. } 

\noindent{\bf Theorem 6.1} \cite{ihp2}. {\it Assume that (CLT') holds, and that Condition ($\widetilde{K}$) (of Section 4) holds ~w.r.t.~$\cB$, $\widetilde{\cB}$, with the additional following conditions: we have (K1) (of Section~1) on $\widetilde{\cB}$, and  
$\|Q(t)-Q\|_{\scriptsize{\cB,\widetilde{\cB}}}=O(|t|)$. Then we have $\displaystyle\Delta_n = O(n^{-\frac{1}{2}})$ for any $\mu\in\widetilde{\cB}'$. }

\noindent {\it Proof (sketch)}. See \cite{ihp2} for details. The conclusions of Theorem (K-L) are satisfied. As in Lemma 4.3, let us set $\ell(t) = \mu(\Pi(t)1_E)$. In order to copy the Fourier techniques used for the i.i.d.~Berry-Esseen theorem (see \cite{fel} \cite{dur}), we have to improve Lemmas 4.3 and 5.2 as follows: \\[0.1cm]
(a) $\displaystyle \sup_{t\in {\cal O}}\frac{|\ell(t)-1|}{|t|} < +\infty\ \ $ and 
$\ \ \displaystyle \sup_{t\in {\cal O}}\, \frac{1}{|t|}\, \big|\E_\mu[e^{i t S_n}] -  \lambda(t)^n\, \ell(t)\big| =  O(\kappa^{n})$ \\[0.1cm]
(b) $\lambda(u) = 1 -  \sigma^2\frac{u^2}{2} + O(u^3)$ near $u=0$. \\[0.1cm]
\noindent Assertion (a) cannot be derived from the Keller-Liverani theorem (even by using the precise statements of \cite{keli}). However one can proceed as follows. As in the standard perturbation theory 
\cite{ds}, the perturbed eigen-projection $\Pi(t)$ in Theorem (K-L) can be expressed as the line integral of $(z-Q(t))^{-1}$ over a suitable oriented circle centered at $\lambda=1$ (see Section 7.2). By using the formula \\[0.2cm]
\indent $\displaystyle \ \ \ \ \ \ \ \ \ \ \ \ \ \ \ \ \ \ (z-Q(t))^{-1} -  (z-Q)^{-1} = (z-Q)^{-1}\, \big[Q(t)-Q\big]\, (z-Q(t))^{-1}$, \\[0.2cm]
the last assertion in Theorem (K-L), the assumption $\|Q(t)-Q\|_{\scriptsize{\cB,\widetilde{\cB}}}=O(|t|)$, and finally the fact that 
(K1) holds on $\widetilde{\cB}$, one can then conclude that $\|\Pi(t)-\Pi\|_{\scriptsize{\cB,\widetilde{\cB}}} = O(|t|)$. Hence the desired property for $\ell(t)$. The second assertion in (a) can be established similarly by using Formula (CF) of Section 3 and the second line integral given in Section 7.2. \\
To get (b), one may repeat the short proof of Lemma 5.2 
by starting here from the property $\lambda(\frac{t}{\sqrt n})^n - e^{-\frac{\sigma^2}{2}t^2} = O(\frac{|t|}{\sqrt n})$ 
which follows from (CLT') and (a). One then obtains  
$(\frac{\sqrt n}{t})^2 (\lambda(\frac{t}{\sqrt n}) - 1) + \frac{\sigma^2}{2} =  O(\frac{|t|}{\sqrt n})$, and 
setting $u=\frac{t}{\sqrt n}$, this leads to the expansion (b) (see Lem.~IV.2 in \cite{ihp2}). \fdem 

\noindent {\bf A sufficient condition for (CLT').} 

Actually, one of the difficulties in the previous theorem is to show Hypothesis (CLT'). By the use of martingale techniques 
derived from \cite{jan}, the first named author showed in \cite{ihp2} the next statement. 

\noindent{\bf Proposition 6.2 } \cite{ihp2}. {\it We have (CLT') when the two following conditions hold: \\[0.15cm]
(G1) $\breve\xi = \sum_{n=0}^{+\infty} Q^n\xi$ absolutely 
converges  in $\L^3(\pi)$. \\[0.15cm]
(G2) $\sum_{p=0}^{+\infty} Q^p\psi$ absolutely converges in $\L^{\frac{3}{2}}$, where 
$\psi = Q(\breve\xi^2) - (Q\breve\xi)^2 - (\pi(\breve\xi^2) - \pi((Q\breve\xi)^2)\, 1_E$. }

\noindent Let us notice that $\breve\xi$ is the solution of the Poisson equation $\breve\xi - Q\breve\xi = \xi$, already 
introduced in Gordin's theorem (Section 2). Also observe that the above function $\psi$ can be expressed as 
$\psi = Q(\breve\xi^2) - (Q\breve\xi)^2 - \sigma^2\, 1_E$, where $\sigma^2$ is the asymptotic variance of Gordin's theorem. 

\noindent {\it About the practical verification of (G1) (G2).} \\
In practice, one often proceeds as follows to verify the two above conditions. 
Since $\pi(\xi) = 0$, Condition (G1)  holds if $Q$ is strongly ergodic w.r.t.~some $\cB\hookrightarrow\L^3(\pi)$ and 
if $\xi\in\cB$. If moreover $Q$ is strongly ergodic w.r.t.~some $\cB_2 \hookrightarrow\L^{\frac{3}{2}}(\pi)$ 
containing all the 
functions $g^2$ with $g\in\cB$, then  Condition (G2) holds. Indeed, under these hypotheses, $\breve\xi\in\cB$, thus $\psi\in\cB_2$, 
and, since $\pi(\psi)=0$, the series $\sum_{p=0}^{+\infty} Q^p\psi$ absolutely converges in $\cB_2$, thus in 
$\L^{\frac{3}{2}}(\pi)$. 

\noindent Condition~(G2) is the functional version  of the projective assumption 
$\sum_{n\geq0}\big\|\E[Z_n^2\, |\, \cF_0] - \E[Z_0^2]\big\|_{\L^{\frac{3}{2}}} < +\infty$ used for stationary martingale difference  sequences $(Z_n)_n$: under this condition, the uniform Berry-Esseen theorem at rate $n^{-\frac{1}{4}}$ is established in \cite{jan} (Chap.~3) for such bounded sequences. In \cite{ded-rio, ded-mer-rio}, this projective assumption (extended to $\L^p$ in \cite{ded-mer-rio}) provides   the expected Berry-Esseen theorem in term of Wasserstein's distances. 

\noindent {\bf Application to the strongly ergodic Markov chains on $\L^2(\pi)$.} 

Let us assume that $(X_n)_{n\geq0}$ is a strongly ergodic Markov chain on $\L^2(\pi)$ (Ex.~1 of Section~1). In the stationary case (i.e.~$\mu=\pi$), since $(X_n)_{n\geq0}$ is strongly mixing (see \cite{rosen}), Bolthausen's theorem \cite{bolt} yields the estimate $\Delta_n = O(n^{-\frac{1}{2}})$ if $\pi(|\xi|^p)<+\infty$ for some $p>3$. In the special case of uniform ergodicity, Nagaev's work \cite{nag2}, and some of its extensions (see e.g \cite{datta}), provide the previous estimate in the non-stationary case, but under the strong moment condition $\sup_{x\in E} \int_E|\xi(y)|^3\, Q(x,dy) < +\infty$. The next statement only requires the expected third-order moment condition. 

\noindent{\bf Corollary 6.3.} {\it If $\pi(|\xi|^3)<+\infty$ and $\mu = \phi\, d\pi$, with some $\phi\in\L^3(\pi)$, 
then $\Delta_n = O(n^{-\frac{1}{2}})$. }

\noindent{\it Proof.} Set $\L^p=\L^p(\pi)$. We have (K1) on $\L^3$ and $\L^{\frac{3}{2}}$, see \cite{rosen}. So Conditions (G1) (G2), hence (CLT'), are fulfilled (use the above remark with $\cB_2 = \L^{\frac{3}{2}}$). Besides, we have ($\widetilde{K3})$ with $\cB=\L^3$ and $\widetilde{\cB} = \L^{\frac{3}{2}}$ (Prop.~4.1). Finally we have $\|Q(t)-Q\|_{\scriptsize{\L^3,\L^{\frac{3}{2}}}}=O(|t|)$. Indeed, let $f\in\L^3$. Using $|e^{ia}-1|\leq |a|$, one gets \\[0.12cm]
\indent $\pi (\, |Q(t)f-Qf|^{\frac{3}{2}}\, ) \leq   \pi (\, |Q(|e^{it\xi}-1|\, |f|)|^{\frac{3}{2}}\, ) 
\leq  |t|^{\frac{3}{2}}\, \pi (\, Q(|\xi|^{\frac{3}{2}}\, |f|^{\frac{3}{2}})\, ) 
=   |t|^{\frac{3}{2}}\, \pi (|\xi|^{\frac{3}{2}}\, |f|^{\frac{3}{2}})$, \\[0.12cm]
and the Schwarz inequality yields 
$\|Q(t)f-Qf\|_{\frac{3}{2}}  \leq |t|\, (\pi (|\xi|^{\frac{3}{2}}\, |f|^{\frac{3}{2}}))^{\frac{2}{3}} \leq |t|\, \|\xi\|_{_3}\, \|f\|_{_3}$. 
We have proved that the hypotheses of Theorem 6.1 are fulfilled with $\cB=\L^3$ and $\widetilde{\cB} = \L^{\frac{3}{2}}$. \fdem 
%==================================================================================
%\noindent{\bf 7. Regularity of the eigen-elements of $Q(t)$.} 
\section{Regularity of the eigen-elements of the Fourier kernels}
The goal of this section is to present an abstract operator-type Hypothesis, called $\cC(m)$, ensuring
that the dominating eigenvalue $\lambda(t)$ and the associated eigen-elements of $Q(t)$ have $m$ continuous derivatives 
on some neighbourhood ${\cal O}$ of $0$. The usual spectral method already exploited this idea 
by considering the action of $Q(t)$ on a single space, but as illustrated in Section 3, 
the resulting operator-moment conditions may be very restrictive in practice.  
The use of a ``chain'' of spaces developed here 
enables to greatly weaken these assumptions.  \\[0.1cm]
\indent As a first example we shall see in Section~7.3 
that, for the strongly ergodic Markov chains on $\L^2(\pi)$, Hypothesis $\cC(m)$ reduces to $\pi(|\xi|_2^\alpha) < +\infty$   
for some $\alpha>m$. This condition is slightly stronger than the assumption $\pi(|\xi|_2^m) < +\infty$ of the i.i.d.~case ensuring that the common characteristic function has $m$ continuous derivatives. 
But it is much weaker than the condition $\sup_{x\in E}(Q|\xi|_2^m)(x) < +\infty$ of the usual spectral method (see Section 3).  Other simple reductions of  $\cC(m)$ will be obtained in Sections 10-11 for Examples 2-3 of Section 1. \\[0.1cm]
\indent Roughly speaking one can say that Hypothesis $\cC(m)$ below 
(together with possibly the non-arithmeticity condition) allows to extend to 
strongly ergodic Markov chains the classical i.i.d.~limit theorems established with Fourier techniques. This will be illustrated in Sections~8-9 by a one-dimensional Edgeworth expansion and a multidimensional Berry-Esseen type theorem in the sense of the Prohorov metric. This is also exploited in \cite{guiher} to prove a multidimensional renewal theorem. \\[0.1cm]
\indent Before dealing with the regularity of the eigen-elements of $Q(t)$, we investigate that of the function $t\mapsto(z-Q(t))^{-1}$, where $(z-Q(t))^{-1}$ is seen as an element of $\cL(\cB,\widetilde{\cB})$ for suitable spaces $\cB$ and $\widetilde{\cB}$. 
%===============================================
%\noindent{\bf 7.1. Regularity of $(z-Q(\cdot))^{-1}$.} 
\subsection{Regularity of $(z-Q(\cdot))^{-1}$}
\noindent Let ${\cal O}$ be an open subset of $\R^d$, let $X$ be a vector normed space.
Then, for $m\in\N$,  
we shall say that $U\in\cC^m({\cal O},X)$ if 
$U$ is a function from ${\cal O}$ to $X$ which admits $m$ continuous derivatives. 
For convenience, $\cC^\ell({\cal O},\cB_1,\cB_2)$ will stand for
$\cC^\ell({\cal O},\cL(\cB_1,\cB_{2}))$. In view of the probabilistic applications of Sections~8-9, $Q(t)$ still denotes the Fourier kernels defined in Section~1, and the Banach spaces $\cB$, $\tilde {\cB}$, $\cB_\theta$ considered below satisfy the conditions stated before (K1) in Section~1.  
Let $\cB\hookrightarrow\widetilde{\cB}$, and let $m\in\N^*$.  

\noindent  {\bf Hypothesis $\cC(m)$.} 
{\it There exist a subset $I$ of $\mathbb R$ and 
a family of spaces $({\cB }_\theta,\ \theta\in I)$ containing $\cB$, 
$\widetilde{\cB}$, such that $\cB_\theta\hookrightarrow\widetilde{\cB}$ for all $\theta\in I$, and there exist two functions $T_0:I\rightarrow \R$ and $T_1:I\rightarrow \R$ such that, for all 
$\theta\in I$, there exists a neighbourhood 
$\cV_\theta$ of $0$ in $\R^d$ such that we have for 
$j=1,...,m$: \\[0.15cm]
\noindent {\bf (0)} $ [T_0(\theta)\in I \ 
\Rightarrow\ \ \cB_\theta\hookrightarrow \cB_{T_0(\theta)}]\ \ $ and
$\ \ [T_1(\theta)\in I \ \Rightarrow\ \ \cB_\theta\hookrightarrow \cB_{T_1(\theta)}]$ 
\\[0.15cm]
\noindent {\bf (1)} If $T_0(\theta)\in I$, then
$Q(\cdot)\in\cC^0(\cV_\theta,\cB_{\scriptsize\theta},\cB_{\scriptsize T_0(\theta)})$ \\[0.15cm]
\noindent {\bf (2)} If $\theta_j:=T_1(T_0T_1)^{j-1}
(\theta)\in I$, then 
$\ Q(\cdot)\in \cC^j(\cV_\theta,\cB_{\scriptsize\theta},\cB_{\scriptsize \theta_j})$ \\[0.15cm]
\noindent {\bf (3)}  $\ Q(\cdot)$ 
satisfies Hypothesis (K) of Section 4 on $\cB_\theta$ \\[0.15cm]
\noindent {\bf (4)}  There exists 
$a\in\bigcap_{k=0}^m\left[T_0^{-1}(T_0T_1)^{-k}(I)
\cap (T_1T_0)^{-k}(I)\right]$ such that we have $\cB=\cB_a\ $ and $\ \widetilde{\cB} = \cB_{(T_0T_1)^mT_0(a)}$. }

\noindent To fix ideas, let us introduce a more restrictive but simpler hypothesis~:

\noindent  {\bf Hypothesis $\cC'(m)$.} {\it  
There exist $A>m$ and a family of spaces $({\cB }_\theta,\ \theta\in[0,A])$
such that ${\cB}_0=\cB$, ${\cB}_A=\widetilde{\cB}$ and, for all $\theta,\theta'\in[0,A]$
with $0\le \theta<\theta'\le A$, we have~:\\[0.15cm]
\noindent {\bf (a)}\ \ ${\cB}_\theta \hookrightarrow  {\cB}_{\theta'}\hookrightarrow\widetilde{\cB}$,\\[0.15cm]
\noindent {\bf (b)}\ there exists a neighbourhood
$\cV=\cV_{\theta,\theta'}$ of $0$ in $\R^d$ such that,
for any $j\in\{0,...,m\}$ with $j<\theta'-\theta$,\ 
we have $Q\in\cC^j  \left({\cV},\cB_\theta,{\cB}_{\theta'}\right)$,\\[0.15cm]
\noindent {\bf (c)} $\ Q(\cdot)$ 
satisfies Hypothesis (K) of Section 4 on $\cB_\theta$. }

\noindent It is easy to see that Hypothesis $\cC'(m)$ implies Hypothesis $\cC(m)$
(by taking $a=0$, $T_0(x)=x+\varepsilon$
and $T_1(x)=x+1+\varepsilon$ for some well chosen $\varepsilon>0$). Actually Hypothesis $\cC'(m)$ will be satisfied in all our examples, but 
Hypothesis $\cC(m)$ is more general and, despite its apparent complexity, might be more natural to establish than hypothesis $\cC'(m)$ 
(see the end of Section 7.3).

\noindent Let us come back to Hypothesis $\cC(m)$.
The condition on $a$ in (4) means that $a$, $T_0a$, $T_1T_0a$, 
$T_0T_1T_0a$,...,$(T_0T_1)^mT_0(a)$ belong to $I$, and from (0), it follows that the corresponding family of $\cB_\theta$'s is increasing with respect to the continuous  embedding. 
In particular, $\theta:=T_0(a)$ and $\theta_m := T_1(T_0T_1)^{m-1}(\theta)$ are in $I$,  
therefore we have $Q(\cdot)\in\cC^m \left(\cV_{\theta},\cB_{\theta},
\cB_{\theta_m}\right)$ by (2).
It then follows that $Q(\cdot)\in\cC^m \big(\cV_{\theta},\cB
,\widetilde{\cB}\big)$. In practice, we may have 
$Q(\cdot)\in\cC^m \big(\cV_{\theta},\cB,\cB_{T_1^m(a)}\big)$, but the introduction of 
$T_0$ will enable us to get $(z-Q(\cdot))^{-1}\in
\cC^m \big({\cal O},\cB,\widetilde{\cB}\big)$ for some neighbourhood $\cal O$ of $t=0$ and for suitable $z\in\C$.

\noindent {\it Notation.} Recall that we set 
$\cD_{\kappa} = \{z\in\C:\vert z\vert\ge\kappa,\ 
      \vert z-1\vert\ge(1-\kappa)/2\}$ for any 
$\kappa\in(0,1)$. Under Hypothesis $\cC(m)$, 
we have (K) on $\cB$, so from Theorem (K-L) of 
Section 4, if $t$ belongs to some neighbourhood 
$\cU_a$ of $0$ in $\R^d$ and if 
$z\in\cD_{\kappa_a}$ for some $\kappa_a\in(0,1)$, 
then $(z-Q(t))^{-1}$ is a bounded operator on $\cB$, 
and we shall set 
$R_z(t) = (z-Q(t))^{-1}$. It is worth noticing that 
we also have $R_z(t)\in\cL(\cB,\tilde {\cB})$ for all $t\in \cU_a$ and $z\in\cD_{\kappa_a}$. 
In the case $d\geq 2$, for $t=(t_1,\ldots,t_d) $, $R_z^{(\ell)}(t)$ will stand for 
any partial derivative of the form $\displaystyle \frac{\partial^\ell R_z}{\partial t_{i_1}\cdots\partial t_{i_\ell}}(t)$. 

\noindent {\bf Proposition 7.1.} {\it  Under 
Hypothesis $\cC(m)$,   
there exist a neighbourhood $\cal O\subset\, $$\cU_a$ of $0$ in $\R^d$ and $\tilde\kappa\in (\kappa_a,1)$
such that $R_z(\cdot)\in\cC^m({\cal O},\cB,\tilde {\cB})$ 
for all $z\in{\cal D}_{\tilde\kappa}$, and }
$$\cR_\ell := \sup\{\|\, R_z^{(\ell)}(t)\|_{\scriptsize{\cB,\widetilde{\cB}}},\, z\in{\cal D}_{\tilde\kappa},\, t\in {\cal O}\, \}<+\infty,\ \ \ \ell=0,\ldots,m.$$
The proof of Proposition 7.1 is presented in Appendix A under a little bit more abstract setting. It is based on general and elementary derivation arguments. Similar statements concerning the Taylor expansions of $(z-Q(\cdot))^{-1}$ at $t=0$ are developed in \cite{aap,gouliv,seb-08}.

\noindent {\it Remarks.} \\ 
\noindent (a) In hypothesis ${\cC}(m)$,
the set $I$ can be reduced to the following finite set~:
$$\big\{ a,\, T_0a,\, T_1T_0a,\, 
T_0T_1T_0a,\, \ldots, (T_0T_1)^mT_0(a)\big\}.$$
This remark will be of no relevance for checking $\cC(m)$ in our examples, but it will be important in the proof of Proposition~7.1 in order to define the set ${\cal O}$, the real number $\tilde\kappa$ , and finally the bounds $\cR_\ell$ (see the remark following Proposition~A in Appendix~A). 

\noindent (b) In our examples, the derivative condition (2) of Hypothesis $\cC(m)$ can be investigated by using the  partial derivatives  $\frac{\partial^j Q}{\partial t_{p_1}\cdots\partial t_{p_j}}(t)$, defined by means of the kernel  \\
\indent $\displaystyle \ \ \ \ \ \ \ \ \ \ \ \ \ \ \ \ \ \ 
Q_{(p_1,\ldots,p_j)}(t)(x,dy) = i^j\left(\prod_{s = 1}^j \xi_{p_s}(y)\right)e^{i\langle t,\xi(y)\rangle}\, Q(x,dy)$. \\[0.15cm]
Actually, in our examples, we shall verify $\cC(m)$ in the case $d=1$ (for the sake of simplicity), and we shall simply denote by $Q^{(k)}$ the $k$-th derivative of $Q(\cdot)$ occurring in $\cC(m)$, which is defined for $k=0,\ldots,m$ by the kernel  \\[0.15cm]
\indent $\displaystyle \ \ \ \ \ \ \ \ \  \ \ \ \ \ \ \ \ \ 
Q_k(t)(x,dy) = i^k\xi(y)^ke^{it\xi(y)}\, Q(x,dy)\ \ $ ($t\in\R,\ x\in E$).  

\noindent (c) By $\cC(m)$, we know that $\frac{\partial^m Q}{\partial t_k^m}(0)\in\cL(\cB,
\widetilde{\cB})$ ($\, k=1,\ldots,d$). From $1_E\in\cB$, $\pi\in\widetilde{\cB}'$, it follows that $\pi(\frac{\partial^m Q}{\partial t_k^m}(0)1_E) = i^m\, \pi(Q\xi_k^m) = i^m\, \pi(\xi_k^m)$ is defined. So, in substance, Hypothesis $\cC(m)$ implies $\pi(|\xi|_2^m)<+\infty$ (this is actually true if $m$ is even). However, in our examples, we shall need some slightly more restrictive moment conditions to be able to prove $\cC(m)$. 
%==================================================
%\noindent{\bf 7.2. Regularity of the eigen-elements of $Q(\cdot)$.} 
\subsection{Regularity of the eigen-elements of $Q(\cdot)$} 
Suppose that Hypothesis $\cC(m)\ $ holds for some $m\in\N^*$, and as above let us use the notations of Proposition~7.1 and of Theorem (K-L) of Section 4 for $Q(t)$ acting on $\cB$:  
if $t\in{\cal U}_a$, $\lambda(t)$ is the dominating 
eigenvalue of $Q(t)$ and $\Pi(t)$ is the associated 
rank-one eigenprojection. Besides let us define in $\cL(\cB)$: 
$N(t) = Q(t)-\lambda(t)\Pi(t)$. Since $\Pi(t)Q(t) = Q(t)\Pi(t) = \lambda(t)\Pi(t)$, we have 
$$\forall n\geq1,\ \ N(t)^n =  Q(t)^n-\lambda(t)^n\Pi(t).$$ 
It follows from Theorem (K-L) that 
$Q(t)^n = \lambda(t)^n\Pi(t) + N(t)^n$, with 
$\|N(t)^n \|_{\scriptsize{\cB}} \leq C \kappa_a^n.$ \\[0.2cm] 
The operators $Q(t)$, $R_z(t)$, $\Pi(t)$ and $N(t)^n$ 
are viewed as elements of $\cL(\cB)$ when we appeal 
to the spectral theory, and as elements of 
$\cL(\cB,\widetilde{\cB})$ for stating our results of derivation. 

\noindent {\bf Corollary 7.2.} {\it Under Hypothesis $\cC(m)$,  
there exists a neighbourhood $\cV$ of $0$ in $\R^d$ such that~:\\[0.12cm]
(i) $\Pi(\cdot)\in\cC^m(\cV,\cB,\widetilde{\cB})$ \\[0.12cm]
(ii) for all $n\geq 1$, $N_n(\cdot):= 
 N(\cdot)^n\in\cC^m(\cV,\cB,\widetilde{\cB})$, and 
$$ \exists C>0,\ \forall n\geq1,\ \forall\ell=0,
\ldots,m : \ \ \ \sup_{t\in {\cal V}}\|N_n^{(\ell)}(t)\|_{\scriptsize{\cB,\widetilde{\cB}}} \leq C\tilde\kappa^n,$$ 
where $\tilde\kappa\in(0,1)$ is the real number of Proposition~7.1.  \\[0.12cm]
(iii) $\lambda(\cdot)\in\cC^m(\cV,\C)$. }

\noindent {\it Proof.} Let $t\in {\cal O}$, with ${\cal O}$ introduced in 
Proposition 7.1. \\
(i) As in the standard perturbation theory, the eigenprojection $\Pi(t)$ is defined in \cite{keli} by \\[0.1cm]
\indent $\displaystyle \ \ \ \ \  \ \ \ \ \  \ \ \ \ \  \ \ \ \ \ \ \ \  \ \ \ \ \  \ \ \ \ \  \ \ \
\Pi(t) = \frac{1}{2i\pi} \oint_{\Gamma_1} R_z(t)\, dz$, \\[0.1cm]
where this line integral is considered on the oriented circle $\Gamma_1$ centered at $z=1$, with radius $(1-\tilde\kappa)/2$ (thus $\Gamma_1\subset{\cal D_{\tilde\kappa}}$). Then, by Proposition 7.1, $\Pi(\cdot)\in\cC^m({\cal O},\cB,\widetilde{\cB})$. \\
(ii) In the same way, one can write \\[0.1cm]
\indent $\displaystyle \ \ \ \ \  \ \ \ \ \  \ \ \ \ \  \ \ \ \ \ \ \ \  \ \ \ \ \  \ \ \ \ \  \ \ \
N(t)^n = \frac{1}{2i\pi} \oint_{\Gamma_0} z^n\, R_z(t)\, dz$,\\[0.1cm]
where $\Gamma_0$ is here the oriented circle, centered at $z=0$, with radius $\tilde\kappa$ (thus $\Gamma_0\subset{\cal D}_{\tilde\kappa}$). 
By Proposition 7.1, we have $N_n(\cdot)\in
\cC^m({\cal O},\cB,\widetilde{\cB})$ with 
$N_n^{(\ell)}(t) = \frac{1}{2i\pi} \int_{\Gamma_0} z^n\, R_z^{(\ell)}(t)\, dz$ for $\ell=1,\ldots,m$. 
Hence the stated inequalities.  \\
\noindent (iii) Since $\lim_{t\r0}\pi(\Pi(t)
{ 1_E}) = \pi(\Pi{ 1_E}) = 1$ (by Th.~(K-L)), 
there exists a neighbourhood $\cV$ of $0$ contained in $\cal O$ such that
$\pi(\Pi(t)1_E)\neq0$ for any $t\in \cV$.
From $Q(t) = \lambda(t)\Pi(t) + N(t)$, it follows that 
$$\lambda(t) = \frac{\pi\big(Q(t){ 1_E} - N(t){1_E}\big)}{\pi(\Pi(t){1_E})}.$$
From the remark following the statement of 
$\cC(m)$, we have $Q(\cdot)\in\cC^m(\cV,\cB,
\widetilde{\cB})$ (with possibly $\cV$ reduced). Now, 
since ${ 1_E}\in\cB$ and  $N(\cdot)$, $\Pi(\cdot)$ 
are in 
$\cC^m(\cV,\cB,\widetilde{\cB})$, the functions 
$Q(\cdot){ 1_E}$, $N(\cdot){1_E}$, 
$\Pi(\cdot){1_E}$ are in $\cC^m(\cV,\widetilde{\cB})$. 
Finally, since $\pi\in\widetilde{\cB}'$, this gives (iii). \fdem
%==================================
\subsection{Hypothesis $\cC(m)$ for the strongly ergodic Markov chains on $\L^2(\pi)$}
Let us suppose that $(X_n)_{n\geq0}$ is a strongly ergodic Markov chain on $\L^2(\pi)$. Let  $m\in\N^*$, and let us investigate Hypothesis  
$\cC(m)$ by using a family $\{\cB_\theta = \L^\theta(\pi),\, r\leq\theta\leq s\}$ for some suitable $1<r<s$.  

\noindent {\bf Proposition 7.3.} {\it  If $\pi(|\xi|_2^\alpha) < +\infty$ with $\alpha>m$, then 
$\cC(m)$ holds with $\cB=\L^s(\pi)$ and $\widetilde{\cB}=\L^r(\pi)$ 
for any $s>{\alpha \over \alpha-m}$ and $1<r<{\alpha s\over \alpha+ms}$.}

\noindent We give the proof for $d=1$. The extension to $d\geq2$ 
is obvious by the use of partial derivatives. 

\noindent{\it Proof.} Let us notice that the condition on $s$ implies that ${\alpha s\over \alpha+ms} > 1$, so one may choose $r$ as stated, and we have $r<s$. 
Let $\varepsilon>0$ be such that $ r = {\alpha s\over \alpha+ms+\varepsilon(m+1)s}$. 
Let us prove $\cC(m)$ with $\cB_\theta = \L^\theta(\pi)$, $I=[r;s]$, $a=s$, and finally 
$T_0(\theta)={\alpha\theta\over \alpha+\varepsilon\theta}$ and
$T_1(\theta)={\alpha\theta\over \alpha+\theta}$. Since $T_0T_1=T_1T_0$, one gets ${T_0}^k{T_1}^j(\theta)={\alpha\theta\over \alpha+(j+\varepsilon k)\theta}$, in particular $(T_0T_1)^mT_0(s) = r$, so the space $\widetilde{\cB}$ introduced in $\cC(m)$ is 
$\widetilde{\cB} = L^r(\pi)$. Since 
$T_0(\theta)<\theta$ and $T_1(\theta)<\theta$, we have (0), and Lemma 4.2 gives (1) of $\cC(m)$. To prove (2), let $j\in\{1,\ldots,m\}$, and let 
$\theta\in I$ such that $\theta_j := T_1(T_0T_1)
^{j-1}(\theta)\in I$. We have 
$\theta_j< T_1^j(\theta)$, thus $\L^{T_1^j(\theta)}(\pi) \hookrightarrow \L^{\theta_j}(\pi)$, so the regularity property in (2) follows from the following lemma where $Q_k(t)$ stands for the kernel defined in Remark~(b) of Section 7.1. 

\noindent {\bf Lemma 7.4.} {\it 
Let $1\le j\leq m$. Then $Q(\cdot)\in\cC^j\left(\R,\L^{\theta}(\pi),\L^{T_1^j(\theta)}(\pi)\right)$ with $Q^{(k)} = Q_k\ $ ($k=0,\ldots,j$). } 

\noindent {\it Proof.} We denote $\|\cdot\|_p$ for 
$\|\cdot\|_{\L^p(\pi)}$, and 
$\|\cdot\|_{p,q}$ for $\|\cdot\|_{\L^p(\pi),\L^q(\pi)}$. Let us first show that $Q_k(\cdot)\in\cC^0\left(\R,\cB_\theta,\cB_{T_1^{j}(\theta)}\right)$ 
for any $k=0,\ldots,j$. 
The case $k=0$ follows from Lemma 4.2. For $1\leq k\leq j$, 
we have for $t_0,h\in \R$ and $f\in\L^\theta$, 
\begin{eqnarray*}
\Vert Q_k(t_0+h)f-Q_k(t_0)f\Vert_{{T_1}^j(\theta)} 
&=&\Vert Q_k(t_0+h)f-Q_k(t_0)f\Vert_{\scriptsize\alpha\theta\over \alpha +j\theta} \\
&\leq& 2\left\Vert \xi^{k} \min\{1,|h\xi|\} f\right\Vert
    _{\scriptsize\alpha\theta\over \alpha +j\theta} \\
&\leq& 2\left\Vert \xi^{k} \min\{1,|h\xi|\} \right\Vert
    _{\alpha\over j} 
   \Vert f\Vert_\theta
\end{eqnarray*}
with $\left\Vert \xi^{k} \min\{1,|h\xi|\} \right\Vert
    _{\alpha\over j} \r0$ when $h\r0$ by Lebesgue's theorem. 
Now let us prove $Q_k' = Q_{k+1}$ in $\cL\left(\cB_\theta,\cB_{{T_1}^j(\theta)}\right)$ 
for $k=0,\ldots,j-1$. 
Using $|e^{ia}-1-ia|\leq 2|a|\min\{1,|a|\}$, one gets for $t_0,h\in \R$ and 
$f\in\cB_\theta$~:
\begin{eqnarray*}
\left\Vert Q_k(t_0+h)f -Q_k(t_0)f - hQ_{k+1}(t_0)f\right\Vert_{
   \alpha\theta\over \alpha+j\theta} &\leq& 
 \left\Vert  Q\big(|\xi|^{k}\, \vert e^{ih\xi}-1-ih\xi\vert\, \vert f \vert\big)\right\Vert_{
   \alpha\theta\over \alpha+j\theta} \\
&\leq& 2 |h|\left\Vert |\xi|^{(k+1)} \min\{1,|h|\, |\xi|\} |f|\right\Vert
   _{\alpha\theta\over \alpha+j\theta}, 
\end{eqnarray*}
and the  previous computations yield $\|Q_k(t_0+h) -Q_k(t_0) - hQ_{k+1}(t_0)\|_{\theta
 ,{T_1}^j(\theta)} = o(|h|)$.  \fdem 

\noindent We know that $Q$ satisfies $(K)$ on $L^p(\pi)$ for every $p\in]1;+\infty[$ (Prop.~4.1).
Hence we have (3) of $\cC(m)$, and (4) is obvious from the definition of $T_0$, $T_1$ and $r$. \fdem

\noindent In this example, one can also use Lemma 4.2 and Lemma 7.4
to prove that Hypothesis $\cC'(m)$ is satisfied 
by taking $A>m$ such that $r={\alpha s\over\alpha+As}$ and by setting
$\cB_\theta:=L^{ \frac{\alpha s}
{\alpha+\theta s} }(\pi)$. 
%=================================================================================
%\noindent{\bf 8. A one-dimensional first-order Edgeworth expansion.} 
\section{A one-dimensional first-order Edgeworth expansion}
\noindent In this section we assume that $d=1$ (i.e.~$\xi$ is a real-valued measurable function on $E$). When $(X_n)_n$ is Harris recurrent,  the regenerative method provides 
Edgeworth expansions under some ``block'' moment conditions \cite{malin} \cite{jensen}. 
Here we do not assume Harris recurrence, and we present 
an alternative statement. 
To that effect, we shall appeal to Hypothesis $\cC(3)$ of Section 7.1 which ensures (Corollary 7.2) that 
the dominating eigenvalue $\lambda(t)$ of $Q(t)$ is three times continuously differentiable: 
then one shall be able to repeat the arguments of the i.i.d.~first-order Edgeworth expansion of \cite{fel} (Th.~1 p.~506). \\[0.1cm]
\noindent We denote by $\eta$ the density function of $\cN(0,1)$ and 
by $\cal N$ its distribution function.  
The next theorem extends the first-order Edgeworth expansion of the i.i.d.~case, 
with an additional asymptotic bias, namely $b_\mu = \lim_n\E_\mu[S_n]$ which depends on the initial distribution $\mu$. 
As for i.i.d.r.v., this bias is zero in the stationary case (i.e.~$b_\pi=0$). 

\noindent {\bf Theorem 8.1.} {\it 
Suppose that $\pi(|\xi|^3) < +\infty$, that Hypothesis $\cC(3)$ of Section 7.1 holds with $\cB\hookrightarrow\widetilde{\cB}\hookrightarrow\L^1(\pi)$, that the non-arithmeticity condition (S) of Section 5.1 holds on $\cB$,  
and finally that the initial distribution $\mu$ is in $\widetilde{\cB}'$. 
Then the real numbers \\[0.15cm]
\indent  $\displaystyle \ \ \sigma^2 = \lim_n\frac{1}{n}\, \E_\mu[S_n^2] = 
\lim_n\frac{1}{n}\, \E_\pi[S_n^2],\ \ \ \ \ 
m_3 = \lim_n\frac{1}{n}\, \E_\pi[S_n^3],\ \ \ \ \ b_\mu = \lim_n\E_\mu[S_n],$ \\[0.15cm]
are well-defined, and if $\sigma>0$, the following expansion holds uniformly in $u\in\R$ }
$$(E)\ \ \ \P_\mu\left(\frac{S_n}{\sigma\sqrt n}
\leq u\right) = \cN(u) + \frac{m_3}{6\sigma^3\sqrt n} (1-u^2)\, \eta(u) - 
\frac{b_\mu}{\sigma\sqrt n}\, \eta(u) + 
o(\frac{1}{\sqrt n}).$$
\noindent It will be seen in the proof of Lemma~8.4 below that \\[0.15cm]
\indent $\displaystyle \  \ \ \ \ \ \  \ \ \  \ \ \  \ \ \ \ \ \  
\big|\sigma^2 - \frac{1}{n}\, \E_\mu[S_n^2]\, \big|= O\big(\frac{1}{n}\big)\ \ $ and 
$\displaystyle\ \ \big|m_3 - \frac{1}{n}\, \E_\pi[S_n^3]\, \big| =  O\big(\frac{1}{n}\big)$.  

\noindent {\it Case of the strongly ergodic Markov chains on $\L^2(\pi)$.} \\[0.1cm]
In the special case of uniform ergodicity, the expansion (E) was established in  \cite{nag2} for any initial distribution, 
under some hypothesis on the absolute continuous component of $Q(x,dy)$ w.r.t.~$\pi$ and under the following restrictive 
operator-moment condition:   
there exists $g : \R\r\R$ such that $g(u)\r +\infty$ when $|u|\r+\infty$ and 
$\sup_{x\in E} \int_E|\xi(y)|^3\, g(|\xi(y)|)\, Q(x,dy) < +\infty$. 
In \cite{datta}, this result is slightly improved, more precisely (E) is established under the weaker (but still restrictive) 
moment condition $\sup_{x\in E} \int_E|\xi(y)|^3\, Q(x,dy) < +\infty$ and under some refinements of the nonlattice condition given in \cite{nag2}. In the stationary case (i.e.~under $\P_\pi$), the general asymptotic expansions established in \cite{gotze} apply to the uniformly ergodic Markov chains: they yield (E) when $\pi(|\xi|^4) < +\infty$, but under the so-called Cram\'er condition that is much stronger than the nonlattice one. \\
%the r.v.~$\xi(X_0)$ satisfies the so-called Cram\'er's condition. 
With the help of Theorem 8.1, one obtains here the following improvement which is moreover valid for the more general context of the strong ergodicity on $\L^2(\pi)$.  

\noindent{\bf Corollary 8.2.} {\it 
Let us suppose that $(X_n)_{n\geq0}$ is a strongly ergodic Markov chain on $\L^2(\pi)$, that $\pi(|\xi|^\alpha) < +\infty$ with some $\alpha>3$, and that $\xi$ is nonlattice (Prop.~5.4). Then we have (E) for any initial distribution of the form $d\mu=\phi\, d\pi$, where $\phi\in\L^{r'}(\pi)$ for some $r'>\frac{\alpha}{\alpha-3}$.}

\noindent{\it Proof.} Let $r'$ be fixed as above and let $r$ be such that $\frac{1}{r} + \frac{1}{r'} = 1$. Then $1<r<\frac{\alpha}{3}$, and since $\frac{\alpha s}{\alpha+3s}\nearrow\frac{\alpha}{3}$ when $s\r+\infty$, on can choose 
$s$ such that $s>{\alpha \over \alpha-3}$ and $\frac{\alpha s}{\alpha+3s} > r$. We have $\cC(3)$ with $\cB = \L^s(\pi)$, $\widetilde{\cB}=\L^r(\pi)$ (Prop.~7.3), and (S) on $\L^s(\pi)$ (see the proof of Cor.~5.5). \fdem

\noindent{\it Proof of Theorem 8.1.}  We shall appeal repeatedly to the notations and the 
conclusions of Theorem (K-L) (cf.~Section 4) 
and of Corollary 7.2 (case $m=3$). The existence of $\sigma^2$, $m_3$ and $b_\mu$ follows from the two next lemmas. 

\noindent {\bf Lemma 8.3.} {\it We have $\lambda'(0) = 0$ and 
$\mu(\Pi'(0)1_E) = i\, \sum_{k\geq1}\mu(Q^k\xi) = i\,  \lim_n\E_\mu[S_n]$.  }

\noindent {\it Proof.} By deriving the equality $Q(\cdot)\Pi(\cdot)1_E = \lambda(\cdot)\, \Pi(\cdot)1_E$, one gets \\[0.15cm] 
\indent $\displaystyle \ \ \ \ \ \ \ \ \ \ \ \ \ \ \ \ \ \ \ \ \ 
Q'(0)1_E + Q\Pi'(0)1_E = \lambda'(0)\, 1_E + \Pi'(0)1_E\ $ in  $\widetilde{\cB}$. \\[0.15cm]
Thus $\pi(Q'(0)1_E) + \pi(\Pi'(0)1_E) = \lambda'(0) + \pi(\Pi'(0)1_E)$. This 
gives $\lambda'(0) =  i\, \pi(Q\xi) =  i\, \pi(\xi)= 0$, and $iQ(\xi) + Q\Pi'(0)1_E = \Pi'(0)1_E$ in $\widetilde{\cB}$. 
Therefore we have $\Pi'(0)1_E-\pi(\Pi'(0)1_E) 
= i\, \sum_{k\geq1}Q^k\xi$.
This series is absolutely convergent in $\widetilde{\cB}$ since $\pi(Q\xi)=0$, $Q\xi = -iQ'(0)1_E
\in\widetilde{\cB}$ and 
$Q$ is strongly ergodic on $\widetilde{\cB}$. 
Moreover, we have $\pi(\Pi'(0)1_E)=0$.
Indeed, by deriving $\Pi(t)^2=\Pi(t)$,
we get $2\pi(\Pi'(0)1_E)=\pi(\Pi(0)\Pi'(0)1_E+\Pi'(0)\Pi(0)1_E)
  =\pi(\Pi'(0)1_E)$.
Since $\mu\in\widetilde{\cB}'$, this yields the first equality of the second assertion. The second one is obvious.  \fdem

\noindent {\bf Lemma 8.4.} {\it We have $\lim_n\frac{1}{n}\, \E_\mu[S_n^2] = - \lambda''(0)\ $ and 
$\ \lim_n\frac{1}{n}\E_\pi[S_n^3] = i\lambda^{(3)}(0)$. }

\noindent {\it Proof of Lemma 8.4.} For convenience, let us assume that $\mu=\pi$ and prove the two equalities of Lemma 8.4 at once (see Rk.~below). 
Since $\E_\pi[\, |\xi(X_k)|^3\, ] = \pi(|\xi|^3) < +\infty$, we have 
$\E_\pi[\, |S_n^3|\, ]< +\infty$, so \\[0.15cm]
\indent  $\displaystyle \ \ \ \ \ \ \ \ \ \ \ \ \ \ 
\E_\pi[e^{itS_n}] = 1 - \E_\pi[S_n^2]\, \frac{t^2}{2} - i\,\E_\pi[S_n^3]\, \frac{t^3}{6} + o_n(t^3).$  \\[0.15cm]
Besides, Formula (CF) (cf.~Section 3) and the equality $ Q(t)^n =  \lambda(t)^n\Pi(t) +  N(t)^n$ (see Section~7.2) give \\[0.12cm]
\indent  $\displaystyle \ \ \ \ \ \ \ \ \ \ \ \ \ \ \ \ \ \ \ \ \ \ \ \ \ \ \ \ 
\E_\pi[e^{itS_n}] = \lambda(t)^n\, \pi(\Pi(t)1_E) + \pi(N(t)^n1_E),$\\[0.15cm]
and, since $\lambda'(0) = 0$ and $\pi(\Pi'(0)1_E) = 0$ (Lemma 8.3), it follows from Hypothesis  $\cC(3)$ and Corollary 7.2 that \\[0.1cm]
\indent  $\displaystyle \ \ \ \ \ \ \ \ \ 
\lambda(t)^n =  1+ n\frac{\lambda''(0)}{2}t^2 + n\frac{\lambda^{(3)}(0)}{6}t^3 + o_n(t^3), \ \ \ 
\pi(\Pi(t)1_E) = 1 + ct^2 + dt^3 + o(t^3),$\\[0.15cm]
with some $c,d\in\C$, and since $N(0)1_E=0$, we have $\pi(N(t)^n1_E) = e_nt+f_nt^2+g_nt^3+o_n(t^3)$ for all $n\geq 1$, 
with some $e_n, f_n, g_n\in\C$. Moreover, from Assertion (ii) in Corollary 7.2, it follows that the sequences 
$(e_n)_n$, $(f_n)_n$ and $(g_n)_n$ are bounded. From the previous expansions, one can write 
another third order Taylor expansion for $\E_\pi[e^{itS_n}]$, 
from which we easily deduce the following equalities (and so Lemma 8.4): \\[0.15cm]
\indent  $\ \ \ \ \ \ \ \ \ n\lambda''(0)+2c+2f_n = -\E_\pi[S_n^2]\ \ \ \ \mbox{and}\ \ \ \ \ 
n\lambda^{(3)}(0) + 6d + 6g_n =  - i\, \E_\pi[S_n^3].$ \fdem 

\noindent {\it Remark.} By using the above arguments with second-order Taylor expansions, it can be easily proved that the 
first equality of Lemma 8.4 is valid under 
Hypothesis $\cC(2)$ for any 
$\mu\in\widetilde{\cB}'$. To prove $\E_\mu[S_n^2]< 
+\infty$ under Hypothesis $\cC(2)$ and for 
$\mu\in\widetilde{\cB}'$, 
we notice that $Q''(0)1_E = -Q(\xi^2)\in\widetilde{\cB}$, 
so $Q^k(\xi^2)\in\widetilde{\cB}$ for $k\geq1$, and 
$\E_\mu[\xi(X_k)^2] = \mu(Q^k\xi^2)< +\infty$. 

\noindent The proof of the Edgeworth expansion (E) 
is close to that of the i.i.d.~case \cite{fel} 
(XVI.4). 
For convenience, one may assume, without any 
loss of generality, that $\sigma=1$ (of course this reduction also leads to alter the constants 
$m_3$ and $b_\mu$). Set 
$$G_n(u) = \cN(u) + \frac{m_3}{6\sqrt n} (1-u^2)\, \eta(u) - \frac{b_\mu}{\sqrt n} \eta(u)\ \ \ (u\in\R).$$
Then $G_n$ has a bounded derivative $g_n$ on $\R$ whose Fourier transform $\gamma_n$ is given by 
$$\gamma_n(t) =  \gamma_{0,n}(t) + \gamma_{\mu,n}(t),\ \ \mbox{where}\ \ 
\gamma_{0,n}(t) = e^{-\frac{1}{2}t^2}\bigg(1 + \frac{m_3}{6\sqrt n}(it)^3\bigg)
\ \ \ \mbox{and}\ \ \ \gamma_{\mu,n}(t) =  e^{-\frac{1}{2}t^2}\bigg(i\, \frac{b_\mu}{\sqrt n}\, t\bigg).$$
Let us notice that the part $\gamma_{0,n}(t)$ has the same form as in the i.i.d.~context. Let us set 
$$\forall t\in \R,\ \ \phi_n(t) = \E_\mu[e^{itS_n}].$$ 
The first question is to prove the so-called Berry-Esseen inequality 
$$\sup_{u\in\R}\bigg|\, 
\P_\mu(\frac{S_n}{\sqrt n}\leq u) - G_n(u)\bigg| 
\leq \frac{1}{\pi} \int_{-T}^T \bigg|\, 
\frac{\phi_n(\frac{t}{\sqrt n}) - \gamma_n(t)}{t}\, \bigg|\, dt\ +\ \frac{24m}{\pi T},$$
where $m=\sup\{|G'_n(u)|,\, n\geq1,\, u\in\R\}$. 
To do this, let us observe that all the hypotheses 
of Lemma 2 in Section XVI.3 of \cite{fel}, which 
provides this inequality, are satisfied, 
except  $\gamma'_n(0)=0$ because of the additional 
term $\gamma_{\mu,n}(t)$ in  $\gamma_n(t)$. 
However it can be easily seen that  the above cited 
lemma of  \cite{fel} still holds under the condition 
that 
$\frac{\gamma_n(t)-1}{t}$ is continuous at 
the origin. Indeed the argument  in \cite{fel}  
(p.~511)  deriving from the Riemann-Lebesgue theorem 
then remains valid. Obviously the previous condition 
on $\gamma_n$ is fulfilled since 
$\frac{\gamma_{\mu,n}(t)}{t} =  i\, \frac{b_\mu}{\sqrt n}\, e^{-\frac{1}{2}t^2}$. 
Thus we have the desired Berry-Esseen inequality and we can now proceed as in \cite{fel}: 
let $\varepsilon>0$,  let $T=a\sqrt n$ with $a$ such that $\frac{24m}{\pi a} < \varepsilon$. So 
$\frac{24m}{\pi T}\leq  \frac{\varepsilon}{\sqrt n}$. \\ 
Let $0<\delta<a$ such that 
$[-\delta,\delta]$ is contained in the interval ${\cal O}$ of Theorem (K-L) applied on $\cB$, and let us write 
$$\int_{-a\sqrt n}^{a\sqrt n} \bigg|\, 
\frac{\phi_n(\frac{t}{\sqrt n}) - \gamma_n(t)}{t}\, \bigg|\, dt = 
\int_{\delta\sqrt n\leq |t|\leq a\sqrt n}\ + \ \int_{|t| \leq \delta\sqrt n}\ \ := A_n + B_n.$$
The property (E) then follows from the two next lemmas. \fdem 

\noindent {\bf Lemma 8.5.} {\it There exists  $N_0\in\N^*$ such that $A_n \leq \frac{\varepsilon}{\sqrt n}$ for all $n\geq N_0$. }

\noindent {\it Proof.} From Formula (CF) (cf.~Section 3), 
Condition (S) (cf.~Section 5.1) on $\cB$ applied with 
$K_0 = [-a,-\delta]\cup [\delta,a]$, and from 
$\mu\in\widetilde{\cB}'\subset\cB'$, 
there exist $\rho<1$ and $c'\geq 0$ such that we have, for $n\geq1$ and $u\in K_0$: 
$|\phi_n(u)| = |\mu(Q(u)^n1_E)| \leq c'\, \rho^n$. So 
$$\int_{\delta\sqrt n\leq |t|\leq a\sqrt n} 
\frac{|\phi_n(\frac{t}{\sqrt n})|}{|t|}\, dt = 
\int_{\delta\leq |u|\leq a}
 \frac{|\phi_n(u)|}{|u|}\, du \leq 
\frac{2a}{\delta}\, c'\, \rho^n.$$
Moreover, for 
$n$ sufficiently large,  we have 
$\int_{\delta\sqrt n\leq |t|\leq a\sqrt n} 
\frac{|\gamma_n(t)|}{|t|}\, dt \leq  
\int_{|t|\geq\delta\sqrt n } |\gamma_n(t)|dt$. We easily deduce Lemma 8.5 from the two last estimates. \fdem

\noindent {\bf Lemma 8.6.} {\it There exists  $N_0'\in\N^*$ such that $B_n \leq \frac{\varepsilon}{\sqrt n}$ for all $n\geq N_0'$. }

\noindent {\it Proof.} Using $\gamma_n(t) = 
\gamma_{0,n}(t) +\gamma_{\mu,n}(t)$ and the equality 
$\phi_n(t)  = \lambda(t)^n\, \mu(\Pi(t)1_E) + \mu(N(t)^n1_E)$ 
which follows from (CF) and from Theorem (K-L), one can write for any $t$ such that $|t| \leq \delta\sqrt n$ 
\begin{eqnarray*}
\phi_n(\frac{t}{\sqrt n}) - \gamma_n(t) &=& 
\bigg(\lambda(\frac{t}{\sqrt n})^n - \gamma_{0,n}(t)\bigg)\ +\ 
\lambda(\frac{t}{\sqrt n})^n\, \bigg(\mu(\Pi(\frac{t}{\sqrt n})1_E) - 1  - i\, b_\mu\, \frac{t}{\sqrt n}\bigg) \\ 
&\ & \ \ \ \ \ \ + \ \ i\, b_\mu\, \frac{t}{\sqrt n} \bigg(\lambda(\frac{t}{\sqrt n})^n - e^{-\frac{1}{2}t^2}\bigg) \ +\   
\mu(N(\frac{t}{\sqrt n})^n1_E)\\
&:=& i_n(t) + j_n(t) + k_n(t) + \ell_n(t). 
\end{eqnarray*}
Therefore: $\displaystyle B_n \leq \int_{|t| \leq \delta\sqrt n} \bigg(|i_n(t)| + |j_n(t)| + |k_n(t)| + |\ell_n(t)|\bigg)\, 
\frac{dt}{|t|} :=  I_n + J_n + K_n + L_n.$\\[0.18cm]
Then Lemma 8.6 follows from the assertions (i)-(l) below for which, as in the i.i.d.~case, we shall repeatedly appeal to the 
following remark: 
using the Taylor expansion $\lambda(t) = 1-\frac{t^2}{2}+o(t^2)$ near 0 (use Lemmas 8.3-4 and $\sigma^2=1$), 
one can choose the real number $\delta$ such that 
$|\lambda(u)| \leq  1-\frac{u^2}{4} \leq e^{-\frac{u^2}{4}}$ when $|u|\leq \delta$, hence we have 
$\displaystyle|\lambda(\frac{t}{\sqrt n})|^n \leq e^{-\frac{t^2}{4}}$ for any $|t| \leq \delta\sqrt n$. \\[0.18cm]
{\bf (i)} $\exists N_1\in\N^*, \forall n\geq N_1,\ I_{n} \leq  \frac{\varepsilon}{\sqrt n}$. 
This can be proved exactly as in the i.i.d.~case \cite{fel} since we have 
$\lambda(t) =  1 - \frac{t^2}{2} - i\, \frac{m_3}{6}t^3 + o(t^3)$ (Lemmas 8.3-4). \\[0.18cm]
{\bf (j)} $\exists N_2\in\N^*, \forall n\geq N_2,\ J_{n} \leq 
\frac{\varepsilon}{\sqrt n}$. 
Indeed, since $u\mapsto\mu(\Pi(u)1_E)$ has two continuous derivatives on $[-\delta,\delta]$ (Coro.~7.2) 
and $\mu(\Pi'(0)1_E) = i b_\mu$ (Lemma 8.3), there exists $C>0$ such that: 
$\displaystyle J_{n} \leq \int_{|t| \leq \delta\sqrt n} e^{-\frac{t^2}{4}}\, 
\frac{Ct^2}{n}\, \frac{dt}{|t|} \leq \frac{C}{n}\, \int_{-\infty}^{+\infty} e^{-\frac{t^2}{4}}\, |t|\, dt$.  \\[0.18cm]
{\bf (k)} $\exists N_3\in\N^*, \forall n\geq N_3,\ K_n \leq  \frac{\varepsilon}{\sqrt n}$. Indeed we have 
$\displaystyle K_{n} \leq  \frac{|b_\mu|}{\sqrt n}\, 
\int_{|t| \leq \delta\sqrt n}\bigg|\lambda(\frac{t}{\sqrt n})^n - e^{-\frac{1}{2}t^2}\bigg|\, dt$, and from the already mentioned 
second order Taylor expansion of $\lambda(t)$ and  Lebesgue's theorem, it follows that this last integral converges to 0 when 
$n\r+\infty$. \\[0.18cm]
{\bf (l)} $\exists N_4\in\N^*, \forall n\geq N_4,\ L_n \leq  \frac{\varepsilon}{\sqrt n}$.
Indeed, the function $\chi_n : u\mapsto \mu(N(u)^n1_E)$ is continuously 
differentiable on  $[-\delta,\delta]$ and 
there exists $C'>0$ such that we have for all $n\geq1$ and $u\in[-\delta,\delta]$: $|\chi_n'(u)| \leq C'\tilde\kappa^n $ 
(Corollary 7.2(ii)). Since $N(0)1_E=0$, one then obtains $|\mu(N(u)^n1_E)| \leq C'\, \tilde\kappa^n\, |u|$ for $|u|\leq\delta$, so  
$L_n \leq  \frac{C'}{\sqrt n}\, \, \tilde\kappa^n\, 2\delta\sqrt n = 2C'\delta\, \tilde\kappa^n = o( \frac{1}{\sqrt n})$. \fdem 

\noindent {\bf Remark.} In the i.i.d.~case, higher-order Edgeworth expansions can be established, see \cite{fel} (Th.~2 p.~508), but the non-arithmeticity assumption has to be replaced with the so-called more restrictive Cr\`amer condition. 
Notice that, in our context, this condition can be extended to some operator-type Cr\`amer condition, and that the present method could be then employed to prove similar higher-order  Edgeworth expansions. However, the main difficulty is to reduce this operator-type Cr\`amer assumption to some more practical condition. 
%=============================================================
%\noindent{\bf 9. A multidimensional Berry-Esseen theorem. }
\section{A multidimensional Berry-Esseen theorem}
We want to estimate the rate of convergence in the 
central limit theorem
for a ${\R}^d$-valued function $\xi = (\xi_1,\ldots,\xi_d)$. A natural way 
to do this is in
the sense of the Prohorov metric. Let us recall 
the definition of this metric
and some well-known facts about it.
We denote by ${\cal B}({\R}^d)$ the Borel
$\sigma$-algebra of ${\R}^d$ and by
${\cal M}_1({\R}^d)$ the set of probability
measures on $({\R}^d,{\cal B}({\R}^d))$.

\noindent {\bf The Prohorov metric} \cite{Billingsley,Dudley}. {\it For all $P,Q$ in ${\cal M}_1({\R}^d)$, we define: 
$${\cal P}(P,Q):=\inf\left\{\varepsilon>0\ :\ 
      \forall B\in{\cal B}({\R}^d),\ 
(P(B)-Q(B^\varepsilon))\le\varepsilon\right\},$$
where $B^\varepsilon$ is the open  
$\varepsilon$-neighbourhood of $B$.}

\noindent  {\bf The Ky Fan metric for random variables.} {\it If $X$ and $Y$ are two ${\R}^d$-valued random variables
defined on the same probability space $(E_0,{\cal T}_0,{P}_0)$, 
we define~:
$${\cal K}(X,Y):=\inf\left\{\varepsilon>0\ :\ 
  {P}_0\left(\vert X-Y\vert_2>\varepsilon\right)
   <\varepsilon\right\} .$$}
\noindent Let us recall that $\lim_{n\rightarrow +\infty}{\cal K}(X_n,Y)=0$
means that $(X_n)_n$ converges in probability to $Y$.

\noindent {\bf Proposition} (\cite{Dudley} Corollary 11.6.4).   
{\it For all $P,Q$ in ${\cal M}_1({\R}^d)$, the quantity
${\cal P}(P,Q)$ is the infimum of ${\cal K}(X,Y)$ 
over the couples $(X,Y)$ of 
${\R}^d$-valued random variables
defined on the same probability space, whose distributions are respectively $P$ and $Q$.}

\noindent For any $n\ge 1$, $\mu_*\left({S_n\over\sqrt{n}}\right)$ stands for the law of ${S_n\over\sqrt{n}}$ under $\P_\mu$, and 
we denote by $S_n^{\otimes 2}$ the random variable 
with values in the set of
$d\times d$ matrices given by: \\
\indent $\displaystyle \ \ \ \ \ \ \ \ \ \ \ \ \ \ \ \ \ \ \ \ \ \ \ \ \ \ \ \ \ \ \ \ \ \ \ 
\left({S_n}^{\otimes 2}\right)_{i,j} = \sum_{k,\ell=1}^n\xi_i(X_k)\xi_j(X_\ell)$. \\[0.15cm]
\noindent {\bf Theorem 9.1.}
{\it Let us fix 
$m:=max\left(3,\lfloor{d/ 2}\rfloor+1\right)$.
Suppose that Hypothesis ${\cal C}(m)$ (of Section~7.1) holds with $\cB\hookrightarrow\widetilde{\cB}\hookrightarrow\L^1(\pi)$, and that $\mu\in\widetilde{\cB}'$. Then the following limits exist and are equal: \\[0.14cm]
\indent $\ \ \ \ \ \ \ \ \ \ \ \ \ \ \ \  \ \ \ \ \ \ \ \ \ \ 
\displaystyle \Gamma:=\lim_{n\rightarrow +\infty}\, {1\over n}
   \E_\pi[{S_n}^{\otimes 2}] =
   \lim_{n\rightarrow +\infty}\, {1\over n}
    \E_\mu[{S_n}^{\otimes 2}]$.\\[0.14cm]
If $\Gamma$ is invertible, then $\left({S_n\over\sqrt{n}}\right)_n$ converges in 
distribution under $\P_\mu$ to the gaussian 
distribution ${\mathcal N}(0,\Gamma)$, and we have }
$$ {\cal P}\left(\mu_*\left({S_n\over\sqrt{n}}
    \right),{\mathcal N}(0,\Gamma)\right)
  =O(n^{-1/2}).$$
\noindent In the i.i.d.~case, thanks to a smoothing 
inequality (see Proposition 9.3) and to an 
additional truncation argument, 
the conclusion of Theorem 9.1 holds if the random variables admit a moment of 
order $3$. 
For the strongly ergodic Markov chains on $\L^2(\pi)$, one gets the 
following statement which is a consequence of Theorem 
9.1 and 
of Proposition 7.3 (proceed as for Corollary 8.2).  

\noindent{\bf Corollary 9.2.} {\it Let us suppose that $(X_n)_{n\ge 0}$ is a strongly ergodic Markov chain on $\L^2(\pi)$, that $\pi(|\xi|_2^\alpha)<+\infty$
for some $\alpha>m:=\max(3,\lfloor d/2\rfloor+1)$,
and that the initial distribution satisfies
$d\mu=\phi d\pi$ with $\phi\in \L^{r'}(\pi)$
for some $r'>{\alpha\over\alpha-m}$.
Then the conclusion of Theorem 9.1 is true. }

\noindent Concerning the special case of the uniform ergodicity, notice that \cite{gotze} provides 
a multidimensional uniform Berry-Esseen type estimate when $\pi(|\xi|_2^4)<+\infty$. 
However, the hypothesis $\mu=\pi$ (i.e.~$(X_n)_n$ is stationary), and the Cram\'er condition for $\xi(X_0)$, are required in \cite{gotze}, while the (Prohorov) estimate in Corollary 9.2, and more generally in 
Theorem 9.1, is valid in the non-stationary case and without any lattice-type condition. 

\noindent Let us mention that Theorem 9.1 remains true when $\Gamma$ is non invertible 
if, for every $\beta\in{\R}^d$ such that $\langle \beta,\Gamma\beta\rangle=0$, we are able
to prove that $\sup_n\Vert\langle \beta,S_n\rangle\Vert_\infty<+\infty$. In this case, up to a linear
change of coordinates and to a possible change of $d$, we are led to the invertible case
(see Section 2.4.2 of \cite{FPAAP}). This remark applies to the  Knudsen gas model (see Section 1).

\noindent When $d=1$, Theorem 9.1 gives the uniform Berry-Esseen result under Condition  ${\cal C}(3)$ if 
the asymptotic variance $\sigma^2$ is 
nonzero. This is an easy consequence of the definition  of ${\cal P}$ by taking $B=(-\infty,x]$ and $B=(x,+\infty)$. However, as already 
mentioned, ${\cal C}(3)$ is in practice a little more restrictive than the conditions of Section 6; for instance, compare 
the expected condition $\pi(|\xi|^3)<+\infty$ of Corollary 6.3 with that of Corollary 9.2 (case $d=1$). 

The proof of Theorem 9.1 is based on Corollary 7.2,
on lemmas 8.3 and 8.4 
and on the following smoothing inequality due to Yurinskii
\cite{Yurinskii}~:

\noindent {\bf Proposition 9.3.}
{\it 
Let $Q$ be some non degenerate $d$-dimensional normal distribution.
There exists a real number $c_0>0$ such that,
for any real number $T>0$ and 
for any Borel probability measure $P$ admitting moments
of order $\lfloor {d\over 2}
          \rfloor+1$, we have: }
$${\cal P}\left(P,Q\right)\le c_0
\left[{1 \over T}
  +\left(\int_{|t|_2<T}
    \sum_{k=0}^{\lfloor {d\over 2}
          \rfloor+1}
   \sum_{\{i_1,...,i_k\}\in\{1,...,d\}^k}
    \left\vert
    {\partial^k\over \partial t_{i_1}...\partial t_{i_k}}
     \left( P(e^{i\langle t,\cdot\rangle})
         -Q(e^{i\langle t,\cdot\rangle})\right)
      \right\vert^2\, dt\right)^{1\over 2}\right].$$
\noindent {\it Proof of Theorem 9.1.}
The proof uses Corollary 7.2 which is applied here 
under Hypothesis ${\cal C}(m)$ with $m$ defined in Theorem 9.1. In particular we have 
$m\geq3$, and we shall use repeatedly the fact 
that $1_E\in\cB$ and $\pi,\mu\in\widetilde{\cB}'$. Since the proof has common points with 
the proof given in Section 2.4.1 of \cite{FPAAP}, we do not give all
the details. We shall refer to \cite{FPAAP}
for some technical points. \\
\noindent The existence of the asymptotic covariance matrix $\Gamma$ as defined in Theorem 9.1 follows from the next lemma in which 
$\nabla$ and $Hess$ denote the gradient and  the Hessian matrix.   

\noindent {\bf Lemma 9.4.} {\it We have
$\nabla \lambda(0) = 0$ and  
$\lim_n\frac{1}{n}\, \E_\mu[S_n^{\otimes 2}] = 
- Hess\lambda(0)$. }

\noindent {\it Proof.} These properties have 
been proved in the case $d=1$ (Lemmas 8.3-4).
We deduce from them the multidimensional version
by considering, for any $\alpha\in{\R}^d$, the function $t\mapsto Q(t\alpha)$ defined on $\R$. \fdem  

\noindent Without any loss of generality, up to a 
linear change of variables, we may suppose
that the covariance matrix $\Gamma$ is the 
identity matrix.

\noindent Let $\beta>0$ be such that the closed ball  
$\{u\in\R^d : |u|_2 \le \beta \}$ is contained in the set 
$\cal O$ of Corollary 7.2. In the following, the couple $(t,n)$ ($t\in{\R}^d$, $n\ge 2$) will always satisfy 
the condition $|t|_2<{\beta\sqrt{n}}$. For such a couple, we have: ${ t\over\sqrt{n}}\in \cal O$. 
For any function $F$ defined on an open set of $\R^d$, $F^{(k)}$ 
will merely denote  any partial derivative of order $k$ of $F(\cdot)$. 

\noindent Set $\Xi_n(t) := \E_\mu[e^{i\langle t, 
\frac{S_n}{\sqrt n}\rangle}] - e^{-\frac{|t|_2^2}{2}}$. 
According to Proposition 9.3, 
it is enough to prove that we have for 
$k=0,\ldots,[\frac{d}{2}]+1$\\
\indent $\displaystyle \ \ (I) \ \ \ \ \  \ \ \  \ \ \ \ \ \  \  \ \ \  \ \ \ \ \ \ 
\left(\int_{|t|_2\leq{\beta\sqrt{n}}}|\Xi_n^{(k)}(t)|^2\, dt\right)^{1\over 2} = 
O\left({1\over\sqrt{n}}\right)$. \\[0.2cm] 
From the decomposition $\E_\mu[e^{i\langle u, S_n\rangle}] = \lambda(u)^n\, \mu(\Pi(u)1_E) + \mu(N(u)^n1_E)$ which is valid for 
$u\in{\cal O}$, it follows that  
\begin{eqnarray*}
\Xi_n^{(k)} &=&   \bigg(\lambda(\frac{\cdot}
{\sqrt n})^n - e^{-\frac{|\cdot|_2^2}{2}} \bigg)^{(k)} + 
\bigg\{\lambda(\frac{\cdot}{\sqrt n})^n\, \bigg(\mu(\Pi(\frac{\cdot}{\sqrt n})1_E) -1\bigg)\bigg\}^{(k)} + 
\bigg(\mu(N(\frac{\cdot}{\sqrt n})^n1_E)\bigg)^{(k)}\\
&:=& A_n^{(k)} + B_n^{(k)} + C_n^{(k)}
\end{eqnarray*}
where the functions $ A_n$, $B_n$ and $C_n$, defined 
on the set $\{t : |t|_2 < \beta\sqrt n\}$, are implicitly given by 
the above equality. 
In the sequel, we merely use the notation $F_n(t) = O(G_n(t))$ to express that $|F_n(t)| \leq C\, |G_n(t)|$ for some $C\in\R_+$ 
independent of $(t,n)$ such that $|t|_2\leq \beta\sqrt n$. \\[0.12cm]
Setting $N_n(\cdot) = N(\cdot)^n$, 
Corollary 7.2(ii) yields \\[0.12cm]
\indent $\displaystyle \ \ \  \ \ \  \ \ \  \ \ \  \ \ \ 
|C_n^{(k)}(t)| = n^{-\frac{k}{2}}\, |\mu(N_n^{(k)}(\frac{t}{\sqrt n})1_E)| =  O( n^{-\frac{k}{2}}\,\tilde\kappa^n)$. \\[0.12cm] 
So $\displaystyle \int_{|t|_2\leq\beta\sqrt n} |C_n^{(k)}(t)|^2\, dt = O(n^{\frac{d}{2}-k}\, \tilde\kappa^{2n}) 
= O(\frac{1}{\sqrt n})$. 
Now (I) will be clearly valid provided that we have, for some square Lebesgue-integrable function $\chi(\cdot)$ on $\R^d$: \\[0.1cm]
\indent $\displaystyle \ \ \  (II) \ \ \  \ \ \  \ \ \  \ \ \ \ \ \ \ \  \ \ \  \ \ \ \ \ \  \  \ 
 |A_n^{(k)}(t)|  + |B_n^{(k)}(t)| = O\bigg(\frac{1}{\sqrt n}\, \chi(t)\bigg)$. \\[0.1cm]
To prove this estimate for the term $A_n^{(k)}$, one 
can proceed as in the i.i.d.~case. Indeed, 
according to the previous lemma, 
the function $\lambda(\cdot)$ then satisfies the 
same properties and plays exactly the same role, 
as the common characteristic function of the i.i.d.~case
(see Section 3 of \cite{Yurinskii} and Lemma 8 of \cite{Rotar70} or \cite{FPAAP} pages 2349--2350).\\[0.15cm]
\noindent To study $B_n^{(k)}(t)$, set $\lambda_n(t) =  \lambda(\frac{t}{\sqrt n})^n$ for any $(t,n)$ such that 
$|t|_2 \leq \beta\sqrt n$ and,  for $|u|_2 
\leq \beta$, set $\alpha(u) = \mu(\Pi(u)1_E) -1$. 
With these notations, we have $B_n(t) = \lambda_n(t)\, \alpha(\frac{t}{\sqrt n})$, and 
any partial derivative $B_n^{(k)}(t)$ is a finite sum of terms of the form \\[0.15cm]
\indent $\displaystyle \ \ \  \ \ \  \ \ \  \ \ \ \ \  \ \ \  \ \ \  \ \ 
B_{n,p,q}^{(k)}(t) := \lambda_n^{(p)}(t)\, n^{-\frac{q}{2}}\, \alpha^{(q)}(\frac{t}{\sqrt n})\ $ with $\ p+q = k$.  \\[0.2cm]
\noindent{\bf Lemma 9.5.} {\it For $p = 0,\ldots,m$, we have $|\lambda_n^{(p)}(t)| = 
O\bigg((1+|t|_2^p)\,  e^{-\frac{|t|_2^2}{4}}\bigg)$. } 

\noindent Assume this lemma for the moment. 
Since we have, by Corollary 7.2(i), 
$\alpha(\frac{t}{\sqrt n}) = O(\frac{|t|_2}{\sqrt n})$ and $\alpha^{(q)}(\frac{t}{\sqrt n}) = O(1)$ for 
$1\leq q \leq m$,  
this lemma gives for $q=0$ \\[0.12cm]
\indent $\displaystyle \ \ \  \ \ \  B_{n,k,0}^{(k)}(t) = 
O\bigg((1+|t|_2^k)\,  e^{-\frac{|t|_2^2}{4}}\bigg)
\, O(\frac{|t|_2}{\sqrt n}) = 
O\bigg(\frac{1}{\sqrt n}\, (1+|t|_2^{k+1})\, 
 e^{-\frac{|t|_2^2}{4}}\bigg)$, \\[0.2cm]
\noindent and for $q\geq1$: $\displaystyle \ \  B_{n,p,q}^{(k)}(t) = 
O\bigg((1+|t|_2^p)\,  e^{-\frac{|t|_2^2}{4}}\bigg)
\, O(n^{-\frac{q}{2}}) = 
O\bigg(\frac{1}{\sqrt n}\, (1+|t|_2^{k+1})\, 
 e^{-\frac{|t|_2^2}{4}}\bigg)$. \\[0.12cm]
So all the $B_{n,p,q}^{(k)}(t)$'s are 
$O\big(\frac{1}{\sqrt n}\, \chi(t)\big)$ with 
$\chi(t) = (1+|t|_2^{k+1})\,  e^{-\frac{|t|_2^2}{4}}$, 
and this gives the estimate (II) for $B_n^{(k)}(t)$, and finally the proof of (I) is complete. 

\noindent{\it Proof of Lemma 9.5.} Recall $\Gamma$ is by hypothesis the identity matrix, so 
$\lambda(u) = 1 - \frac{|u|_2^2}{2} + o(|u|_2^2)$ 
as $u$ goes to $0$ (use Lemma 9.4).
Hence, for 
$|u|_2\leq \beta$ with $\beta$ possibly reduced, \\[0.12cm]
\indent $\displaystyle \ \ 
|\lambda(u)| \leq \left|\lambda(u) - 1 + 
\frac{|u|_2^2}{2}\right| + \left|1 - \frac{|u|_2^2}{2}\right| 
\leq \frac{|u|_2^2}{4} + 
(1 - \frac{|u|_2^2}{2}) \leq  1-\frac{|u|_2^2}{4} 
\leq e^{-\frac{|u|_2^2}{4}}$,\\[0.12cm]
so $|\lambda(\frac{t}{\sqrt n})| \leq 
e^{-\frac{|t|_2^2}{4n}}$ and 
$|\lambda(\frac{t}{\sqrt n})^n| \leq 
(e^{-\frac{|t|_2^2}{4n}})^{n} = 
e^{-\frac{|t|_2^2}{4}}$. This gives 
the estimate of the lemma for $p=0$. Now, in the case $p\geq1$, one can prove by a straightforward induction that 
$\lambda_n^{(p)}(t)$ is a finite sum of terms of the form \\[0.12cm]
\indent $\displaystyle \ \ \  \ \ \ \ \ 
\gamma(t,n) := 
n(n-1)\cdots(n-j+1)\, n^{-\frac{p}{2}}\, 
\lambda^{(s_1)}(\frac{t}{\sqrt n})\cdots\lambda^{(s_j)}(\frac{t}{\sqrt n})\, \lambda(\frac{t}{\sqrt n})^{n-j}$, \\[0.12cm]
with $j\in\{1,\ldots,p\}$, $s_i\geq1$, and $s_1+\cdots+s_j = p$ (for convenience, $j,s_1,\ldots,s_j$ have been neglected 
in the above notation $\gamma(t,n)$).  
So we must prove that, given such fixed 
$j,s_1,\ldots,s_j$, we have 
$\gamma(t,n) = O\left((1+|t|_2^p)\,  
e^{-\frac{|t|_2^2}{4}}\right)$. 
To that effect, let us observe that 
$\lambda^{(1)}(\frac{t}{\sqrt n})= O
(\frac{|t|_2}{\sqrt n})$ 
since $\lambda^{(1)}(0) = 0$, and 
that $\lambda^{(s)}(\frac{t}{\sqrt n})= O(1)$ for any $s=2,\ldots,m$. This leads to define $a = Card\{i : s_i=1\}$. Then we have   
\begin{eqnarray*}
\gamma(t,n) = O\bigg(n^{j-\frac{p}{2}}\, 
  \frac{|t|_2^a}{n^{\frac{a}{2}}}\, 
   (e^{-\frac{|t|_2^2}{4n}})^{n-j}\bigg) &=& 
O\bigg(e^{\frac{j}{4}\left|\frac{t}
    {\sqrt n}\right|_2^2} n^{j-\frac{p}{2} 
    -\frac{a}{2}}\, |t|_2^a \,
     e^{-\frac{|t|_2^2}{4}}\bigg) \\
&=&  O\bigg(n^{\frac{1}{2}(2j-p-a)}\, 
    (1+|t|_2^p) \, e^{-\frac{|t|_2^2}{4}}\bigg). 
\end{eqnarray*}
For the last estimate, we used the fact that 
$\left|\frac{t}{\sqrt n}\right|_2 \leq \beta$ and $a \leq p$. Finally observe that we have 
$p = s_1+\cdots+s_j \geq a+2(j-a)$ by definition of the number $a$, thus $2j-p-a \leq 0$, so that the desired estimate on 
$\gamma(t,n)$ follows from the previous one. \fdem 
%============================================================================================
%\noindent{\bf 10. Application to $v$-geometrically ergodic Markov chains.} 
\section{Application to $v$-geometrically ergodic Markov chains} 
For the moment, the abstract results 
of the previous sections have been only applied to the (somewhat restrictive) strongly ergodic Markov chains on $\L^2(\pi)$. This section and the next one present applications to other practicable Markov models, namely the so-called $v$-geometrically ergodic Markov chains and the random iterative models (see Examples~2-3 in Section~1). The interest of 
these models for statistical applications and for stochastic algorithms is fully described in \cite{mey} \cite{duf}, and of course, 
the rate of convergence in the c.l.t.~and the Edgeworth expansions are of great importance in practice, see e.g \cite{mira} 
\cite{bertail}. For these models, all the previously studied limit theorems will be stated under general and simple moment conditions. 

\noindent  Throughout this section, we suppose that the $\sigma$-field $\cE$ is countably generated, that $(X_n)_{n\geq0}$ is 
aperiodic and $\psi$-irreducible w.r.t.~a certain positive $\sigma$-finite measure $\psi$ on  $E$. 

\noindent Moreover, given an unbounded function $v : E\r [1,+\infty[$, we assume 
that $(X_n)_{n\geq0}$ is $v$-geometrically ergodic, that is $\pi(v) < +\infty$ and 
there exist real numbers $\kappa_0<1$ and $C\geq 0$ such that we have, for all $n\geq1$ and $x\in E$, \\[0.1cm]
\indent $\displaystyle \ \ \ \ \ \ 
\sup\bigg\{\left|Q^nf(x)-\pi(f)\right|,\ f\, :\, E\r \C \ measurable,\ |f|\leq v\ \bigg\}
\leq C\, \kappa_0^n\, v(x)$. \\[0.1cm]
\noindent If $w$ is an unbounded function defined on  $E$ and taking values in $[1,+\infty[$, we denote by
$(\cB_w,\|\cdot\|_w)$ the  weighted supremum-normed space of measurable complex-valued functions $f$ on $E$ such that 
$$\|f\|_{w}  = \sup_{x\in E} \frac{|f(x)|}{w(x)} < +\infty.$$
\noindent Let us observe that $\mu\in\cB'_{w}$ if $\mu(w)<+\infty$. 
In particular we have $\pi\in\cB'_{v}$ by hypothesis. Clearly, $v$-geometrical ergodicity means that 
$Q$ is strongly ergodic w.r.t.~$\cB_{v}$.  

\noindent Let $0<\theta\leq1$. For the sake of simplicity, we slightly abuse notation below by writing 
$\cB_{\theta} = \cB_{{v^\theta}}$ and $\|\cdot\|_{\theta} = \|\cdot\|_{v^\theta}$. In particular $\cB_{1} = \cB_{v}$ and $\|\cdot\|_{1} = \|\cdot\|_{v}$. \\
The next lemma will be repeatedly used below (here $\xi$ is only supposed to be measurable). 

\noindent{\bf Lemma 10.1.} {\it  Condition ($\widehat{K}$) of Section 5.2 holds on $\cB=\cB_{\theta}$, with $\widehat{\cB} = \L^1(\pi)$. }

\noindent{\it Proof.} The property (K1) of Section~1 on $\cB_{\theta}$ (i.e.~$(X_n)_{n\geq0}$ is $v^\theta$-geometrically ergodic)  follows from the well-known link between $v$-geometric ergodicity and the so-called drift criterion \cite{mey}. 
More precisely, under the aperiodicity and $\psi$-irreducibility hypotheses, the $w$-geometric ergodicity 
for some $w : E\r [1,+\infty[$ is equivalent to the following condition: 
there exist $r<1$, $M\geq0$ and a petite set $C\in\cE$ such that $Qw_0 \leq r w_0 + M 1_C$, where $w_0$ is a function equivalent 
to $w$ in the sense that $c^{-1}\, w \leq w_0 \leq c\, w$ for some $c\in\R^*_+$. 
From that and since the function $t\mapsto t^\theta$ is concave on $\R_+$, $v$-geometric ergodicity implies, 
by virtue of Jensen's inequality, that \\[0.15cm]
\indent $\displaystyle \ \ \ \ \ \ \ \  \ \ \ \ \ \ \ \ \ \ \ \ \ \ \ \ 
Q(v_0^\theta) \leq (r v_0 + M 1_C)^\theta \leq r^\theta v_0^\theta + M^\theta 1_C$, \\[0.15cm]
where $v_0$ stands for some function equivalent to $v$. Thus $(X_n)_{n\geq0}$ is $v^\theta$-geometrically ergodic. \\
Besides, since $\pi(|e^{i\langle h,\, \xi\rangle}-1|\, |f|) 
\leq  \|f\|_{\theta}\, \pi(|e^{i\langle h,\, 
\xi\rangle}-1|\, v^\theta)$
for $f\in\cB_{\theta}$, we have ($\widehat{K2}$) (use Lebesgue's theorem and Remark (a) of Section 4). Besides we have  ($\widehat{K3}$) by Remark (b) of Section 4. Since $\cB_{\theta}$ is a Banach lattice, the property  ($\widehat{K4}$) w.r.t.~$\cB_{\theta}$ can be deduced from the abstract statement \cite[Cor. 1.6]{rab-wolf}. A simpler proof based on \cite{hen4} is presented in \cite{ihp2}. \fdem

\noindent If $|\xi|_2^2\leq C\, v$ for some $C>0$, then $(\frac{S_n}{\sqrt n})_n$ 
converges to a normal distribution $\cN(0,\Gamma)$ for any initial distribution. This is a classical result \cite{mey} which 
can be also deduced, in the stationary case, from the statements of Section 2. Indeed, the condition $|\xi|_2^2\leq C\, v$ 
implies that the 
coordinate functions $\xi_i$ of $\xi$ belong to the space $\cB_{\frac{1}{2}}$. Since $\pi(v) < +\infty$, we have 
$\cB_{\frac{1}{2}}\hookrightarrow\L^2(\pi)$, and the previous lemma shows that $Q$ is strongly ergodic on $\cB_{\frac{1}{2}}$. 
So the desired c.l.t.~follows from Proposition 2.2 and Corollary 2.1 both applied with 
$\cB = \cB_{\frac{1}{2}}$. 

\noindent Recall that, without 
additional assumptions, 
this central limit theorem does not hold under the weaker condition $\pi(|\xi|_2^2)<+\infty$ 
(see \cite{jones}). In the same way, the limit theorems below will hold under moment conditions of the type 
$|\xi|_2^\alpha\leq C\, v$ with some suitable exponent $\alpha\geq2$, and some 
positive constant $C$. So $\alpha$ will measure the 
order in these moment conditions, and we are going to 
see that, except for the multidimensional Berry-Esseen
theorem, it is 
similar (possibly up to $\varepsilon>0$) to that of the i.i.d.~case.
 
\noindent {\it The hypotheses of Assertions (a)-(d) 
below will imply that the above cited c.l.t.~holds, and it will be then understood that 
$\Gamma$ is non-singular (this means $\sigma^2>0$ 
in case $d=1$), hence we have (CLT) of Section~5.1. The nonlattice condition below is that of Proposition 5.4. Finally we suppose that the initial distribution $\mu$ is 
such that $\mu(v)<+\infty$. }

\noindent \noindent {\bf Corollary 10.2.} {\it \\
{\bf (a)} If $|\xi|_2^2\leq C\, v$ and $\xi$ is 
nonlattice, then we have (LLT) of Theorem 5.1 with $\cB = \cB_{\frac{1}{2}}$
and $\tilde \cB = \cB_{v}$. \\
{\bf (b)} (Case $d=1$)  If $|\xi|^3\leq C\, v$, then the uniform Berry-Esseen estimate holds: 
$\displaystyle\Delta_n = O(n^{-\frac{1}{2}})$. \\[0.2cm]
{\bf (c)} (Case $d = 1$)  If $|\xi|^\alpha\leq C\, v$ with some $\alpha>3$ and $\xi$ is nonlattice, 
then the first-order Edgeworth expansion (E) of  Theorem 8.1 holds. \\[0.15cm]
{\bf (d)} If $\vert \xi\vert_2^\alpha\le Cv$ with some 
$\alpha>\max\left(3,\lfloor d/2\rfloor+1\right)$, then the (Prohorov) Berry-Esseen estimate holds: 
${\cal P}\left(\mu_*\left({S_n\over\sqrt{n}}\right),{\cal N}(0,\Gamma)\right) = O(n^{-1/2})$. }

\noindent From the usual spectral method, (a) was established in \cite{sev} for bounded functionals $\xi$. 
Assertion (a)  extends the result of \cite{ihp1} stated under a kernel condition on $Q$. From Bolthausen's theorem 
\cite{bolt}, the one-dimensional uniform Berry-Esseen theorem holds under $\P_\pi$ (stationary case) if $\pi(|\xi|^p)<+\infty$ for some $p>3$. Assertion (b), already presented in \cite{ihp2}, extends this result to the 
non-stationary case under an alternative third-order moment condition. Assertion~(c) was established in  \cite{kontmey} 
for bounded functional $\xi$, and (d) is new to our knowledge. 

\noindent {\it Proof of Corollary 10.2.} Set $\widetilde{\cB} :=\cB_1 = \cB_{v}$. From Lemma~10.1, we have on each $\cB_{\theta}$: (K1), ($\widetilde{K3})$ (see Rk.~(c) in Sect.~4), and we have (S) if and only if $\xi$ is nonlattice (Prop.~5.4). \\[0.1cm]
{\it (a)} Since $g:=|\xi|_2\in \cB_{\frac{1}{2}}$, 
one gets: \\[0.15cm]
\indent $\displaystyle \ \ \ \ \ \ \ \ \ \ \ \ \ 
\forall f\in \cB_{\frac{1}{2}},\ \ |Q(t+h)f-Q(t)f| 
\leq Q\big(|e^{i\langle h,\, \xi\rangle}-1|\, 
|f|\big) \leq |h|_2\, \|g\|_{\frac{1}{2}}\, \|f\|_{\frac{1}{2}}Qv$, \\[0.15cm] 
and since $\frac{Qv}{v}$ is bounded, this proves 
($\widetilde{K2})$, hence ($\widetilde{K}$), with 
$\cB= \cB_{\frac{1}{2}}$. So Theorem 5.1 applies.  \\[0.1cm]
{\it (b)} Since (K1) holds on $\cB_{\frac{1}{3}}$ 
and $\cB_{\frac{2}{3}}$, 
and $\cB_{\frac{1}{3}}\hookrightarrow\L^3(\pi)$, 
$\cB_{\frac{2}{3}}\hookrightarrow\L^{\frac{3}{2}}
(\pi)$, we have (G1) (G2), so (CLT') (Section 6). We have 
$|Q(t)f-Qf| \leq |t|\, \|\xi\|_{\frac{1}{3}}\, \|f\|_{\frac{1}{3}}Qv^{\frac{2}{3}}$ for all $f\in \cB_{\frac{1}{3}}$, and since 
$\frac{Qv^{\frac{2}{3}}}{v} \leq \frac{Qv}{v}$, one gets $\|Q(t)-Q\|_{\scriptsize{\cB_{\frac{1}{3}},\cB_1}}=O(|t|)$. So  Theorem 6.1 applies with $\cB= \cB_{\frac{1}{3}}$. \\[0.1cm] 
Using the next proposition, Assertions (c) and (d) follow from Theorems 8.1 and 9.1. \fdem

\noindent {\bf Proposition 10.3.} {\it If $|\xi|_2^\alpha\leq C\, v$ with $\alpha > m$ ($\ m\in\N^*$), 
then $\cC(m)$ holds with $\cB=\cB_a$, $\widetilde{\cB} = \cB_1$, for any $a>0$ such that:  $a +{m\over\alpha} < 1$.}

\noindent{\it Proof.} For convenience, let us assume
that $d=1$. The extension to $d\geq2$ 
is obvious by the use of partial derivatives. Let $\varepsilon>0$ such that $a + {m+(2m+1)\varepsilon\over\alpha}\leq1$. We take $I=[a,1]$, $\cB_\theta$ ($\theta\in I$) as above defined, and we consider 
$T_0(\theta)=\theta+{\varepsilon\over\alpha}$, $T_1(\theta)=\theta+{1+\varepsilon\over\alpha}$. Recall that we set: $Q_k(t)(x,dy) = i^k\xi(y)^ke^{it\xi(y)}\, Q(x,dy)\, $ ($k\in\N,\ t\in\R,\ x\in E$).
With these notations, the proof of $\cC(m)$ is a consequence of the two following lemmas. \fdem 

\noindent {\bf Lemma 10.4.} {\it For any $k=0,\ldots,m$ and $\theta,\theta'>0$ such that $\theta + \frac{k}{\alpha} < \theta'\leq1$, 
we have $Q_k\in\cC^0(\R,\cB_\theta,\cB_{\theta'})$.  }

\noindent{\it Proof.} Let $0<\delta\leq1$ such that $\theta + \frac{k+\delta}{\alpha} \leq \theta'$. Using the inequality 
$|e^{iu}-1|\leq 2|u|^{\delta}$ ($u\in \R$), one gets for $t,t_0\in \R$ and $f\in\cB_{\theta}$: 
$$|Q_k(t)f-Q_k(t_0)f|\leq Q\big(|\xi|^k\, |e^{i(t-t_0)\xi}-1|\, |f|\big) \leq 2\, C^{\frac{k+\delta}{\alpha}}\, |t-t_0|^{\delta}\, 
\|f\|_{\theta} Q(v^{\frac{k+\delta}{\alpha}+\theta}),$$
hence $\|Q_k(t)f-Q_k(t_0)f\|_{{\theta'}}  \leq  2\, C^{\frac{k+\delta}{\alpha}}\, |t-t_0|^{\delta}\, 
\|f\|_{\theta}\, \|Q(v^{\theta'})\|_{{\theta'}}$. \fdem

\noindent {\bf Lemma 10.5.} {\it   For any $k=0,\ldots,m-1$ and $\theta,\theta'>0$ such that 
$\theta+\frac{k+1}{\alpha} < \theta'\leq 1$, we have $Q_k\in\cC^1(\R,\cB_{\theta},\cB_{\theta'})$  with $Q_k' = Q_{k+1}$. }

\noindent{\it Proof.} Let $0<\delta\leq1$ such that $\theta +\frac{k+1+\delta}{\alpha} \leq \theta'$. Using 
$|e^{iu}-1-iu|\leq 2|u|^{1+\delta}$ and proceeding as above, one gets 
$\|Q_k(t)f - Q_k(t_0)f - (t-t_0)Q_{k+1}(t_0)f\|_{\theta'} \leq 
2\, C^{\frac{k+1+\delta}{\alpha}}  |t-t_0|^{1+\delta} \|f\|_{\theta}\, \|Q(v^{\theta'})\|_{{\theta'}}$ for $t_0,t\in \R$ and 
$f\in\cB_{\theta}$. Since 
$Q_{k+1}\in\cC^0(\R,\theta,\theta')$, 
this yields the desired statement. \fdem 

\noindent {\it Remark.} The above proof shows that Assertion (b) of Corollary 10.2  holds 
under the alternative following hypotheses: 
$(X_n)_{n\geq0}$ is $v^{\frac{2}{3}}$-geometrically ergodic, $\mu(v^{\frac{2}{3}})<+\infty$, 
$|\xi|^3\leq C\, v$, and finally $\pi(v) < +\infty$, in order to have $\cB_{\frac{1}{3}}\hookrightarrow \L^3(\pi)$ and $\cB_{\frac{2}{3}}\hookrightarrow\L^{\frac{3}{2}}(\pi)$ (use $\widetilde{\cB}=\cB_{\frac{2}{3}}$). 
%============================================================================================
%\noindent {\bf 11. Applications to iterative  Lipschitz models.}
\section{Applications to iterative  Lipschitz models}
%===============================
%\noindent {\bf 11.1. Iterative  Lipschitz models.} 
\subsection{Iterative  Lipschitz models} 
Here $(E,d)$ is a non-compact metric space in which every closed ball is compact. We endow it with its
Borel $\sigma$-field $\cE$. Let $(G,\cG)$ be a measurable space,
let $(\theta_n)_{n\geq 1}$ be a sequence of  i.i.d.r.v.~taking values in
$G$. Let $X_0$ be a $E$-valued r.v.~independent of $(\theta_n)_n$, and  
finally let
$F : E\times G\r E$ be a measurable function. We set 
\\[0.15cm]
\indent $\displaystyle\ \ \ \ \ \ \ \ \  \ \ \ \ \ \ \ \ \ \  \ \ \ \  
\ \ \ \ \ \  \ \ \ \  \ \ \ \ \
X_n = F(X_{n-1},\theta_n),\ \ n\geq 1.$ \\[0.15cm]
\noindent For $\theta\in G$, $x\in E$, we set 
$F_\theta x =  
F(x,\theta)$ and we suppose that
$F_\theta : E\r E$ is Lipschitz continuous. Then $(X_n)_{n\geq1}$ is called an
iterative  Lipschitz model \cite{diaco} \cite{duf}. It is a Markov chain and its
transition probability is:\\[0.15cm]
\indent $\displaystyle\ \ \ \ \ \ \ \ \  \ \ \ \ \ \ \ \ \ \  \ \ \ \  
\ \ \ \ \ \
\ \ \ \  \ \ \ \ \ Qf(x) = \E[\, f(F(x,\theta_1))\, ]$. \\[0.15cm]
Let $x_0$ be a fixed point in $E$. As in \cite{duf}, we shall appeal  
to the following r.v:
$$\cC = \sup\bigg\{\frac{d(F_{\theta_1}x,F_{\theta_1}y)}{d(x,y)},\  
x,y\in E,\ x\neq y\bigg\}\ \ \ \ \mbox{and}\ \ \ \
\cM  = 1 + \cC + d\big(F(x_0,\theta_1),x_0\big).$$
As a preliminary, let us present a sufficient condition for the  
existence and the uniqueness of an
invariant distribution. The following proposition is proved in  
\cite{aap} (Th.~I).

\noindent {\bf Proposition 11.1.} {\it  Let $\alpha\in(0,1]$,  
$\eta\in\R_+$. Under the moment condition
$\E[\cM^{\alpha(\eta+1)}] < +\infty$ and the mean contraction  
condition $\E[\cC^\alpha\, \max\{\cC,1\}^{\alpha\eta}]< 1$,
there exists a unique stationary distribution, $\pi$, and we have
$\pi(d(\cdot,x_0)^{\alpha(\eta+1)})<+\infty$. }

\noindent More precise statements can be found in the literature (see  
e.g \cite{diaco} \cite{duf}). However,
the hypotheses occurring in Proposition 11.1 are convenient in our  
context and are similar to those introduced later.

\noindent Finally, we shall suppose that $\xi$ satisfies the following condition, with given $S,s\geq0$: 

\noindent $\displaystyle {\bf (L)_s} \ \ \ \ \ \ \ \ \ \ \ \ \
\forall (x,y)\in E\times E,\ \ |\xi(x)-\xi(y)|_2 \leq S\, d(x,y)\, \big[1+d(x,x_0)+d(y,x_0)\big]^s$.

\noindent  For convenience, Condition $\bf (L)_s$ has been stated as a  
weighted-Lipschitz condition w.r.t.~the distance
$d(\cdot,\cdot)$ on $E$. However, by replacing $d(\cdot,\cdot)$ with  
the distance $d(\cdot,\cdot)^a$ ($0<a\leq1$),
Condition $\bf (L)_s$ then corresponds to the general weighted-H\"older  
condition of \cite{duf}. \\[0.2cm]
\noindent  Section 11.2 below will introduce weighted H\"older-type  
spaces and investigate
all the hypotheses of the previous sections.
Using these preliminary statements, we shall see in Section 11.3 that  
the limit theorems of the preceding sections then apply
to $(\xi(X_n))_n$ under some mean contraction and moment conditions.  
These conditions will focus on the random variables $\cC$, $\cM$ and
will depend on the real number $s$ of Condition $\bf (L)_s$.

\noindent To compare with the i.i.d.~case, let us summarize the results  
obtained in Section 11.3 in the following special setting~:
$(X_n)_n$ is a ${\R}^d$-valued iterative Lipschitz sequence such that  
$\cC<1$ a.s.. For convenience we also assume that $(X_n)_n$ is
stationary, with stationary distribution $\pi$, and we consider the  
random walk  associated to
$\xi(x) = x - \E_\pi[X_0]$, that is: \\[0.1cm]
\indent $\displaystyle\ \ \ \ \ \ \ \ \  \ \ \ \ \ \ \ \  \  \ \ \ \ \  
\ \ \  \  \ \ \ \ \ \ \ \
S_n = X_1+\ldots+X_n-n\E_\pi[X_0]$. \\[0.1cm]
Finally suppose that  $\E[\, \cM^{2}] < +\infty$. Then the sequence  
$(\frac{S_n}{\sqrt n})_n$ converges to
$\cN(0,\Gamma)$ \cite{benda}, and we assume  that $\Gamma$ is invertible.
Corollaries of Section 11.3 will then provide the following results~:\\[0.2cm]
{\it (i) Local limit theorem~:} $\xi$ nonlattice $\ \Rightarrow\ $  
(LLT) of Section 5.1 with for instance $f=h=1_E$,  \\[0.2cm]
{\it (ii) ($d=1$) Uniform Berry-Esseen type theorem~:} $\E[\cM^3] < +\infty\  
\Rightarrow\ \Delta_n = O(n^{-\frac{1}{2}})$, \\[0.25cm]
{\it (iii) ($d=1$) First-order Edgeworth expansion~:} 
$\E[\cM^{3 +  \varepsilon}] <+\infty$, $\xi$ nonlattice $\ \Rightarrow\ $ (E) of Section 8,
\\[0.25cm]
{\it (iv) multidimensional  Berry-Esseen theorem (with Prohorov metric)~:} 
$\E[\cM^{m +  \varepsilon}] <+\infty$ with $m=\max\left(3,\lfloor d/2\rfloor+1\right)$ $\ \Rightarrow\ $ the conclusion of Theorem 9.1 holds.

\noindent More generally, the previous assertions apply to  
$(\xi(X_n))_n$ whenever
$\xi$ is a  Lipschitz continuous function on $E$ (i.e.~$\bf (L)_s$ holds with $s=0$). 

\noindent {\it Example. The autoregressive models. } \\
A simple and typical example is the autoregressive chain  
defined in $\R^d$ by \\[0.1cm]
\indent $\displaystyle\ \ \ \ \ \ \ \ \  \ \ \ \ \ \ \ \  \  \ \ \ \ \  
\ \ \  \  \ \ \ \ \ \ \ \
X_n = A_n\, X_{n-1} + \theta_n$ ($n\in\N^*$), \\[0.1cm]
where $(A_n,\theta_n)_{n\geq1}$ is a i.i.d.~sequence of r.v.~taking values in $\cM_d(\R)\times\R^d$, independent of $X_0$. 
($\cM_d(\R)$ denotes the set of real $d\times d$-matrices.) Assume that we have $|A_1| < 1\, $ a.s., where $|\cdot|$  denotes here both some norm on $\R^d$ and the associated matrix norm. Taking the distance $d(x,y) = |x-y|$ on $\R^d$, we  have $\cC=|A_1|$ and $\cM \leq 2 + |\theta_1|$. So the above moment conditions in (i)-(iv) only concern $|\theta_1|$. \\[0.15cm]
\noindent The special value $A_n=0$ corresponds  to the i.i.d.~case ($S_n = \theta_1+\ldots+\theta_n-n\E[\theta_1]$), and we can see that the moment conditions in (i)-(iii) are then 
optimal for Statements (i) (ii), and optimal up 
to $\varepsilon>0$ for Statement (iii). 

\noindent Let us mention that \cite{gui-lepage} investigates the convergence to stable laws for the random walk 
associated to the above autoregressive model $(X_n)_n$ (case $d=1$) and to $\xi(x)=x$. By using the Keller-Liverani theorem, \cite{gui-lepage} presents very precise statements, similar to the i.i.d.~case, in function of the "heavy tail" property of the stationary distribution of $(X_n)_n$.   
%===========================================
%\noindent {\bf 11.2. Preliminary results.} 
\subsection{Preliminary results} 
The weighted H\"older-type spaces, introduced in \cite{lep83}, 
have been used by several authors for proving quasi-compactness under some contracting property \cite{mira,pei}. 
Here we slightly modify the definition of these spaces by considering two positive parameters $\beta$ and $\gamma$ in the weights. This new definition is due to D.~Guibourg.  

\noindent Let us consider $0<\alpha\leq 1$ and $0<\beta\leq\gamma$. For $x\in E$, we set $p(x) = 1+ \, d(x,x_0)$, and for $(x,y)\in E^2$, we set  
$$\Delta_{\alpha,\beta,\gamma}(x,y) = p(x)^{\alpha\gamma}\, p(y)^{\alpha\beta} + p(x)^{\alpha\beta}\,   p(y)^{\alpha\gamma}.$$
Then $\cB_{\alpha,\beta,\gamma}$ denotes the space of $\C$-valued functions on $E$ satisfying the following condition
$$m_{\alpha,\beta,\gamma}(f) =  
\sup\bigg\{\frac{|f(x)-f(y)|}{d(x,y)^\alpha\, \Delta_{\alpha,\beta,\gamma}(x,y)},\ x,y\in E,\  
x\neq y\bigg\}\, <\, +\infty.$$
Set  $\displaystyle \ |f|_{\alpha,\gamma}  = \sup_{x\in E}\  
\frac{|f(x)|}{p(x)^{\alpha(\gamma+1)}}$ and
$\|f\|_{\alpha,\beta,\gamma} = m_{\alpha,\beta,\gamma}(f) + |f|_{\alpha,\gamma}$. Then $(\cB_{\alpha,\beta,\gamma},\|\cdot\|_{\alpha,\beta,\gamma})$  is a Banach space. In the special case $\gamma = \beta$, we shall simply denote $\cB_{\alpha,\gamma} = \cB_{\alpha,\beta,\gamma}$. 

\noindent The next result which concerns Condition (K1) on $\cB_{\alpha,\beta,\gamma}$ is established in \cite{aap} [Th.~5.5] in the case $\beta=\gamma$. Since the extension to the case $0<\beta\leq\gamma$ is very easy, we give the following result without proof. 

\noindent {\bf Proposition 11.2.} {\it If $\E[\cM^{\alpha(\gamma+1)} + \cC^\alpha\, \cM^{\alpha(\gamma+\beta)}]<  
+\infty$, $\ \E[\cC^\alpha\, \max\{\cC,1\}^{\alpha(\gamma+\beta)}]< 1$, then $Q$ is  
strongly ergodic on  $\cB_{\alpha,\beta,\gamma}$. }

\noindent Now we give a sufficient condition for the central limit  
theorem in the stationary case.
Similar statements are presented in \cite{duf}, and in \cite{benda}  
when $\xi$ is Lipschitz continuous
(i.e.~$s=0$ in $\bf (L)_s$).

\noindent {\bf Proposition 11.3.} {\it If $\E[\, \cM^{2s+2} +  
\cC^{\frac{1}{2}}\, \cM^{2s+1}\, ] < +\infty$ and
$\E[\cC^{\frac{1}{2}}\, \max\{\cC,1\}^{2s+\frac{3}{2}}]< 1$,
then, under $\P_\pi$, $(\frac{S_n}{\sqrt n})_n$ converges to a normal  
distribution $\cN(0,\Gamma)$. }

\noindent{\it Proof.} We apply Proposition 11.1 with $\alpha=\frac{1}{2}$ and $\eta = 4s+3$. This yields the existence and the uniqueness of $\pi$, and $\pi(d(\cdot,x_0)^{2s+2})<+\infty$. Here we consider $\gamma = \beta = 2s+1$ and the corresponding space $\cB = \cB_{\frac{1}{2},\gamma}$. For $f\in\cB$, 
we have $|f|\leq |f|_{\frac{1}{2},2s+1}\, p(x)^{s+1}$. Thus $\cB\hookrightarrow\L^2(\pi)$.  
Besides, from $\bf (L)_s$, it can be easily seen that the coordinate functions of $\xi$ belong to $\cB$, and 
by  Proposition 11.2, $Q$ is strongly ergodic on $\cB$. We conclude by applying Proposition 2.2 and Corollary 2.1 
with $\cB$ as above defined. \fdem

\noindent The possibility of considering $\alpha<1$ as above is  
important. To see that, consider for instance the case $s=0$
(i.e.~$\xi$ is Lipschitz continuous on $E$). Then $\xi\in\cB_{1,\gamma}$ for any
$\gamma>0$, and we could also consider $\cB = \cB_{1,\gamma}$ in the  
previous proof, but it
is worth noticing that the condition $\cB_{1,\gamma}\hookrightarrow\L^2(\pi)$
would then require the moment condition   
$\pi(d(\cdot,x_0)^{2(1+\gamma)})<+\infty$ which is stronger than
$\pi(d(\cdot,x_0)^{2})<+\infty$ used above. Anyway, we shall often  
appeal below to the conditions $s+1\leq \beta \leq \gamma$ and  
$\E[\cM^{\alpha(\gamma+1)}] < +\infty$. If $s=0$ and $\beta=\gamma=1$, then the previous moment condition
is  $\E[\cM] < +\infty$ if $\alpha=\frac{1}{2}$, while it is  $\E[\cM^{2}]  < +\infty$  if $\alpha=1$. 

\noindent Now we investigate the action of the Fourier kernels $Q(t)$  
on  the space $\cB_{\alpha,\beta,\gamma}$.
The proofs of Propositions 11.4-8 below present no theoretical problem.  
However the presence of Lipschitz coefficients
in the definition of  $\cB_{\alpha,\beta,\gamma}$ makes the computations  
quite more technical than  those seen for the $v$-geometrically ergodic Markov chains.
For convenience, these proofs are presented in Appendix B. 
The arguments will be derived
from \cite{aap}. However, the next four statements 
improve the corresponding ones in \cite{aap} (See Remark below).

\noindent {\bf Proposition 11.4.}  {\it Condition (K) of Section 4 holds on $\cB_{\alpha,\beta,\gamma}$ if we have $s+1\leq\beta\leq\gamma$ and $\E\left[\, \cM^{\alpha(\gamma+1)} + \cC^\alpha\, \cM^{\alpha(\gamma+\beta)}\, \right] <  +\infty$, 
$\displaystyle \ \E\left[\, \cC^\alpha\, \max\{\cC,1\}^{\alpha(\gamma+\beta)}\, \right] < 1$.}

\noindent {\bf Proposition 11.5}. {\it We have  
$\|Q(t+h)-Q(t)\|_{_{\scriptsize{\cB_{\alpha,\beta,\gamma}\, ,\, \cB_{\alpha,\beta,\gamma'}}}}\r0$
when  $t\r0$ if the following conditions hold: $\ s+1 \leq \beta \leq \gamma < \gamma'$ and $\E\left[\, \cM^{\alpha(\gamma'+1)}\,  
+\, \cC^\alpha\,\cM^{\alpha(\gamma' + \beta)}\, \right] < +\infty$.}

\noindent {\bf Proposition 11.6.}  {\it We have  $\|Q(t)-Q\|_{_{\scriptsize{\cB_{\alpha,\beta,\gamma}\, ,\, \cB_{\alpha,\beta,\gamma'}}}} = O(|t|)$ if the following conditions hold: 
$\ s+1\le \beta\leq\gamma$, $\, \gamma' \geq \gamma + \frac{s+1}{\alpha}$, and $\displaystyle \E\left[\,  \cM^{\alpha(\gamma'+1)} + \cC^\alpha\, \cM^{\alpha(\gamma'+\beta)} \,\right] <+\infty$. }

\noindent {\bf Proposition 11.7.}  {\it We have $\cC(m)$ of Section 7.1 ($m\in\N^*$) with $\cB = \cB_{\alpha,\beta,\gamma}$ and $\widetilde{\cB} = \cB_{\alpha,\beta,\gamma'}$ if we have $s+1\leq\beta\leq\gamma$, $\, \gamma' > \gamma + \frac{m(s+1)}{\alpha}$, and }
$$\E\left[\, \cM^{\alpha(\gamma'+1)}\, +\, \cC^\alpha\cM^{\alpha(\gamma'+ \beta)}\,\right] <+\infty\ \ \ \ \ \ 
\E\left[\, \cC^\alpha\, \max\{\cC,1\}^{\alpha(\gamma'+\beta)}\right] < 1.$$
\noindent Concerning the spectral condition (S) of Section 5.1, we now study the  
possibility of applying
the results of Section 5.2. Observe that this cannot be done with the  
help of Proposition 11.4 because Condition (K) only concerns $Q(t)$ for $t$ near $0$. By considering another  
auxiliary semi-norm  on $\cB_{\alpha,\beta,\gamma}$,
we shall prove in Appendix B.5 the following result for which
the hypotheses are somewhat more restrictive than those of Proposition 11.4.

\noindent {\bf Proposition 11.8.}  {\it Assume $s+1 < \beta \leq \gamma<\gamma'$, $\, \displaystyle \E\left[\cM^{\alpha(\gamma'+1)}+\cC^\alpha\,\cM^{\alpha(\gamma'+\beta)}\right] < +\infty$ and
$\displaystyle \E\big[\cC^{\alpha}\, \max\{\cC,1\}^{\alpha(\gamma+\beta)}\big]< 1$. 
Then Condition (S) holds on $\cB_{\alpha,\beta,\gamma}$ if and only if $\xi$ is non-arithmetic w.r.t.~$\cB_{\alpha,\beta,\gamma}$. If $\xi$ is nonlattice, the two previous equivalent conditions hold. }

\noindent {\bf Remark.} The possibility of considering the spaces $\cB_{\alpha,\beta,\gamma}$ with $\beta\neq \gamma$ 
is important, in particular to apply Proposition 11.7. Indeed, let us assume $\cC<1$ a.s. and consider the case $s=0$ to simplify. 
Then the condition for $\cC(m)$ is 
$\E[\cM^{\alpha(\frac{m}{\alpha} + \gamma + \beta)+\varepsilon}] <+\infty$ 
(for some $\varepsilon>0$), where $\beta$ and $\gamma$ are such that $1\leq\beta\leq\gamma$.  This condition can be 
rewritten as $\E[\cM^{m + \alpha(\gamma + \beta)
+\varepsilon}] <+\infty$. Consequently, under a moment assumption of the form $\E[\cM^{m + \varepsilon_0}] <+\infty$ for some $\varepsilon_0>0$, we can choose $\alpha$ sufficiently small in order to ensure Condition $\cC(m)$.  \\ 
Actually, the condition $\E[\cM^{\alpha(\frac{m}{\alpha} + \gamma + \beta)+\varepsilon}] <+\infty$ is useful for proving (K1) on the biggest space occurring in $\cC(m)$. It is worth noticing that, when working with the weights defined in \cite{mira,pei,aap} (which corresponds to our weights in the special case $\beta=\gamma$), then Condition (K1) must be satisfied on $\cB_{\alpha,\gamma',\gamma'}$ with $\gamma' > \gamma + \frac{m}{\alpha}$: this then requires the moment 
condition $\E[\cM^{2\alpha(\frac{m}{\alpha} + 
\gamma)+\varepsilon}]  = \E[\cM^{2m  + 2\alpha\gamma
+\varepsilon}] <+\infty$ (apply Prop.~11.2 on $\cB_{\alpha,\gamma',\gamma'}$), whose order is greater than $2m$. 
Our parameter $\beta$ enables us to avoid this drawback. 
%=========================================
%\noindent {\bf 11.3. Limit theorems for $(\xi(X_n))_n$.} 
\subsection{Limit theorems for $(\xi(X_n))_n$} 
\noindent {\it  The hypotheses of Corollaries 11.9-12 below will imply
those of Proposition 11.3. Consequently the c.l.t.~stated in this  
proposition will hold automatically, and it will be understood that
$\Gamma$ is non-singular.  }

\noindent Concerning the next conditions imposed on the initial 
distribution $\mu$, it is worth noticing that,
if $\mu(d(\cdot,x_0)^{\alpha(1+\gamma)})<+\infty$, then  
$\mu\in\cB_{\alpha,\beta,\gamma}'$. The conditions imposed on $\mu$ below will be always satisfied
for $\mu=\pi$ or $\mu=\delta_x$ ($x\in E$) (for $\pi$ it comes
from Proposition~11.1). 

\noindent {\bf Local limit theorem} ($d\geq1$).  

To present a simple application of  Theorem 5.1, let us simply investigate
Statement (LLT) of Section 5.1 with $f=h=1_E$. We want to prove that, for any  
compactly supported
continuous function $g : \R^d\r\R$, we have \\[0.14cm]
\indent $\displaystyle {(LLT')} \ \ \ \lim_n 
\sup_{a\in\R^d} \bigg\vert\sqrt{\det\Gamma}\, (2\pi n)^{\frac{d}{2}}\,
\E_\mu[\, g(S_n-a)\, ] - 
e^{-\frac{1}{2n}\langle  
\Gamma^{-1}a,a\rangle}\, \int_\R g(x)dx\bigg\vert = 0$. \\[0.2cm]
\noindent {\bf Corollary 11.9.}  {\it Suppose that $\E\big[\, \cM^{2s+2} +  
\cC^{\frac{1}{2}}\, \cM^{2s+1+\delta}\big] < +\infty$ for some $\delta>0$, that 
$\E\big[\cC^{\frac{1}{2}}\, \max\{\cC,1\}^{2s+\frac{3}{2}}\big]< 1$, that  
$\xi$ satisfies $\bf (L)_s$ and is nonlattice, and 
finally that we have $\mu(d(\cdot,x_0)^{\frac{2+s+\varepsilon_0}{2}})<+\infty$ for some $\varepsilon_0>0$. Then we have (LLT'). }

\noindent {\it Proof.} By using the above preliminary 
statements, let us prove that the hypotheses of 
Theorem~5.1  hold. We have (CLT) (Prop.~11.3). Let 
$\alpha = \frac{1}{2}$, $0<\varepsilon\leq\min
\{\frac{1}{2},\frac{2\delta}{3},\frac{\varepsilon_0}
{2}\}$, $\beta = \gamma =  s+1+\varepsilon$, and 
$\gamma' = \gamma + \varepsilon = s+1+2\varepsilon$. We set $\cB = \cB_{\frac{1}{2},\gamma}$. We have (S) and (K1) on 
$\cB$ (Prop.~11.8, 11.2). 
Besides, with  $\widetilde{\cB} = \cB_{\frac{1}{2},
\beta,\gamma'}$, we have ($\widetilde{K2}$) 
(Prop.~11.5), and ($\widetilde{K3}$) 
(use Prop.~11.4, $\widetilde{\cB}\hookrightarrow\L^1
(\pi)$ and Rk.~(c) in Section~4). Hence 
($\widetilde{K}$) holds. Finally, 
our assumption on $\mu$ implies $\mu\in\widetilde{\cB}'$. 
\fdem

\noindent According to the previous proof, the property (LLT) may be also investigated with functions $f\in\cB_{\frac{1}{2},\beta,\gamma}$ (for some suitable $s+1 < \beta \leq\gamma$), and the sufficient nonlattice condition can be replaced by the more precise non-arithmeticity condition (w.r.t.~$\cB_{\frac{1}{2},\beta,\gamma}$) of Proposition~5.3. Finally observe that, if $s=0$ (i.e.~$\xi$ is Lipschitz continuous on $E$), as for  
example $\xi(x) = \|x\|$, and if we have $\cC<1$ a.s.,
then (LLT') is valid under the expected moment condition $\E[\cM^2]  <  
+\infty$.

\noindent {\bf One-dimensional uniform Berry-Esseen theorem} ($d=1$). 

\noindent {\bf Corollary 11.10.} {\it Suppose $\E\big[\, \cM^{3(s+1)} +  
\cC^{\frac{1}{2}}\cM^{3s+2}\, \big]< +\infty$
and $\E\big[\cC^{\frac{1}{2}}\max\{\cC,1\}^{3s+\frac{5}{2}}\big]< 1$, that  $\xi$  
satisfies $\bf (L)_s$, and $\mu(d(\cdot,x_0)^{2(s+1)})<+\infty$. Then $\displaystyle\Delta_n =  O(n^{-\frac{1}{2}})$. }

\noindent {\it Proof.}  To apply Theorem 6.1, we have to prove (CLT')
of Section 6 and to find some spaces $\cB$ and  $\widetilde{\cB}$ on which
($\widetilde{K}$) holds with the additional condition  
$\|Q(t)-Q\|_{\scriptsize{\cB,\widetilde{\cB}}}=O(|t|)$. To investigate
(CLT'), we shall use the procedure based on conditions (G1)-(G2) (of Section 6). In particular this procedure
requires that $\xi\in\cB\hookrightarrow\L^3(\pi)$. Since $\bf (L)_s$  
implies $\xi\in\cB_{\frac{1}{2},2s+1}$, let us consider
$\cB = \cB_{\frac{1}{2},2s+1}$ (so here $\beta=\gamma=2s+1$). For $f\in\cB$, we have
$|f|\leq |f|_{\frac{1}{2},2s+1}\, p^{s+1}$, and since  $\pi(d(\cdot,x_0)^{3(s+1)})<+\infty$ (use Prop.~11.1 with
$\alpha = \frac{1}{2}$, $\eta = 6s+5$), one gets  $\cB\hookrightarrow\L^3(\pi)$. Now set $\widetilde{\cB} = \cB_{\frac{1}{2},\beta,4s+3}$. It can be easily seen that $\widetilde{\cB}$ 
contains all the functions $g^2$ with $g\in\cB$, and since  
each  $f\in\widetilde{\cB}$ satisfies $|f| \leq |f|_{\frac{1}{2},4s+3}\, p^{2(s+1)}$, one obtains $\widetilde{\cB}\hookrightarrow\L^{\frac{3}{2}}(\pi)$. We have (K1) on $\cB$ and $\widetilde{\cB}$ (Prop.~11.2). This gives (G1) (G2), hence (CLT'). Besides we have ($\widetilde{K3}$) (use Prop.~11.4, $\widetilde{\cB}\hookrightarrow\L^1(\pi)$ and Rk.~(c) in Section~4), and we have ($\widetilde{K2}$) (Prop.~11.5). Hence ($\widetilde{K}$). Finally, Proposition~11.6 yields $\|Q(t)-Q\|_{\scriptsize{\cB,\widetilde{\cB}}} = O(|t|)$, and $\mu(d(\cdot,x_0)^{2(s+1)})<+\infty$ implies that $\mu\in\widetilde{\cB}'$. \fdem

\noindent {\bf The first-order Edgeworth expansion} ($d=1$).  

For convenience, we investigate the property (E) of Theorem 8.1 under the hypothesis  
that $\cC$ is strictly contractive a.s..

\noindent {\bf Corollary 11.11.}  {\it Suppose that  $\cC<1$ a.s., that 
$\E[\cM^{3(s+1) + \varepsilon_0}] <+\infty$
for some $\varepsilon_0>0$, that  $\xi$ satisfies $\bf (L)_s$ and is  
nonlattice, and
$\mu(d(\cdot,x_0)^{3(s+1) + \varepsilon_0}) <+\infty $. Then we have (E).}

\noindent {\it Proof.}  To check the hypotheses of Theorem 8.1, first observe that the hypothesis $\cC<1$ a.s. implies  
$\E[\cC^{\alpha}\max\{\cC,1\}^b]< 1$ for
any $\alpha\in(0,1]$ and $b\geq0$. We have  
$\pi(d(\cdot,x_0)^{3(s+1)})  <+\infty$ (Prop.~11.1). From
$|\xi(x)| \leq p(x)^{s+1}$, it follows that $\pi(|\xi|^3) < +\infty$.
Let us prove that $\cC(3)$ holds w.r.t.~$\cB= \cB_{\alpha,\beta,\gamma}$ and $\widetilde{\cB} =  \cB_{\alpha,\beta,\gamma'}$ for suitable $\alpha,\beta,\gamma,\gamma'$. 
Let $\delta>0$, $\beta = \gamma = s+1 +\delta$, and let us choose $0<\alpha\leq1$ such that  
$\alpha(\gamma+2\delta+s+1)\leq \varepsilon_0 $. Let $\gamma' = \gamma +  \frac{3(s+1)}{\alpha} + \delta$. Then Proposition 11.7 yields the desired property. 
To study Condition (S) on $\cB_{\alpha,\beta,\gamma}$, use Proposition 11.8. 
Finally, we have $\alpha(\gamma'+1) =  3(s+1) + \alpha(\gamma + \delta +1) \leq  3(s+1) + \varepsilon_0$, so $\mu(d(\cdot,x_0)^{\alpha(\gamma'+1)}) <+\infty $. This proves that $\mu\in\widetilde{\cB}'$. \fdem

\noindent Other similar statements may be derived by proceeding as above.
For instance, let us consider $0<\alpha\leq1$ (fixed here),
$\beta = \gamma = s+1 +\delta$, and $\gamma' = \gamma + \frac{3(s+1)}{\alpha} +  
\delta$ with some small $\delta>0$, and suppose that we have
$\E[\cM^{\alpha(\gamma'+1)}\, +\, \cC^\alpha\cM^{\alpha(\gamma'+\beta)}] <+\infty$,
$\E[\cC^\alpha\, \max\{\cC,1\}^{\alpha(\gamma'+\beta)}] < 1$, and  
$\mu(d(\cdot,x_0)^{\alpha(\gamma'+1)}) <+\infty$. Then
we have (E) if $\xi$ is non-arithmetic w.r.t.~$\cB_{\alpha,\beta,\gamma}$.

\noindent {\bf The multidimensional Berry-Esseen theorem with the Prohorov  
distance} ($d\geq1$). 

Again we give a statement in the particular case when ${\cC}<1$ a.s.. From Theorem 9.1, we get the following.

\noindent {\bf Corollary 11.12.}  {\it Suppose  $\cC<1$ a.s.~and  
$\E[\cM^{m(s+1) + \varepsilon_0}] <+\infty$
for some $\varepsilon_0>0$ and with $m:=\max\left(3,\lfloor d/2\rfloor+1\right)$, that  
$\xi$ satisfies $\bf (L)_s$ and
$\mu(d(\cdot,x_0)^{m(s+1) + \varepsilon_0}) <+\infty $. Then the conclusion of Theorem
9.1 holds.}

\noindent {\it Proof.} Set $\beta = \gamma=1+s$. Let $\delta>0$, and $0<\alpha\leq1$ be such that $\alpha(\gamma+\delta+s+1)\leq \varepsilon_0$, and set $\gamma' = \frac{m(s+1)}{\alpha} + \gamma + \delta$. Then we have  $\cC(m)$ with $\cB={\cB}_{\alpha,\beta,\gamma}$ and $\widetilde{\cB}={\cB}_{\alpha,\beta,\gamma'}$ 
(by Proposition~11.7), and the hypothesis on $\mu$ gives $\mu\in \widetilde{\cB}'$. \fdem

\noindent {\bf Extension.}  Mention that all the previous statements  
remain valid when, in the hypotheses, the r.v.~$\cC$ is
replaced with the following one~:
$$\cC^{(n_0)} = \sup\bigg\{\frac{d(F_{\theta_1}\cdots F_{\theta_{n_0}}x\, ,\,
F_{\theta_1}\cdots F_{\theta_{n_0}}y)}{d(x,y)},\ x,y\in E,\ x\neq  
y\bigg\}\ \ (n_0\in\N^*).$$
The proofs of the preliminary statements of Section 11.2 are then similar. 
%=================================================
\section{More on non-arithmeticity and nonlattice conditions}
This section presents some complements concerning the spectral condition (S) of Section 5.1, in particular we 
prove Proposition 5.3 and specify Proposition 5.4. 
%====================================
\subsection{Proof of Proposition 5.3.} 
We assume that the assumptions ($\widehat{K}$) and (P) of Section 5.2 hold. Recall that Condition (S) on $\cB$ states that, for each compact set $K_0$ in $\R^d\setminus\{0\}$, there exist $\rho<1$ and $c\geq 0$ such that we have, for all $n\geq1$ and $t\in K_0$, $\|Q(t)^n\|_{_{\scriptsize \cB}} \leq c\, \rho^n$. \\[0.2cm]
We have to prove that (S)
is not true if and only if there exist $t\in\R^d$, $t\neq0$, $\lambda\in\C$, $|\lambda|=1$, a $\pi$-full $Q$-absorbing set $A\in\cE$, and a bounded element $w$ in $\cB$ such that 
$|w|$ is nonzero constant on $A$, satisfying: \\[0.2cm]
\indent $(*)\ \ \ \ \ \ \ \ \ \ \ \ \ \ \ \ \ \  
\displaystyle \forall x\in A,\ \ 
e^{i\langle  t,\xi(y)\rangle} w(y) = \lambda w(x)\ \ Q(x,dy)-\mbox{\it a.s.}$. 

\noindent {\bf Lemma 12.1.} {\it Let $t\in\R^d$ such that $r(Q(t)) \geq 1$. Then \\
(i) $r(Q(t)) = 1$ and $Q(t)$ is quasi-compact. \\
(ii) We have ($*$) with $\lambda$, $A$ and $w$ as above stated. }

\noindent{\it Proof of Assertion (i).} By ($\widehat{K4}$), 
we have $r_{ess}(Q(t)) < 1 \leq r(Q(t))$, thus $Q(t)$ is quasi-compact on $\cB$. Now let 
$\lambda$ be any  eigenvalue of modulus $r(Q(t))$, and let $f\neq0$ be an associated eigenfunction in $\cB$.    
Then $|\lambda|^n|f| = |Q(t)^nf| \leq Q^n|f|$, and (P) yields $|\lambda| \leq 1$. \fdem

\noindent By (i), there exist $\lambda\in\C$, $|\lambda|=1$ and $w\in\cB$, $w\neq 0$, such that $Q(t)w=\lambda w$. From $Q(t)^nw=\lambda^n w$, one gets $|w| \leq Q^n|w|$, and (P) then implies that 
$|w|\leq\pi(|w|)$, either everywhere on $E$, or $\pi$-a.s. on $E$, according that $\cB\subset\cL^1(\pi)$ or $\cB\subset\L^1(\pi)$. From now, if $\cB\subset\L^1(\pi)$, $w$ is replaced with any measurable function of its class, and for convenience, this function is still denoted by $w$. 
Since 
$v = \pi(|w|) - \vert w\vert \geq 0$ and $\pi(v)=0$, we have $|w| = \pi(|w|)$ $\pi$-a.s. 
Let us define the set \\[0.15cm]
\indent $\ \ \ \ \ \ \ \ \ \ \ \ \ \ \ \ \ \  \ \ \ \ \ \ \ \ \ \ \ \ \ \
A_0 = \{z\in E : |w(z)| = \pi(|w|)\}$. \\[0.15cm]
Then we have $\pi(A_0)=1$ (i.e.~$A_0$ is $\pi$-full). 

\noindent {\it Remark.} In the special case
when $\delta_x\in\cB'$ for all $x\in E$
(and when $\cB$ is stable under complex modulus), 
the proof of (ii) is presented in \cite{hulo} (Prop.~V.2), with the more precise conclusion: we have ($*$) with $w\in\cB\cap\cB^{^\infty}$ and $A=A_0$. Let us briefly recall the main arguments. From (K1), one can here deduce from the inequality $|w| \leq Q^n|w|$ that $|w|\leq \pi(|w|)$ everywhere on $E$. Thus $w\in\cB^{^\infty}$. 
Besides, the equality 
$Q(t)w(x) = \int_E e^{i\langle  t,\xi(y)\rangle} w(y)\, Q(x,dy) = \lambda w(x)$ is valid for all $x\in E$. Let $x\in A_0$. 
Then this equality and the previous inequality give ($*$). Finally ($*$) shows that $A_0$ is $Q$-absorbing. 

\noindent If $Q(t)w = \lambda\, w$ almost surely, the previous arguments must be slightly modified as follows. 

\noindent {\it Proof of (ii).}  First, by proceeding as in the proof of Proposition 2.4, one can easily get a $\pi$-full $Q$-absorbing set $B\subset A_0$. 
Besides the following set is clearly $\pi$-full: \\[0.15cm]
\indent $\ \ \ \ \ \ \ \ \ \ \ \ \ \ \ \ \ \ \ \ \ \ \ \ \ \
C = \{z\in E : \forall n\geq1,\, Q(t)^nw(z) = \lambda^n w(z)\}$. \\[0.15cm]
So the set $A=B\cap C$ is also $\pi$-full. Let $x\in A$. We have \\[0.15cm]
\indent $\ \ \ \ \ \ \ \ \ \ \ \ \ \ \ \ \ \  \ \ \ \ \ \ \ 
Q(t)w(x) = \int_E e^{i\langle  t,\xi(y)\rangle} w(y)\, Q(x,dy) = \lambda w(x)$. \\[0.15cm]
Since $Q(x,B) = 1$ ($B$ is $Q$-absorbing), one can replace $E$ by $B$ in the previous integral,  
and since $|\lambda^{-1} w(x)^{-1}\, e^{i\langle  t,\xi(y)\rangle} w(y)| = 1$ for all $y\in B$, 
we then obtain the equality ($*$). It remains to prove that $A$ is $Q$-absorbing. To that effect, we must just prove that $Q(x,C)=1$ for any $x\in A$. Set $D_x = \{y\in E : e^{i\langle  t,\xi(y)\rangle} w(y) = \lambda w(x)\}$. 
We know that $Q(x,D_x)=1$, and from $\lambda^{n+1} w(x) = 
\int_{D_x} e^{i\langle  t,\xi(y)\rangle} Q(t)^nw(y)\, Q(x,dy)\ \ $ ($n\geq 1$), we deduce that \\[0.15cm]
\indent $\ \ \ \  \ \ \ \ \ \ \ \ \ \ \ \ \ \ \ \  \ \ \ \ \ \ \ \ \ 
\lambda^n =  \int_{D_x} w(y)^{-1}\, Q(t)^nw(y)\, Q(x,dy)$. \\[0.15cm]
Since $Q(x,B)=1$, this equality holds also with $B$ instead of $D_x$. Besides, for any $y\in B$, we have 
$|Q(t)^nw(y)| \leq Q^n|w|(y) = \int_B |w(z)|Q^n(y,dz) = \pi(|w|)$, so that $|w(y)^{-1}\, Q(t)^nw(y)| \leq 1$. So, for some $D_{x,n}\in\cE$ such that $Q(x,D_{x,n})=1$, we have $Q(t)^nw(y) = \lambda^n\, w(y)$ for each $y\in D_{x,n}$. From $\cap_{n\geq 1}D_{x,n} \subset C$, one gets $Q(x,C)=1$ as claimed. \fdem

\noindent {\bf Lemma 12.2.} {\it Let $t\in\R^d$. If the equality ($*$) holds with $\lambda$, $A$ and $w$ as stated at the beginning of this section, then we have $\, r(Q(t)) \geq 1$. }

\noindent{\it Proof.} By integrating ($*$), one gets $Q(t)w=\lambda w$ on $A$, and since $A$ is $Q$-absorbing, this gives $Q(t)^nw=\lambda^n w$ on $A$ for all $n\geq 1$. Suppose $r(Q(t)) < 1$. Then 
$\lim_n Q(t)^nw = 0$ in $\cB$, and since $\cB\hookrightarrow\L^1(\pi)$, we have $\lim_n\pi(|Q(t)^nw|)=0$, but this is impossible because 
$|Q(t)^nw| = |w|$ on $A$, and by hypothesis $|w|$ is a nonzero constant on $A$ and $\pi(A)=1$. \fdem 

\noindent The previous lemmas show that, for any fixed $t\in\R^d$, we have 
$r(Q(t)) \geq 1$ iff the equality ($*$) holds for some $\lambda$, $A$ and $w$ as stated at the beginning of this section. Consequently, in order to prove the equivalence of Proposition 5.3, it remains to establish the following lemma whose proof is based on the use of the spectral results of \cite{keli}.  

\noindent{\bf Lemma 12.3.} {\it We have:  (S) $\Leftrightarrow\  \forall t\in\R^d,\, t\neq 0,\,  r(Q(t)) < 1$.} 

\noindent{\it Proof.} The direct implication is 
obvious. For the converse, let us consider a compact 
set $K_0$ 
in $\R^d\setminus\{0\}$. Let us first prove that \\[0.1cm] 
\indent $\displaystyle \ \ \ \ \ \ \ \ \ \ \ \ \ \ \ \ \ \ \ \ \ \ \ \ \  \ \ \ \  r_{K_0} = \sup \{r(Q(t)),\ t\in K_0\} < 1$. \\[0.1cm]
For that, let us assume that $r_{K_0}=1$. Then there exists a 
subsequence $(\tau_k)_k$ in $K_0$ such that we have 
$\lim_k r(Q(\tau_k)) = 1$. For $k\geq1$, let $\lambda_k$ be a spectral value of $Q(\tau_k)$ 
such that $|\lambda_k| = r(Q(\tau_k))$. By compactness, one may assume that the sequences $(\tau_k)_k$ and $(\lambda_k)_k$ converge. 
Let $\tau=\lim_k \tau_k$ and $\lambda=\lim_k\lambda_k$; observe that $\tau\in K_0$, thus $\tau\neq0$, and $|\lambda|=1$. 
Besides, by ($\widehat{K2}$) ($\widehat{K3}$)  ($\widehat{K4}$), the $Q(t)$'s satisfy the conditions 
of \cite{keli} near $\tau$. From \cite{keli} (p.~145), it follows that $\lambda$ is a spectral value of $Q(\tau)$, 
but this is impossible since, by hypothesis, $r(Q(\tau)) < 1$. This shows the claimed statement. \\
Let $\rho\in(r_{K_0},1)$. By applying \cite{keli} to $Q(\cdot)$ near any point $t_0\in K_0$, there exists 
a neighbourhood ${\cal O}_{t_0}$ of $t_0$ such that 
$\sup\{\|(z-Q(t))^{-1}\|_{\scriptsize \cB},\ t\in {\cal O}_{t_0},\ |z|=\rho\} < +\infty$. 
Since $K_0$ is compact, one gets 
$\sup\{\|(z-Q(t))^{-1}\|_{\scriptsize \cB}
,\ t\in K_0,\ |z|=\rho\} < +\infty$.  
Finally let $\Gamma$ be the oriented circle defined by $\{|z|=\rho\}$. 
Then the inequality stated in (S) follows from the 
following usual spectral formula \\[0.1cm]
\indent $\displaystyle \ \ \ \ \ \ \ \ \ \ \ \ \ \ \ \ \ \ \ \ \ \ \ \ \ \ \forall t\in K_0,\ \ 
Q(t)^n = \frac{1}{2i\pi}\int_{\Gamma} z^n (z-Q(t))^{-1}dz.$ 
\fdem 
%==================================
\subsection{Study of the set $G = \big\{t\in\R^d : r(Q(t))=1\big\}$} 
Here we still assume that Conditions ($\widehat{K}$) and (P) of Section 5.2 are fulfilled. We then know that (S) is equivalent to $G = \{0\}$ (Lem.~12.3). 
We assume moreover that the set of bounded elements
of $\cB$ is stable under complex conjugation and
under product. 
The next proposition specifies the statements of Proposition 5.4. 

\noindent{\bf Proposition 12.4.} {\it The set $G = \big\{t\in\R^d : r(Q(t))=1\big\}$ is a closed subgroup of  
$(\R^d,+)$. \\[0.12cm] 
Moreover, if the space $\widehat{\cB}$ of ($\widehat{K}$) verifies $\widehat{\cB}\hookrightarrow\L^1(\pi)$ and if Condition (CLT) of Section 5.1 holds, then $G$ is discrete,
and we have then the following properties. \\[0.12cm]
(i) If $G \neq \{0\}$, then there exist a point $a\in\R^d$, a closed  
subgroup $H$ in $\R^d$ of the form
$H=(vect\, G)^{\perp}\oplus\Delta$, where $\Delta$ is a discrete  
subgroup of $\R^d$, a $\pi$-full $Q$-absorbing set
$A\in\cE$, and a bounded measurable function $\theta\, :\, E\r\R^d$ such that
$$(**)\ \  \ \ \ \forall x\in A, \ \ \xi(y)+\theta(y)-\theta(x)\in  
a+H\ \ Q(x,dy)-a.s..$$
(ii) If ($**$) holds with  a $\pi$-full $Q$-absorbing set $A\in\cE$, a  
subgroup $H\neq\R^d$, and  a measurable function
$\theta\, :\, E\r\R^d$ such that $e^{i\langle t,\theta\rangle}\in\cB$  
for all $t\in\R^d$, then $G \neq \{0\}$.  }

\noindent {\it Proof.} Let $g_1,\, g_2\in G$, and for $k=1,2$, using Lemma 12.1, let  
$\lambda_k$,  $A_k$, and $w_k$ be the elements
associated with $g_k$ in ($*$). Then $A=A_1\cap A_2$ is a $\pi$-full  
$Q$-absorbing set, and $g_1-g_2$ satisfies ($*$)
with $A$, $\lambda = \lambda_1\overline{\lambda_2}$, and with $w=w_1\overline{w_2}\in\cB$. Thus $g_1-g_2\in G$ by Lemmas 12.1-2. Besides $0\in G$ since $Q1_E=1_E$. So $G$ is a  
subgroup of $(\R^d,+)$. To prove that $G$ is closed, let us consider any
sequence $(t_n)_n\in G^\N$ such that $\lim t_n = t$ in $\R^d$. By quasi-compactness (Lemma 12.1), each $Q(t_n)$
admits an eigenvalue, say $\lambda_n$, of modulus one. Now let  
$\lambda$ be a limit point of the sequence $(\lambda_n)_n$. Then
$|\lambda|=1$, and from \cite{keli} (p.~145), it follows that  
$\lambda$ is a spectral value of $Q(t)$, so $r(Q(t))\geq 1$, and $t\in G$ by Lemma 12.1. \\[0.15cm]
Now we assume that $\widehat{\cB}\hookrightarrow\L^1(\pi)$ (so ($\widetilde K$) of Sect.~4 is fulfilled) and that (CLT) holds. \\[0.15cm]
{\it $G$ is discrete.} From Lemma 5.2, we have $\lambda(t) = 1-\frac{1}{2}\langle\Gamma t,t\rangle + o(\|t\|^2)$ for  
$t$ near 0, where $\lambda(t)$ denotes the dominating eigenvalue  of $Q(t)$. Hence we  
have $r(Q(t)) = |\lambda(t)| < 1$ for $t$ near 0,
$t\neq0$. This proves that $0$ is an isolated point in $G$, hence $G$  
is discrete.  \\[0.1cm]
\noindent {\it Proof of (i).}
Set $G = \Z a_1\oplus\ldots\oplus\Z a_p$ with  $p\leq d$,
and let $\lambda_k$,  $A_k$, and $w_k$ be the elements associated with  
$a_k$ in ($*$).
Then  $A=\cap_{k=1}^p A_k$ is a $\pi$-full $Q$-absorbing set, and if  
$x\in A$ and $g=n_1a_1+\ldots+n_pa_p$ is any element of $G$,
we deduce from ($*$)
applied to each $a_k$, and by product that: \\[0.12cm]
\indent $\displaystyle\ \ \ \ \  \ \ \ \  \ \ \ \ \ \ \ \
\forall x\in A,\ \ e^{i\langle g,\xi(y)\rangle} \prod_{k=1}^p  
w_k(y)^{n_k} = \prod_{k=1}^p\lambda_k^{n_k}
\prod_{k=1}^pw_k(x)^{n_k}\ \ Q(x,dy)-a.s.$. \\[0.12cm]
Since $|w_k|$ is a nonzero constant function on $A$, one may assume  
without loss of generality that $|w_{k|A}| = 1_A$, so that there
exists a measurable function $\alpha_k : E\r [0,2\pi[$ such that we  
have, for all $z\in A$: $w_k(z) = e^{i\alpha_k(z)}$.
For $z\in A$, we set $V(z) = (\alpha_1(z),\ldots,\alpha_p(z))$ in $\R^p$.
Since the linear map $\chi : h \mapsto (\langle a_1,h  
\rangle,\ldots,\langle a_p,h \rangle)$ is clearly bijective
from $vect(G)$ into $\R^p$, one can define the element  
$\chi^{-1}(V(z))$ which satisfies
$\langle a_k, \chi^{-1}(V(z) \rangle =\alpha_k(z)$ for each $k=1,\ldots,p$.
Finally let  $\theta\, :\, E\r\R^d$ be a bounded measurable function   
such that $\theta(z) =  \chi^{-1}(V(z))$ for all $z\in A$.
Then we have $w_k(z) = e^{i\langle a_k,\theta(z) \rangle}$ for any  
$z\in A$ and $k=1,\ldots,p$. Consequently one gets
$\prod_{k=1}^pw_k(z)^{n_k} = e^{i\langle g,\theta(z) \rangle}$ for  
$z\in A$, and the above equality becomes, by setting
$\lambda_g = \prod_{k=1}^p\lambda_k^{n_k}$,
$$\forall x\in A,\ \ e^{i\langle g,\xi(y)+\theta(y)-\theta(x)\rangle}  
= \lambda_g\ \ Q(x,dy)-a.s..$$
For any $g\in G$, let us define $\beta_g\in\R$ such that
$\lambda_g = e^{i\beta_g}$, and for $x\in\R^d$, set $T_g(x) = \langle  
g,x \rangle$. The previous property yields \\[0.15cm]
\indent $\displaystyle\ \ \ \ \  \ \ \ \  \ \ \ \ \ \ \ \
\forall x\in A,\ \ \xi(y)+\theta(y)-\theta(x)\in\cap_{g\in G}\,  
T_g^{-1}(\beta_g+2\pi\Z)\ \ Q(x,dy)-a.s.$. \\[0.15cm]
Now let us define $\displaystyle H = \cap_{g\in G}\,  
T_g^{-1}(2\pi\Z)$. Then $H$ is a subgroup of $\R^d$, and the elements of
$\cap_{g\in G}\, (T_g^{-1}(\beta_g + 2\pi\Z))$ are in the same class  
modulo $H$. That is: \\[0.15cm]
\indent $\ \ \ \ \ \ \ \ \ \ \ \ \ \ \ \ \ \ \ \ \  \ \ \ \ \ \ \ \ \ \
\exists a\in\R^d,\ \cap_{g\in G}\, (T_g^{-1}(\beta_g + 2\pi\Z))\subset   
a+H$. \\[0.15cm]
This proves ($**$), and it remains to establish
that $H$ has the stated form. Actually, since $H$ is closed, $H$ is of  
the form $H=F \oplus\Delta$, where $F$
and $\Delta$ are respectively a subspace and a discrete subgroup in  
$\R^d$. So we have to prove that  $F=(vect\, G)^{\perp}$. \\
Let $x\in (vect\, G)^{\perp}$. Since $(vect\, G)^{\perp} = \cap_{g\in  
G}\, T_g^{-1}(\{0\})\subset H$, we have
$x = f + d$ for some $f\in F$, $d\in \Delta$, and for $\alpha\in\R$,  
the fact that $\alpha x\in  (vect\, G)^{\perp}\subset H$ yields
$\alpha x = f_\alpha + d_\alpha$ with some $f_\alpha\in F$ and  
$d_\alpha\in \Delta$. But we also have the unique decomposition
$\alpha x = \alpha f + \alpha d$ in $F\oplus\, vect\, \Delta$. Hence  
we have $\alpha d = d_\alpha\in\Delta $, and since
$\Delta$ is discrete and $\alpha$ can take any real value, we have  
necessary $d=0$. That is, $x\in F$. \\
Conversely, let $f\in F$ and let $g\in G$. Since $F\subset H$, we have  
$\langle g, f\rangle \in2\pi\Z$. Now let
$\alpha$ be any fixed nonzero irrational number. Since $\alpha f\in  
F\subset H$, we have
$ \alpha\, \langle g, f\rangle = \langle g, \alpha f\rangle  
\in2\pi\Z$. Hence $\langle g, f\rangle =0$. This gives
$f\in (vect\, G)^{\perp}$.

\noindent {\it Proof of (ii).} Let $t\in H^{\perp}$, $t\neq 0$. Then
$\langle t, \xi(y)\rangle + \langle t,\theta(y)\rangle  - \langle  
t,\theta(x)\rangle
=\langle t, a\rangle\ \ Q(x,dy)-$a.s. for all $x\in A$.  Setting  
$w(\cdot) = e^{i\langle t,\theta(\cdot)\rangle}$ and
$\lambda =  e^{i\langle t,a\rangle}$, this yields for all  $x\in A$ \\[0.12cm]
\indent $\ \ \ \ \ \ \ \ \ \ \ \ \ \ \ \ \ \ \ \ \  \ \ \ \ \ \ \ \ \ \
e^{i\langle t, \xi(y)\rangle}w(y) = \lambda w(x)\ \ Q(x,dy)-a.s.$. \\[0.15cm]
Since $w\in\cB$ by hypothesis, this gives ($*$), and Lemmas 12.1-2 implies that $t\in G$. \fdem

\newpage

\noindent{\bf Appendix A. Proof of Proposition 7.1.}

\noindent Proposition 7.1 will follow from the slightly more general Proposition below. The derivative arguments are presented here in the case $d=1$, but the  extension to $d\geq2$ is obvious by the use of the partial derivatives.  

\noindent Let $I$ be any subset of $\mathbb R$, 
let $T_0:I\rightarrow\mathbb R$ and $T_1:I\rightarrow\mathbb R$, 
let $({\cB }_\theta,\ \theta\in I)$ be a family of general 
Banach spaces. We shall write $\Vert\cdot\Vert_{\theta,\theta'}$ for $\Vert\cdot\Vert_{\scriptsize \cB_\theta,\cB_{\theta'}}$ 
and $\Vert\cdot\Vert_{\theta}$ for $\Vert\cdot\Vert_{\scriptsize\cB_\theta}$. Recall that we set 
$\cD_{\kappa} = \{z\in\C:\vert z\vert\ge\kappa,\ \vert z-1\vert\ge(1-\kappa)/2\}$ for any $\kappa\in(0,1)$. The notation $\cB_{\theta}\hookrightarrow\cB_{\theta'}$ means that $\cB_{\theta}\subset\cB_{\theta'}$ and that the identity map from $\cB_{\theta}$ into $\cB_{\theta'}$ is continuous. 

\noindent Let $\cB$ and $\widetilde{\cB}$ be some spaces of the previous family, and assume that $\cB_\theta\hookrightarrow\widetilde{\cB}$ for all $\theta\in I$. Finally let ${\cal U}$ be an open 
neighbourhood of $0$ in $\R^d$, and let $(Q(t),\ t\in{\cal U})$ be any family of operators in $\cL(\tilde\cB)$ such that we have $Q(t)_{|\cB_\theta}\in\cL(\cB_\theta)$ for all $t\in{\cal U}$ and $\theta\in I$. Let us introduce the following hypothesis. 

\noindent  {\bf Hypothesis $\cD(m)\ $} ($m\in\N^*$). {\it  
For all $\theta\in I$ there exists a neighbourhood
$\cV_\theta\subset{\cal U}$  of 
$0$ in $\R^d$ such that, for all $j=1,...,m$, 
we have: \\[0.12cm]
\noindent {\bf (0)} $\ [T_0(\theta)\in I \ 
\Rightarrow\ \ \cB_\theta\hookrightarrow \cB_{T_0(\theta)}]\ \ $ and
$\ \ [T_1(\theta)\in I \ 
\Rightarrow\ \ \cB_\theta\hookrightarrow\cB_{T_1(\theta)}]$ \\[0.12cm]
\noindent {\bf (1)} $\ \ T_0(\theta)\in I$ 
implies that 
$Q(\cdot)\in\cC^0(
\cV_\theta,\cB_{\theta},\cB_{T_0(\theta)})$ \\[0.12cm]
\noindent {\bf (2)} $\ \ \theta_j:=T_1(T_0T_1)^{j-1}(\theta)\in I$ 
implies that   
$Q(\cdot)\in \cC^j(\cV_\theta,\cB_{\theta},\cB_{\theta_j})$ \\[0.12cm]
\noindent {\bf (3')} There exists a real number 
$\kappa_\theta\in(0,1)$ such that, for all
$\kappa\in[\kappa_\theta,1)$, there exists
a neighbourhood
$\cV_{\theta,\kappa}\subseteq \cV_\theta$  of $0$ in $\R^d$ such that $R_z(t):=(z-Q(t))^{-1}\in\cL(\cB_\theta)$ for all 
$z\in\cD_{\kappa}$ and all $t\in\cV_{\theta,\kappa}$, and we have 
$$M_{\theta,\kappa} := \sup\big\{\|R_z(t)\|_{\theta},\ 
t\in \cV_{\theta,\kappa},\ z\in\cD_{\kappa}\big\} < +\infty$$ 
\noindent {\bf (4)} There exists 
$a\in\bigcap_{k=0}^m\left[T_0^{-1}(T_0T_1)^{-k}(I)
\cap (T_1T_0)^{-k}(I)\right]$ such that we have $\cB=\cB_a$ and $\widetilde{\cB} = \cB_{(T_0T_1)^mT_0(a)}$. }

\noindent When applied to the Fourier kernels, the above conditions (0) (1) (2) and (4) are exactly those of Hypothesis $\cC(m)$ in Section 7.1, and according to Theorem (K-L) of Section~4, Condition~(3') of $\cD(m)$ is implied by (3) of $\cC(m)$. Hence $\cC(m)$ implies that the Fourier kernels satisfy $\cD(m)$, so Proposition 7.1 follows from the next proposition. Let us notice that, from (4), we have \\[0.15cm]
\indent $\displaystyle \ \ \ \ \ \ \ \ \ \ \ \  \ \ \ \ 
\Theta_a = \{a,\, T_0a,\, T_1T_0a,\, T_0T_1T_0a,\, \ldots, (T_0T_1)^mT_0(a)\} \subset I$. \\[0.15cm]
Let us define $\tilde\kappa = \max_{\theta\in\Theta_a}\kappa_\theta\in(0,1)$, 
and $\widetilde{\cal O} = \bigcap_{\theta\in\Theta_a}\cV_{\theta,\tilde\kappa}$. 

\noindent {\bf Proposition A.} {\it Under 
Hypothesis $\cD(m)$, we have
$R_z(\cdot)\in\cC^m(\widetilde{\cal O},\cB,\tilde {\cB})$ 
for all $z\in{\cal D}_{\tilde\kappa}$, and for any compact subset ${\cal O}$ of $\widetilde{\cal O}$, we have ${\cal R}_{\ell} := \sup\{\|\, R_z^{(\ell)}(t)\|_{\scriptsize{\cB,\widetilde{\cB}}},\, z\in{\cal D}_{\tilde\kappa},\, t\in {\cal O}\, \} < +\infty$ for each $\ell=0,\ldots,m$. }

\noindent {\bf Remark. } {\it Let ${\cal O}$ be a compact subset of $\widetilde{\cal O}$. By Conditions (1) (2), we have for any $\theta\in\Theta_a$: \\[0.12cm]
$T_0(\theta)\in\Theta_a\ \Rightarrow\ {\cal Q}_{0,\theta} := \sup_{t\in{\cal O}}\|Q(t)\|_{\scriptsize{\cB_{\theta},\cB_{T_0(\theta)}}} < +\infty$ \\[0.1cm]
$\theta_j =T_1(T_0T_1)^{j-1}(\theta)\in\Theta_a\ \Rightarrow\ {\cal Q}_{j,\theta} := \sup_{t\in{\cal O}}\|Q^{(j)}(t)\|_{\scriptsize{\cB_{\theta},\cB_{\theta_j}}} < +\infty\ $ ($j=1,\ldots,m$). \\[0.1cm] 
The proof below shows that ${\cal R}_{\ell}$ in Proposition~A can be bounded by a polynomial expression involving the (finite) constants $\cM := \max_{\theta\in\Theta_a} M_{\theta,\tilde\kappa}$ and 
${\cal Q}_j := \max_{\theta\in\Theta_a\cap\tau_{j}^{-1}(\Theta_a)} {\cal Q}_{j,\theta}\ $ ($j=0,\ldots,\ell)$,
with $\tau_0 := T_0$, and $\tau_j:=T_1(T_0T_1)^{j-1}$ if $j\geq1$. }

\noindent The proof below involves the derivatives of some operator-valued maps defined as the composition of $Q(t)$ (or its derivatives) and  $R_z(t)$ (or its derivatives obtained by induction), where these operators are seen as elements of $\cL(\cB_{\theta_1},\cB_{\theta_2})$ and $\cL(\cB_{\theta_2},\cB_{\theta_3})$ for suitable $\theta_i\in I$. To that effect, it will be convenient to use the next notations. \\[0.15cm]
\noindent {\bf Notation.} {\it 
Let $\theta_1,\theta'_1\in I$.
An element of ${\mathcal L}(\theta_1,\theta'_1)$ 
is a family
$f=(f_z(t))_{z,t}$ of elements of ${\mathcal L}
({\mathcal B}_{\theta_1},{\mathcal B}_{\theta'_1})$ 
indexed by $(z,t)\in J$
(for some $J\subseteq {\mathbb C}\times{\mathbb R}^d$) satisfying the following condition: 
there exists $\hat\kappa_0\in(0,1)$ such that, 
for all $\kappa\in[\hat\kappa_0,1)$, there
exists a neighbourhood $\widetilde{\cU}_\kappa$ of $0$ 
in $\R^d$ such that $\cD_{\kappa}\times 
\widetilde{\cU}_\kappa\subseteq J$.

\noindent Let $\theta,\theta',\theta''\in I$. 
Given $V=(V_z(t))_{(z,t)\in J_U}\in 
{\mathcal L}(\theta,\theta')$ 
and $U=(U_z(t))_{(z,t)\in J_V}\in 
{\mathcal L}(\theta',\theta'')$, we define
$UV=\left(U_z(t)V_z(t)\right)_{(z,t)\in J_U\cap
 J_V}\in {\mathcal L}(\theta,\theta'')$.

\noindent Let $\ell\in\N$, let $\theta$ and $\theta_1$ in I be such that 
$\cB_{\theta} \hookrightarrow \cB_{\theta_1}$, and let $\theta'$ and $\theta'_1$ be in I. 
An element $f=(f_z(t))_{z,t}$ of ${\mathcal L}
(\theta_1,\theta'_1)$ is said to be in
${\mathcal C}^{\ell}(\theta,\theta')$
if the following condition holds: there exists $\hat\kappa\in(\hat\kappa_0,1)$ 
such that, for all $\kappa\in[\hat\kappa,1)$, 
there exists a neighbourhood 
${\mathcal U}_\kappa \subseteq \widetilde\cU_{\kappa}$ 
of $0$ in $\R^d$ such that, for all 
$z\in{\cal D}_\kappa$ and all $t\in \cU_\kappa$, we have 
$f_z(t)(\cB_\theta)\subseteq \cB_{\theta'}$, $\, f_z(\cdot)_{\vert \cB_\theta}
\in\cC^\ell(\cU_\kappa,\cB_\theta,\cB_{\theta'})$ and
$$\sup_{z\in{\cal D}_\kappa,\ t\in {\cal U}_\kappa,
\ j=0,...,\ell}\Vert f_z^{(j)}(t)\Vert_{\theta,\theta'} <+\infty.$$
When $f = (f_z(t))_{z,t}\in\cC^\ell(\theta,\theta')$, we set $f^{(\ell)} = (f^{(\ell)}_z(t))_{z,t}$. 
} 

\noindent Let us observe that (2) in $\cD(m)$ implies that $Q:=(Q(t))_{z,t}\in\cC^j(\theta,\theta_j)$ when we have $\theta\in I$ 
and $\theta_j:=T_1(T_0T_1)^{j-1}(\theta)\in I$. Now, we are in a position to state the next obvious (but important) facts (I)-(III), which will be repeatedly used in the proof of Proposition A.  Let $\theta_1$, $\theta_2$, $\theta_3$ and $\theta_4$ be in $I$.  

\noindent {\it 
(I) Assume that $\cB_{\theta_1}
\hookrightarrow\cB_{\theta_2}$, that
$\cB_{\theta_3}\hookrightarrow
\cB_{\theta_4}$ and that $V\in\cL(\theta_2,\theta_3)$. If $V\in\cC^k(\theta_2,\theta_3)$, then  
$V$ is in $\cC^k(\theta_2,\theta_4)$, in 
$\cC^k(\theta_1,\theta_3)$ and in
$\cC^k(\theta_1,\theta_4)$. 

\noindent (II) Assume that $V\in\cL(\theta_1,\theta_2)$ and $U\in\cL(\theta_2,\theta_3)$. If 
$V\in\cC^0(\theta_1,\theta_2)$ and $U\in\cC^0(\theta_2,\theta_3)$, then $UV\in\cC^0(\theta_1,\theta_3)$. 

\noindent (III) Let $U\in{\mathcal L}(\theta_3,
\theta_4)$ and $V\in{\mathcal L}(\theta_1,
\theta_2)$. Assume that 
$\cB_{\theta_1} \hookrightarrow \cB_{\theta_2}\hookrightarrow
\cB_{\theta_3}\hookrightarrow
\cB_{\theta_4}$, that 
$V\in\cC^0(\theta_1,\theta_2)\cap 
\cC^1(\theta_1,\theta_3) $ and that
$U\in\cC^1(\theta_2,\theta_4)\cap  \cC^0(
\theta_3,\theta_4) $. Then $UV$ is defined in ${\mathcal L}(\theta_1,\theta_4)$, and we have 
$UV\in\cC^1(\theta_1,\theta_4)$ and $(UV)' = U'V + UV'$. }

\noindent{\it Proof of Proposition A.} Lemmas A.1-2 below will be our basic statements. 

\noindent  {\bf Lemma A.1.} {\it If $\theta, T_0(\theta)\in I$, then $(R_z(t))_{z,t}\in{\cC}^0(\theta,T_0(\theta))$. }

\noindent{\it Proof.}  
Let $\kappa\in[\max(\kappa_\theta,
\kappa_{T_0(\theta)}),1)$.
Let $\cU_{\theta,\kappa}^{(0)}=\cV_{\theta,\kappa}
\cap \cV_{T_0(\theta),\kappa}$.
From the usual operator formula  
$Id-W^{n+1} = \sum_{k=0}^nW^k(Id-W)$, one
easily deduces the following equality, which is valid for any  bounded  
linear operators $S$ and $T$ on a Banach space
such that $S$ and $S-T$ are invertible: \\[0.1cm]
\indent $\displaystyle 
(*)\ \ \ \ \ \ \ \ \ \ \ \ \ \ 
(S-T)^{-1} = \sum_{k=0}^n (S^{-1}T)^kS^{-1} +  
(S^{-1}T)^{n+1}(S-T)^{-1}$. \\[0.1cm]
With $n=0$, $S=z-Q(t_0)$, $T=Q(t)-Q(t_0)$, thus $S-T = z-Q(t)$,  
Formula $(*)$ yields 
$$\forall z\in \cD_\kappa,\ \forall
t\in\cU_{\theta,\kappa}^{(0)},\ \   R_z(t) - R_z(t_0) =  
R_z(t_0)\, (Q(t)-Q(t_0))\, R_z(t).$$
Using the constants $M_{\theta,\kappa}$ and 
$M_{T_0(\theta),\kappa}$,
Condition (1) 
in $\cD(m)$ gives the desired property. \fdem

\noindent {\bf Lemma A.2.} {\it If 
$\theta,T_0(\theta),T_1T_0(\theta),T_0T_1T_0(\theta)\in I$, then we have  $(R_z(t))_{z,t}\in\cC^1\big(\theta,T_0T_1T_0(\theta)\big)$ and $R'=R\, Q'\, R$. }

\noindent{\it Proof.} 
Let us define $\theta_1=T_0(\theta)$, 
$\theta_2=T_1T_0(\theta)$,
$\theta_3=T_0T_1T_0(\theta)$
and $ \kappa_\theta^{(1)} =\max(\kappa_\theta,
\kappa_{\theta_1} ,\kappa_{\theta_2},
\kappa_{\theta_3})$. Let us consider a real number
$\kappa\in[\kappa_\theta^{(1)} ,1)$.
We define
$ \cU_{\theta,\kappa}^{(1)}
  =\cU_{\theta,\kappa}^{(0)}\cap\cU_
{\theta_2,\kappa}^{(0)}$.
Let $t_0,t\in \cU_{\theta,\kappa}^{(1)}$ and 
$z\in {\cal D}_{\kappa}$.
Formula $(*)$ with $n=1$, $S=z-Q(t_0)$, $T=Q(t)-Q(t_0)$ gives 
$$R_z(t) = R_z(t_0) \ +\  R_z(t_0)\, [Q(t)-Q(t_0)]\, R_z(t_0) \ +\  
\vartheta_z(t),$$
with $\vartheta_z(t) := R_z(t_0)\, [Q(t)-Q(t_0)]\,  
R_z(t_0)\, [Q(t)-Q(t_0)]\, R_z(t).$ But we have: 
$${\Vert \vartheta_z(t)\Vert_{\theta,\theta_3}\over \vert t-t_0\vert}
\le\Vert R_z(t_0) \Vert_{\theta_3}\Vert Q(t)-Q(t_0)\Vert
    _{\theta_2,\theta_3}\Vert R_z(t_0)\Vert_{\theta_2}
    {\Vert Q(t)-Q(t_0)\Vert_{{\theta_1},{\theta_2}}
             \over\vert t-t_0\vert}
     \Vert R_z(t)\Vert_{\theta,{\theta_1}}$$     
which goes to 0 as $t$ goes to $t_0$, 
uniformly in $z\in\cD_{\kappa}$
(according to condition (2) and with the use
of $M_{\theta_3,\kappa}$,  
$M_{\theta_2,\kappa}$ and 
$M_{\theta,\kappa}$). 
In the same way, we have: \\[0.15cm]
\indent $\left\Vert R_z(t_0)(Q(t)-Q(t_0))R_z(t_0)-(t-t_0)R_z(t_0)Q'(t_0)R_z(t_0)
\right\Vert_{\theta,\theta_3}  $
$$\le M_{\theta_2,\kappa}
 \Vert Q(t)-Q(t_0)-(t-t_0)Q'(t_0)\Vert_{\theta_1,\theta_2}
  M_{\theta,\kappa} = o(t-t_0). $$
\noindent This shows that $R_z'(t_0)=R_z(t_0)Q'(t_0)R_z(t_0)$ in $\cL(\cB_\theta,\cB_{\theta_3})$. 
Moreover, $(R_z(t))_{z,t}\in\cC^0(\theta,\theta_1)$, $(Q'(t))_{z,t}\in\cC^0(\theta_1,\theta_2)$, 
and $(R_z(t))_{z,t}\in\cC^0(\theta_2,\theta_3)$, 
therefore $(R_z'(t))_{z,t}\in\cC^0(\theta,\theta_3)$. \fdem

\noindent By Lemma A.1, the following assertion holds: \\[0.15cm]
\noindent $(H_0)$ {\it If $\theta\in I$ and if
$T_0(\theta)\in I$, then $R=(R_z(t))_{z,t}\in\cC^0\left(\theta,T_0(\theta)\right)$. } \\[0.2cm]
\noindent For $\ell=1,\ldots,m$, let us set 
$${\cal E}_\ell = \big\{(i,j,k)\in\Z^3 : i\ge 0,j\ge 0, k\ge 0,\ i+j+k=\ell-1\big\},$$
and let us denote by ($H_\ell$) the following assertion: \\[0.1cm]
\noindent ($H_\ell$) {\it If $\displaystyle \, \theta\in \bigcap_{k=0}^\ell \left[T_0^{-1}(T_0T_1)^{-k}(I)\cap 
    (T_1T_0)^{-k}(I)\right]$, then $R=(R_z(t))_{z,t}\in \cC^\ell\left({\theta},
{(T_0T_1)^\ell T_0(\theta)}\right)$ and $\displaystyle R^{(\ell)}=\sum_{(i,j,k)\in {\cal E}_\ell}R^{(i)}
Q^{(1+j)}R^{(k)}$.} \\[0.15cm]
\noindent We want to prove $(H_m)$ by induction. By Lemma A.2, ($H_1$) holds. 

\noindent {\bf Lemma A.3.} {\it Let $1\leq \ell\leq m-1$. If $(H_0), (H_1),\ldots,(H_\ell)$ hold, then we have ($H_{\ell+1}$). }

\noindent{\it Proof.} Let $\theta\in \bigcap_{k=0}^{\ell+1}\left[{T_0}^{-1}(T_0T_1)^{-k}(I)
\cap (T_1T_0)^{-k}(I)\right]$. From $\cB_{(T_0T_1)^\ell T_0(\theta)}\hookrightarrow
\cB_{(T_0T_1)^{\ell+1}T_0(\theta)}$ and ($H_\ell$), we have 
$$R=(R_z(t))_{z,t}\in\cC^\ell\left(\theta,(T_0T_1)^{\ell+1}T_0(\theta)\right)\ \ \mbox{and}\ \  
R^{(\ell)}=\sum_{(i,j,k)\in {\cal E}_\ell}R^{(i)}
Q^{(1+j)}R^{(k)}.$$ 
Let $(i,j,k)\in{\cal E}_\ell$. We have to prove that $R^{(i)}Q^{(1+j)}R^{(k)}\in 
\cC^1\left(\theta,(T_0T_1)^{\ell+1} T_0(\theta)\right)$ 
and that: \\[0.15cm]
\indent $\displaystyle \ \ \ \ \ \ \ \ \ \ \ \ 
\left(R^{(i)}Q^{(1+j)}R^{(k)}\right)'=
R^{(i+1)}Q^{(1+j)}R^{(k)}
   +R^{(i)}Q^{(2+j)}R^{(k)}+R^{(i)}Q^{(1+j)}R^{(k+1)}$. \\[0.15cm]
Since $1\le k+1\le \ell$ and by induction hypothesis, we have:  
$$R^{(k)}\in\cC^1\left(\theta,
    {(T_0T_1)^{(k+1)}T_0(\theta)}\right)\ \ \mbox{and}\ \
 R^{(k)}\in\cC^0\left(\theta,
    {(T_0T_1)^{k}T_0(\theta)}\right).$$
Moreover, since $2+j\leq \ell+1\leq m$ and according to $\cD(m)$, we have~:
$$ Q^{(1+j)}\in\cC^0\left({(T_0T_1)^{(k+1)}T_0(\theta)},{T_1(T_0T_1)^{(k+j+1)}T_0(\theta)}\right)$$
$$\mbox{and}\ \ \ \ \ \ \ \ \ \ Q^{(1+j)}\in\cC^1\left({(T_0T_1)^{k}T_0(\theta)},
       {T_1(T_0T_1)^{k+j+1}T_0(\theta)}\right).$$
From Property (III), we then deduce that we have 
$Q^{(1+j)}R^{(k)}\in \cC^1\left(\theta,T_1(T_0T_1)^{k+j+1} T_0(\theta)\right)$
and $( Q^{(1+j)}R^{(k)} )'=Q^{(2+j)}R^{(k)} 
         +Q^{(1+j)}R^{(k+1)} $. \\
\noindent Analogously we have 
$Q^{(1+j)}R^{(k)}\in \cC^0\left(\theta,T_1(T_0T_1)^{k+j}\, T_0(\theta)\right)$, 
and, since $i+1\le\ell$, we have~:
$$R^{(i)}\in \cC^0\left(T_1(T_0T_1)^{k+j+1}\, T_0(\theta),(T_0T_1)^{k+j+i+2}\, T_0(\theta)\right)$$
$$\mbox{and}\ \ \ \ \ \ \ \ R^{(i)}\in \cC^1\left(T_1(T_0T_1)^{k+j}\, T_0(\theta),
  (T_0T_1)^{k+j+i+2}\, T_0(\theta)\right).$$
Since $k+j+i+2 = \ell+1$, this gives the desired property. \fdem 

\noindent Since, by hypothesis, $a\in\bigcap_{k=0}^m
\left[T_0^{-1}(T_0T_1)^{-k}(I)\cap (T_1T_0)^{-k}(I)\right]$, the properties $(H_0)$,...,$(H_m)$ show that the conclusions of Proposition A are  valid. 
More exactly, the previous induction proves that the neighbourhood $\widetilde{\cal O}$ of $t=0$ and the real number $\tilde\kappa$ may be defined as stated before Proposition A, and that for any compact subset ${\cal O}\subset\widetilde{\cal O}$, the constants ${\cal R}_{\ell}$ are bounded as indicated in the remark following Proposition~A. \fdem

\newpage

\noindent{\bf Appendix B. Proof of Propositions 11.4-8.}

\noindent {\bf B.0. Notations.}
For convenience we present the proofs of Propositions 11.4-8 in the  
case $d=1$. The extension to $d\geq2$ is straightforward
for Proposition 11.4,5,6,8 (just replace the inequality $|t\xi(x)|  
\leq |t|\, |\xi(x)|$ with the Schwarz inequality
$|\langle t,\xi(x)\rangle| \leq \|t\|\, \|\xi(x)\|$). It is easy for
Proposition 11.7 by considering partial derivatives.

\noindent We set $\Theta x = F(x,\theta_1)$. So $\Theta$ is a random  
Lipschitz transformation on $E$, and
the transition probability $Q$ can be expressed as:  $Qf(x) =  
\E[f(\Theta x)]$.

\noindent For any $\lambda\in(0,1]$, we set
$p_\lambda(x) = 1+\lambda\, d(x,x_0)$. For any $0<\alpha\leq 1$, $0<\beta\leq\gamma$, and $(x,y)\in E^2$, let us set 
$$\Delta^{(\lambda)}_{\alpha,\beta,\gamma}(x,y) = p_\lambda(x)^{\alpha\gamma}  p_\lambda(y)^{\alpha\beta} + p_\lambda(x)^{\alpha\beta}  p_\lambda(y)^{\alpha\gamma}.$$
Then the space  
$\cB_{\alpha,\beta,\gamma}$ defined in Section 11 is unchanged when
$m_{\alpha,\beta,\gamma}(f)$ is replaced with
$$m^{(\lambda)}_{\alpha,\beta,\gamma}(f) =  
\sup\bigg\{\frac{|f(x)-f(y)|}{d(x,y)^\alpha\, \Delta^{(\lambda)}_{\alpha,\beta,\gamma}(x,y)},\ x,y\in E,\  
x\neq y\bigg\},$$
and for any $f\in \cB_{\alpha,\beta,\gamma}$, the following quantity \\
\indent $\displaystyle\ \ \ \ \ \ \ \  \ \ \ \ \ \ \ \ \ \  \ \ \ \ \  
\ \ \ \ \ \ \ \ \ \ \ \ \ \
|f|^{(\lambda)}_{\alpha,\gamma}  = \sup_{x\in E}\  
\frac{|f(x)|}{p_\lambda(x)^{\alpha(\gamma+1)}}$, \\[0.15cm]
is finite. The resulting new norm $\displaystyle \|f\|^{(\lambda)}_{\alpha,\beta,\gamma} =  
m^{(\lambda)}_{\alpha,\beta,\gamma}(f) + |f|^{(\lambda)}_{\alpha,\gamma}$
is equivalent to the norm $\|\cdot\|_{\alpha,\beta,\gamma}$ defined in Section 11. Consequently, for  
$(\alpha,\beta,\gamma)$ fixed as above, Propositions 11.4-8 can be established by considering on $\cB_{\alpha,\beta,\gamma}$ the norm $\|f\|^{(\lambda)}_{\alpha,\beta,\gamma}$ (for some value $\lambda\in(0,1]$).  
In most of the next estimates, we shall assume 
$\lambda=1$; the possibility of choosing suitable small $\lambda$ will  
occur in the proof of the Doeblin-Fortet inequaliies (in Prop.~11.4 and Prop.~11.8). Anyway, this already appears in the proof  
of Proposition 11.2, see \cite{aap}.

\noindent Let $\cC_\lambda = \max\{\cC,1\}   + \lambda\,  d(\Theta x_0  
, x_0)$.
In the sequel, we shall use repeatedly the fact that  
$p_\lambda(\cdot)$ and $p(\cdot)$ are equivalent functions,
and that (see \cite{aap} p.~1945)
$$\sup_{x\in E} \frac{p_\lambda(\Theta x)}{p_\lambda(x)} \leq  
\cC_\lambda\leq \cM,$$
from which we deduce that
$$\Delta^{(\lambda)}_{\alpha,\beta,\gamma}(\Theta x,\Theta y) 
\leq \cC_\lambda^{\alpha(\gamma+\beta)}\, \Delta^{(\lambda)}_{\alpha,\beta,\gamma}(x,y).$$
We shall also use the fact that 
$$d(y,x_0) \leq d(x,x_0)\ \Rightarrow\ \Delta^{(\lambda)}_{\alpha,\beta,\gamma}(x,y) \leq 2\, p_\lambda(x)^{\alpha\gamma}  p_\lambda(y)^{\alpha\beta}.$$
Indeed, if $d(y,x_0) \leq d(x,x_0)$, then we have $p_\lambda(y) \leq p_\lambda(x)$, so that $$p_\lambda(x)^{\alpha\beta}  p_\lambda(y)^{\alpha\gamma} = p_\lambda(x)^{\alpha\beta}\, p_\lambda(y)^{\alpha(\gamma-\beta)}\,  p_\lambda(y)^{\alpha\beta}\leq  p_\lambda(x)^{\alpha\beta}\, p_\lambda(x)^{\alpha(\gamma-\beta)}\,  p_\lambda(y)^{\alpha\beta},$$
thus $p_\lambda(x)^{\alpha\beta}  p_\lambda(y)^{\alpha\gamma} \leq p_\lambda(x)^{\alpha\gamma}  p_\lambda(y)^{\alpha\beta}$. 

\noindent {\bf B.1. A preliminary lemma.} The proofs of Propositions 11.4-8 are based on the following lemma.

\noindent {\bf Lemma B.1.} {\it Let $q : E\r\C$ measurable, $Kf(x) =  
\E[q(\Theta x)\, f(\Theta x)\, ]$, and
let $\lambda\in(0,1]$. Suppose that there exist constants
$a,\, A,\, b,\, B$ such that we have for all $x,y\in E$ satisfying
$d(y,x_0) \leq d(x,x_0)$ \\[0.15cm]
\indent $\displaystyle\ \ \ \ \ \ \ \  \ \ \ \ \ \ \ \ \ \ $ {\bf (i)}
$\ \ \ |q(x)|\leq A\, p_\lambda(x)^a\ \ \ $; $\ \ \ ${\bf (ii)}
$\ \ \ |q(x)-q(y)| \leq B\, d(x,y)^\alpha\, p_\lambda(x)^b$.  \\[0.15cm]
Then we have for $f\in\cB_{\alpha,\beta,\gamma}$ and $x,y$ as above stated }
\begin{eqnarray*}
|Kf(x)| &\leq& A\, |f|^{(\lambda)}_{\alpha,\gamma}\,  
p_\lambda(x)^{a+\alpha(\gamma+1)}\, \E[\cM^{a+\alpha(\gamma+1)}]  
\\[0.2cm]
|Kf(x)-Kf(y)| &\leq& A\, m^{(\lambda)}_{\alpha,\beta,\gamma}(f)\, d(x,y)^\alpha\,
p_\lambda(x)^a\, \Delta^{(\lambda)}_{\alpha,\beta,\gamma}(x,y)\, 
\E[\cC^\alpha\, \cC_\lambda^{a+\alpha(\gamma+\beta)}] \\[0.15cm]
&\ & \ \ \ \ \ + \ B\, |f|^{(\lambda)}_{\alpha,\gamma}\, d(x,y)^\alpha\,  
p_\lambda(x)^{b}\,  p_\lambda(y)^{\alpha(\gamma+1)}\,
\E[\cC^\alpha\, \cM^{b+\alpha(\gamma+1)}].
\end{eqnarray*}
\noindent{\it Proof.} We have
\begin{eqnarray*}
|Kf(x)| \leq \E[|q(\Theta x)\, f(\Theta x)|] &\leq&
A\,|f|^{(\lambda)}_{\alpha,\gamma}\, \E[\, p_\lambda(\Theta x)^a\,  
p_\lambda(\Theta x)^{\alpha(\gamma+1)}\, ] \\
&\leq&  A\,|f|^{(\lambda)}_{\alpha,\gamma}\,   
p_\lambda(x)^{a+\alpha(\gamma+1)}\E[\cM^{a+\alpha(\gamma+1)}\, ].
\end{eqnarray*}
Moreover, for $x,y\in E$ satisfying $d(y,x_0) \leq d(x,x_0)$ (thus  
$p_\lambda(y)\leq p_\lambda(x)$), we have
\begin{eqnarray*}
|Kf(x) - Kf(y)| &\leq& \E\bigg[|q(\Theta x)|\, |f(\Theta x)-f(\Theta  
y)|\bigg] \ +\
\E\bigg[ |f(\Theta y)|\, |q(\Theta x) - q(\Theta y)|\bigg] \\
&\leq& A\, m^{(\lambda)}_{\alpha,\beta,\gamma}(f)\, 
\E\bigg[p_\lambda(\Theta x)^a\, d(\Theta x,\Theta y)^\alpha\, 
\Delta^{(\lambda)}_{\alpha,\beta,\gamma}(\Theta x,\Theta y)\bigg] \\
&\ & \ \ \ \ \ +\  |f|^{(\lambda)}_{\alpha,\gamma}\, B\, \E\bigg[p(\Theta  
y)^{\alpha(\gamma+1)}\, d(\Theta x,\Theta y)^\alpha\,
p_\lambda(\Theta x)^b\bigg]\\
&\leq& A\, m^{(\lambda)}_{\alpha,\beta,\gamma}(f)\,  d(x,y)^\alpha\,   
p_\lambda(x)^a\, \Delta^{(\lambda)}_{\alpha,\beta,\gamma}(x,y)\, 
\E\bigg[\cC^\alpha\, \cC_\lambda^{a+\alpha(\gamma+\beta)}\bigg] \\
&\ & \ \ \ \ \ +\ B\,  |f|^{(\lambda)}_{\alpha,\gamma}\, d(x,y)^\alpha \,    
p_\lambda(x)^{b}\,  p_\lambda(y)^{\alpha(\gamma+1)}\,
\E\bigg[\cC^\alpha\, \cM^{b+\alpha(\gamma+1)}\bigg]. 
\end{eqnarray*}
Lemma B.1 is then proved. \fdem 

\noindent For the use of Lemma B.1, it is worth noticing that the  
supremum bound defining the
H\"older constants $m_{\alpha,\beta,\gamma}(f)$ or  
$m^{(\lambda)}_{\alpha,\beta,\gamma}(f)$ can be obviously computed over
the elements $x,y\in E$ such that $d(y,x_0) \leq d(x,x_0)$.
Lemma B.1 will be applied below with $q(\cdot)$ depending on the  
function $\xi$.
Remember that $\xi$ verifies the following  hypothesis:  \\[0.15cm]
\noindent $\displaystyle {\bf (L)_s} \ \ \ \ \ \ \ \ \ \ \ \ \
\forall (x,y)\in E\times E,\ \ |\xi(x)-\xi(y)| \leq S\, d(x,y)\,  
[1+d(x,x_0)+d(y,x_0)]^s$. \\[0.15cm]
 From $\bf (L)_s$, it follows that there exists $C>0$ such that we have  
for $x\in E$ \\[0.15cm]
\indent $\displaystyle\ \ \ \ \ \ \ \  \ \ \ \ \ \ \ \ \ \ \  \ \ \ \  
\ \ \ \ \ \  \ \ \ \ \ \ \ \ \ \
|\xi(x)| \leq C\, p(x)^{s+1},$ \\[0.15cm]
and for $x,y\in E$ satisfying $d(y,x_0) \leq d(x,x_0)$:  \\[0.15cm]
\indent $\displaystyle\ \ \ \ \ \ \ \
|\xi(x)-\xi(y)| \leq C\, d(x,y)\, p(x)^s\ $ and $\ |\xi(x)-\xi(y)|  
\leq C\, d(x,y)^\alpha\, p(x)^{s+1-\alpha}$.

\noindent {\bf B.2. Proof of Proposition 11.4}. {\it This proposition  
states that (K) of Section 4 holds w.r.t.~the space $\cB_{\alpha,\beta,\gamma}$
if we have $s+1\leq\beta\leq\gamma$ and }
$$I = \E[\, \cM^{\alpha(\gamma+1)} + \cC^\alpha\, \cM^{\alpha(\gamma+\beta)}\, ] <  +\infty \ \ \ \mbox{and}\ \ \ 
\E[\cC^\alpha\, \max\{\cC,1\}^{\alpha(\gamma+\beta)}] < 1.$$

\noindent The strong ergodicity condition (K1) of Section 1 holds by Proposition  11.2. Besides we have for $f\in\cB_{\alpha,\beta,\gamma}$ 
$$\pi\big(|e^{i\langle t,\, \xi\rangle}-1|\, |f|\big)
\leq |f|_{\alpha,\gamma}\, \pi\big(|e^{i\langle t,\, \xi\rangle}-1|\,  
p^{\alpha(\gamma+1)}\big).$$
Since $\pi(p^{\alpha(\gamma+1)}) < +\infty$ (Prop.~11.1), the continuity condition of (K) is satisfied: in fact, from Lebesgue's theorem and Remark (a) of Section 4, we have ($\widehat{K2}$) of Section 5.2. To study the  Doeblin-Fortet inequalities of (K), notice that
$Q(t) = K$ where $K$ is associated to $q(x) = e^{i t\xi(x)}$ with the  
notations of Lemma B.1.
By using $\bf (L)_s$ and the inequality $|e^{iT}-1|\leq 2|T|^\alpha$, one  
easily gets (i)-(ii) in Lemma B.1 with $A=1$, $a=0$ and 
$B = D_\lambda\, |t|^\alpha\ $, $b=\alpha s$, where $D_\lambda$ is
a positive constant resulting from  $\bf (L)_s$ and the equivalence between  
$p_\lambda(\cdot)$ and $p(\cdot)$.
Then, from Lemma B.1, we have for any $f\in\cB_{\alpha,\beta,\gamma}$ 
$$|Q(t)f|^{(\lambda)}_{\alpha,\gamma} \leq \E[\, \cM^{\alpha(\gamma+1)}\,]\, |f|^{(\lambda)}_{\alpha,\gamma} \leq I\, |f|^{(\lambda)}_{\alpha,\gamma}$$
and for $x,y\in E$ such that $d(y,x_0) \leq d(x,x_0)$
\begin{eqnarray*}
|Q(t)f(x)-Q(t)f(y)| &\leq& 
m^{(\lambda)}_{\alpha,\beta,\gamma}(f)\, d(x,y)^\alpha\,  \Delta^{(\lambda)}_{\alpha,\beta,\gamma}(x,y) \, \E[\cC^\alpha\, \cC_\lambda^{\alpha(\gamma+\beta)}] \\
&+& \, D_\lambda\, |t|^\alpha \, |f|^{(\lambda)}_{\alpha,\gamma}\,  d(x,y)^\alpha\, 
p_\lambda(x)^{\alpha s}\,   p_\lambda(y)^{\alpha(\gamma+1)}\, \E[\cC^\alpha\, \cM^{\alpha(\gamma+s+1)}].
\end{eqnarray*}
Since $p_\lambda(x)^{\alpha s}\,   p_\lambda(y)^{\alpha(\gamma+1)} \leq 
p_\lambda(x)^{\alpha(s+1)}\,   p_\lambda(y)^{\alpha\gamma} \leq p_\lambda(x)^{\alpha\beta}\, 
p_\lambda(y)^{\alpha\gamma} \leq \Delta^{(\lambda)}_{\alpha,\beta,\gamma}(x,y)$, the previous inequalities prove that $Q(t)$ continuously acts on  
$\cB_{\alpha,\beta,\gamma}$, and setting $E_\lambda =  I\, D_\lambda$, that 
$$m^{(\lambda)}_{\alpha,\beta,\gamma}(Q(t)f) \leq \E[\cC^\alpha\,  
\cC_\lambda^{\alpha(\gamma+\beta)}]\, m^{(\lambda)}_{\alpha,\beta,\gamma}(f)
+ E_\lambda\, |t|^\alpha\, |f|^{(\lambda)}_{\alpha,\gamma}.$$
Now, using the fact that the norms $\|f\|^{(\lambda)}_{\alpha,\beta,\gamma}$  
and $\|f\| = m^{(\lambda)}_{\alpha,\beta,\gamma}(f) + \pi(|f|)$
are equivalent (see \cite{aap} Prop.~5.2),
one obtains with some new constant $E_\lambda'$~:
\begin{eqnarray*}
m^{(\lambda)}_{\alpha,\beta,\gamma}(Q(t)f) &\leq&
\E[\cC^\alpha\, \cC_\lambda^{\alpha(\gamma+\beta)}]\,
m^{(\lambda)}_{\alpha,\beta,\gamma}(f)
+ E_\lambda'\, |t|^\alpha\,
\bigg(m^{(\lambda)}_{\alpha,\beta,\gamma}(f) +
\pi(|f|)\bigg)\\
  &\leq& \bigg(\E[\cC^\alpha\, \cC_\lambda^{\alpha(\gamma+\beta)}] +  
E_\lambda'\, |t|^\alpha\bigg)\, m^{(\lambda)}_{\alpha,\beta,\gamma}(f)
+ E_\lambda'\, |t|^\alpha\, \pi(|f|).
\end{eqnarray*}
Since $\cC_\lambda\leq \cM$ and $\cC_\lambda \r \max\{\cC,1\}$ when  
$\lambda\r0$, it follows from Lebesgue theorem
that one can choose $\lambda$ such that $\E[\cC^\alpha\,  
\cC_\lambda^{\alpha(\gamma+\beta)}] < 1$. Now let $\tau>0$ such that
$$\kappa : = \E[\cC^\alpha\, \cC_\lambda^{\alpha(\gamma+\beta)}] +  
E_\lambda'\, \tau^\alpha <1.$$ 
Then, if $|t|\leq\tau$, we have 
$$m^{(\lambda)}_{\alpha,\beta,\gamma}(Q(t)f)\leq \kappa\,  
m^{(\lambda)}_{\alpha,\beta,\gamma}(f)+ E_\lambda'\, \tau^\alpha\,  
\pi(|f|).$$
Since $\pi(|Q(t)f|) \leq \pi(|Qf|) = \pi(|f|)$, this gives $\|Q(t)f\|  
\leq  \kappa\, \|f\| + (1+E_\lambda'\, \tau^\alpha)\, \pi(|f|)$,
and this easily leads to the  Doeblin-Fortet inequalities of (K), with ${\cal O} = (-\tau,\tau)$.  \fdem

\noindent {\it In the next proofs, except for Proposition 11.8, the  
technical parameter $\lambda$ used above
will be neglected, namely we shall assume $\lambda=1$, and the  
effective computation of
the constants occurring in the proofs will be of no relevance. So, to  
simplify the next estimates,
we shall still denote by  $C$ the constant in the above inequalities  
resulting from $\bf (L)_s$,
even if it is slightly altered through the computations (the effective  
constants will actually depend on parameters
as $\alpha$, $t_0\in\R$ fixed, $k\in\N$ fixed, $s$, $S$ ...). }

\noindent Proposition 11.5 will follow from Lemma B.4 with $k=0$.

\noindent {\bf B.3. Proof of Proposition 11.6}. {\it 
Actually let us prove that  
$\|Q(t)-Q\|_{_{\cB_{\alpha,\beta,\gamma},\cB_{\alpha,\beta,\gamma'}}} = O(|t|)$ if 
$\, 0<\beta\leq\gamma$, $\, \gamma' \geq \gamma + \frac{s+1}{\alpha}$, and }
$\displaystyle I=\E\left[\,  \cM^{s+1+\alpha(\gamma+1)} + \cC^\alpha\, \cM^{s+1+\alpha(\gamma+\beta)} \,\right] <+\infty$. 

\noindent Let $K=Q(t)-Q(0)$. Then $K$ is associated to
$q(x)= e^{it\xi(x)} -  1$. Using $\bf (L)_s$ and the inequality  
$|e^{iT}-1|\leq |T|$, one easily gets (i)-(ii) in Lemma B.1
with $A= C\, |t|$, $a=s+1$, and $B = C\,|t|$ and $b= s+1 - \alpha$. So
$$|Kf(x)| \leq C\, |t|\, |f|_{\alpha,\gamma}\,  
p(x)^{s+1+\alpha(\gamma+1)}\, \E[\,  \cM^{s+1+\alpha(\gamma+1)}\, ] \leq 
I\, C\, |t|\, |f|_{\alpha,\gamma}\,  
p(x)^{\alpha(\gamma'+1)},$$
and, by using the fact that $p(y)\leq p(x)$ (thus $\Delta_{\alpha,\beta,\gamma}(x,y) \leq 2\, p(x)^{\alpha\gamma}  p(y)^{\alpha\beta}$)
\begin{eqnarray*}
|Kf(x)-Kf(y)| &\leq& C\, |t|\, m_{\alpha,\beta,\gamma}(f)\, d(x,y)^\alpha\, 
p(x)^{s+1}\, 2\, p(x)^{\alpha\gamma}  p(y)^{\alpha\beta}\, \E[\cC^\alpha\, \cC_\lambda^{s+1+\alpha(\gamma+\beta)}] 
\, \\[0.15cm]
&\, & + \ C\, |t|\, |f|_{\alpha,\gamma}\, d(x,y)^\alpha\,  
p(x)^{s+1-\alpha}\,p(y)^{\alpha(\gamma+1)}\, \E[\cC^\alpha\, \cM^{s+1 - \alpha+\alpha(\gamma+1)}]. 
\end{eqnarray*}
Since $p(x)^{s+1+\alpha\gamma}  p(y)^{\alpha\beta} \leq p(x)^{\alpha\gamma'} p(y)^{\alpha\beta} \leq \Delta_{\alpha,\beta,\gamma'}(x,y)$ and 
$$p(x)^{s+1-\alpha}\,p(y)^{\alpha(\gamma+1)} \leq p(x)^{s+1+\alpha\gamma} \leq p(x)^{\alpha\gamma'} \leq \Delta_{\alpha,\beta,\gamma'}(x,y)$$ 
it follows that $\displaystyle |Kf(x)-Kf(y)| \leq 2\, I\, C\, |t|\,  \|f\|_{\alpha,\beta,\gamma} \,  
d(x,y)^\alpha\, \Delta_{\alpha,\beta,\gamma'}(x,y)$. \fdem

\noindent {\bf B.4. Proof of Proposition 11.7}. {\it This proposition  
states that $\cC(m)$
($m\in\N^*$) holds with $\cB = \cB_{\alpha,\beta,\gamma}$ and $\widetilde{\cB} =   
\cB_{\alpha,\beta,\gamma'}$ if we have $s+1\leq\beta\leq\gamma$, $\, \gamma' > \gamma + \frac{m(s+1)}{\alpha}$, and }
$$\E[\, \cM^{\alpha(\gamma'+1)}\, +\, \cC^\alpha\cM^{\alpha(\gamma'+ \beta)}\,] <+\infty\ \ \ \ \ \ 
\E[\, \cC^\alpha\, \max\{\cC,1\}^{\alpha(\gamma'+\beta)}] < 1.$$
\noindent Let $k\in\N$. Let us recall that we set 
$Q_k(t)(x,dy) =  
i^k\xi(y)^ke^{it\xi(y)}\, Q(x,dy)\ $ ($x\in E$, $t\in \R$).
For $u\in\R$, we set $e^{iu\xi(\cdot)} = e_u(\cdot)$. 

\noindent {\bf Lemma B.4.} {\it For $k\in\N$, we have  
$Q_k\in\cC^0(\R,\cB_{\alpha,\beta,\gamma},\cB_{\alpha,\beta,\gamma'})$ under the following conditions:  
$s+1\leq\beta\leq\gamma$, $\ \gamma' > \gamma + \frac{(s+1)k}{\alpha}$, and
$I = \E[\, \cM^{\alpha(\gamma'+1)}\, +\,  
\cC^\alpha\,\cM^{\alpha(\gamma' + \beta)}\, ] < +\infty$.}

\noindent{\it Proof.} Let $t,t_0\in\R$, $h=t-t_0$. We suppose that $|h|\leq1$. Let $K=Q_k(t)-Q_k(t_0)$. Then $K$ is associated to
$q(x)=(i\xi(x))^k\, \big(e_t(x) -  e_{t_0}(x)\big)$. Let  
$0<\varepsilon<\alpha$. Using the inequality
$|e^{iT}-1|\leq 2|T|^\varepsilon$, one gets (i) in Lemma B.1 with $A=C\,  
|h|^\varepsilon$ and $a=(s+1)(k+\varepsilon)$. Using also
$|e^{iT}-1|\leq 2|T|^\alpha$, we have for $k\geq1$ and
for $x,y\in E$ such that $d(y,x_0) \leq d(x,x_0)$ (thus $p(y)\leq p(x)$)~:
\begin{eqnarray*}
|q(x)-q(y)| &\leq& |\xi(x)^k-\xi(y)^k|\, |e_t(x) -  e_{t_0}(x)| +  
|\xi(y)|^k\, \left|\, \big(e_t(x) -  e_{t_0}(x)\big) -
\big(e_t(y) -  e_{t_0}(y)\big)\, \right| \\
&\leq& C\, |\xi(x)-\xi(y)|\, p(x)^{(s+1)(k-1)}\, |h|^\varepsilon\,  
p(x)^{(s+1)\varepsilon}\\
&\ & \ \ +\ C\,  p(x)^{(s+1)k}\bigg(|e_h(x) - e_h(y)| + |e_h(y) - 1|\,  
|e_{t_0}(x) -  e_{t_0}(y)|\bigg) \\
&\leq&  C\,  |h|^\varepsilon\, d(x,y)^\alpha\, p(x)^{s+1-\alpha}\,  
p(x)^{(s+1)(k-1+\varepsilon)} \\
&\ & \ \ +\  C\, p(x)^{(s+1)k}\,  \bigg(|h|^\alpha\,  d(x,y)^\alpha\,   
p(x)^{\alpha s} +
|h|^\varepsilon p(x)^{(s+1)\varepsilon}\, |t_0|^\alpha\,   
d(x,y)^\alpha\, p(x)^{\alpha s}\bigg) \\
&\leq&  C\, |h|^\varepsilon\, d(x,y)^\alpha\, p(x)^{(s+1)(k+\varepsilon) -  
\alpha} +
C\, |h|^\varepsilon\,  d(x,y)^\alpha\, p(x)^{(s+1)(k+\varepsilon)+\alpha s}.
\end{eqnarray*}
Hence (ii) in Lemma B.1 holds with $B = C\,|h|^\varepsilon$ and $b=  
(s+1)(k+\varepsilon) + \alpha s$.
If $k=0$, the previous computation, which starts from $|q(x)-q(y)|  
\leq |\, (e_t(x) -  e_{t_0}(x)) - (e_t(y) -  e_{t_0}(y))\, |$,
yields the same conclusion. \\
By hypothesis, one can
choose $\varepsilon$ such that $\gamma' \geq  \gamma +  
\frac{(s+1)(k+\varepsilon)}{\alpha}$, and
Lemma B.1 yields for $f\in\cB_{\alpha,\beta,\gamma}$
\begin{eqnarray*}
|Kf(x)| &\leq& C\, |h|^\varepsilon\, |f|_{\alpha,\gamma}\,  p(x)^{(s+1)(k+\varepsilon) + \alpha(\gamma+1)}
\, \E[\cM^{(s+1)(k+\varepsilon) + \alpha(\gamma+1)}]  \\
&\leq&  I\, C\, |h|^\varepsilon\, |f|_{\alpha,\gamma}\, p(x)^{\alpha(\gamma'+1)}.
\end{eqnarray*}
Next, using $s+1\leq\beta$ and $d(y,x_0) \leq d(x,x_0)$ (thus $\Delta_{\alpha,\beta,\gamma}(x,y) \leq 2\, p(x)^{\alpha\gamma}  p(y)^{\alpha\beta}$) gives 
\begin{eqnarray*}
|Kf(x)-Kf(y)| &\leq&  C\,  |h|^\varepsilon\,  m_{\alpha,\beta,\gamma}(f)\, d(x,y)^\alpha\,
p(x)^{(s+1)(k+\varepsilon)}\, 2\, p(x)^{\alpha\gamma}\, p(y)^{\alpha\beta}\, 
\E[\cC^\alpha\, \cC_\lambda^{\alpha(\gamma' + \beta)}] \\[0.15cm]
&\, & +\  C\, |h|^\varepsilon \, |f|_{\alpha,\gamma}\,  d(x,y)^\alpha\, 
p(x)^{(s+1)(k+\varepsilon) + \alpha s}\,  p(y)^{\alpha(\gamma+1)}\,
\E[\cC^\alpha\, \cM^{\alpha(\gamma' + s+1)}]\\
&\leq& 2\, I\, C\,  |h|^\varepsilon\,  \|f\|_{\alpha,\beta,\gamma}\,  
d(x,y)^\alpha\,p(x)^{\alpha\gamma'}\, p(y)^{\alpha\beta}\\[0.15cm]
&\ &\ \ \ +\ I\  C\, |h|^\varepsilon \,  \|f\|_{\alpha,\beta,\gamma}\, d(x,y)^\alpha\, 
p(x)^{(s+1)(k+\varepsilon)+\alpha s}\, p(y)^{\alpha(\gamma+ 1 - \beta)}\, p(y)^{\alpha\beta}\\
&\leq& 2\, I\, C\,  |h|^\varepsilon\,  \|f\|_{\alpha,\beta,\gamma}\,  
d(x,y)^\alpha\, \Delta_{\alpha,\beta,\gamma'}(x,y) \\[0.15cm]
&\ &\ \ \ +\ I\  C\, |h|^\varepsilon \,  \|f\|_{\alpha,\beta,\gamma}\, d(x,y)^\alpha\, 
p(x)^{\alpha \gamma' + \alpha(s + 1 - \beta)}\, p(y)^{\alpha\beta}, 
\end{eqnarray*}
and we have $p(x)^{\alpha \gamma' + \alpha(s + 1 - \beta)}\, p(y)^{\alpha\beta} \leq p(x)^{\alpha \gamma'}\, p(y)^{\alpha\beta} \leq \Delta_{\alpha,\beta,\gamma'}(x,y)$ because $s+1\leq\beta$. \fdem

\noindent {\bf Lemma B.4'.} {\it For $k\in\N$, we have  
$Q_k\in\cC^1(\R,\cB_{\alpha,\beta,\gamma},\cB_{\alpha,\beta,\gamma'})$ with  
$Q_k'=Q_{k+1}$ under the conditions: $s+1\leq\beta\leq\gamma$, $\, \gamma' > \gamma + \frac{(s+1)(k+1)}{\alpha}$, and
$I = \E[\, \cM^{\alpha(\gamma' +1)} + \cC^\alpha\, \cM^{\alpha(\gamma'  
+ \beta)}\, ] < +\infty$.}

\noindent{\it Proof.} Let $t,t_0\in\R$, $h=t-t_0$, and assume $|h|\leq1$. Let
$K=Q_k(t)-Q_k(t_0)-h\, Q_{k+1}(t_0)$, and
$q(x) = (i\xi(x))^k\, \big(e_t(x) -  e_{t_0}(x) -i\, h\, \xi(x)\, e_{t_0}(x)\big)$.
For $u\in\R$, we set $\phi(u) = e^{iu}-1-iu$. Let  
$0<\varepsilon<\alpha$. We shall use the following usual inequalities   
\\[0.15cm]
\indent $\ \ \ \ \ \  \ \ \ \ \ \ \ \  \ \ \ \ \ \ \ \ \ \ \ |\phi(u)|  
\leq 2\, |u|^{1+\varepsilon},\ \  \ \
|\phi(u)-\phi(v)| \leq 2\, |u-v|\, (|u|^\varepsilon +  
|v|^\varepsilon)$. \\[0.15cm]
Writing $q(x) =  (i\xi(x))^k\, e_{t_0}(x)\, \phi\big(h\xi(x)\big)$, one  
easily gets (i) in Lemma B.1 with
$A=C\, |h|^{1+\varepsilon}$ and $a=(s+1)(k+1+\varepsilon)$. Proceeding  
as in the previous proof, one obtains
for $x,y\in E$ such that $d(y,x_0) \leq d(x,x_0)$
\begin{eqnarray*}
|q(x)-q(y)| &\leq& |\xi(x)^k-\xi(y)^k|\, |\phi(h\xi(x))| +\ |\xi(y)|^k  
\bigg|e_{t_0}(x) \phi(h\xi(x)) - e_{t_0}(y) \phi(h\xi(y))\bigg| \\
&\leq&  C\, d(x,y)^\alpha\ \, p(x)^{s+1-\alpha}\, p(x)^{(s+1)(k-1)}\,  
|h|^{1+\varepsilon}\, p(x)^{(s+1)(1+\varepsilon)} \\
&\ &\ \ \  + \  C\, p(x)^{(s+1)k}\,  \bigg(|\phi(h\xi(x)) -  
\phi(h\xi(y))| + |\phi(h\xi(y))|\, |e_{t_0}(x) - e_{t_0}(y)|\bigg)\\
&\leq&  C\, |h|^{1+\varepsilon}\, d(x,y)^\alpha\, p(x)^{s+1-\alpha +  
(s+1)(k+\varepsilon)}\ +\
C\, |h|^{1+\varepsilon}\, p(x)^{(s+1)k}\, \times \\
&\ & \ \ \ \ \ \ \ \ \ \ \ \
\bigg(|\xi(x)-\xi(y)|\,  p(x)^{(s+1)\varepsilon}  +
p(x)^{(s+1)(1+\varepsilon)} \, |t_0|^\alpha\, d(x,y)^\alpha\,  
p(x)^{\alpha s}\bigg) \\
&\leq&  C\, |h|^{1+\varepsilon}\, d(x,y)^\alpha\, p(x)^{s+1-\alpha  +  
(s+1)(k+\varepsilon)} \\
&\ & \ \ +\ C\, |h|^{1+\varepsilon}\,  d(x,y)^\alpha\, p(x)^{(s+1)k}\,
\bigg(p(x)^{s+1-\alpha + (s+1)\varepsilon}  +   
p(x)^{(s+1)(1+\varepsilon)+\alpha s}\bigg).
\end{eqnarray*}
We have $s+1-\alpha  + (s+1)(k+\varepsilon)  = (s+1)(k+1+\varepsilon)  
- \alpha \leq (s+1)(k+1+\varepsilon)+\alpha s$, and finally
one gets (ii) in Lemma B.1 with $B=C\, |h|^{1+\varepsilon}$ and $b=  
(s+1)(k+1+\varepsilon)+\alpha s$.
To prove that $Q_k\in\cC^1(\R,\cB_{\alpha,\beta,\gamma},\cB_{\alpha,\beta,\gamma'})$,  
one can then apply Lemma B.1 by proceeding exactly as in
the previous proof (replace $|h|^\varepsilon$ with  
$|h|^{1+\varepsilon}$, and $k$ with $k+1$, with $\varepsilon$
such that $\alpha\gamma'\ge \alpha\gamma+(s+1)(k+1+\varepsilon)$). 
\fdem

\noindent Now one can prove Proposition 11.7. Let us assume that 
$s+1\leq\beta\leq\gamma$ and $\gamma' > \gamma + \frac{m(s+1)}{\alpha}$, and let $\varepsilon>0$ be such that
$\gamma+{m(s+1)\over\alpha}+(2m+1)\varepsilon\leq\gamma'$. Let $I=[\gamma,\gamma']$, and for $\theta\in I$, set $\cB_\theta := \cB_{\alpha,\beta,\theta}$, $T_0(\theta)=\theta+\varepsilon$ and $T_1(\theta)=\theta+{s+1\over\alpha}+\varepsilon$.
With these choices, the conditions (0) (4) of $\cC(m)$ are obvious, the regularity conditions (1) (2) of $\cC(m)$ follow  from lemmas B.4-4', and finally Condition (3) follows from Proposition 11.4. \fdem

\noindent {\bf B.5. Proof of Proposition 11.8}.  {\it This proposition  
states that, if $s+1 < \beta \leq \gamma<\gamma'$, 
$\, I = \E[\cM^{\alpha(\gamma'+1)}+\cC^\alpha\,\cM^{\alpha(\gamma'+\beta)}] < +\infty$ and
$\E[\cC^\alpha\, \max\{\cC,1\}^{\alpha(\gamma+\beta)}]< 1$, 
then Condition (S) holds on $\cB_{\alpha,\beta,\gamma}$ if and only if $\xi$ is
non-arithmetic w.r.t.~$\cB_{\alpha,\beta,\gamma}$. If $\xi$ is nonlattice, the two previous equivalent conditions hold. }

\noindent This is a direct consequence of Propositions 5.3-4 and of the  
following lemma (for Condition~(P), see Rk.~at the end of Section~5.2). Notice that one may suppose that $\gamma'$ is fixed such that $s+1+(\gamma'-\gamma) \leq \beta$. 
Let $(\widehat{\cB},|\cdot|_{\alpha,\gamma'})$ be the Banach space of all complex-valued functions $f$ on $E$ such that 
$|f|_{\alpha,\gamma'} = \sup_{x\in E}\ \frac{|f(x)|}{p(x)^{\alpha(\gamma'+1)}} < +\infty$. 

\noindent {\bf Lemma B.5.}  {\it Under the above hypotheses, Condition  
($\widehat{K}$) of Section 5.2 is fulfilled with $\cB = \cB_{\alpha,\beta,\gamma}$ and $\widehat{\cB}$ as above defined. }

\noindent {\it Proof.} Condition (K1) holds by Proposition 11.2.
Since $|\cdot|_{\alpha,\gamma'} \leq \|\cdot\|_{\alpha,\beta,\gamma'}$, we  
have  ($\widehat{K2}$) by Lemma B.4 (case $k=0$).
To prove ($\widehat{K3}$) and ($\widehat{K4}$), let us observe that the norms  
$\|\cdot\|_{\alpha,\beta,\gamma}$ and $|\cdot|_{\alpha,\gamma'}$ may be replaced with
any equivalent norms ; of course ($\widehat{K2}$) then remains valid.
Given a real parameter $\lambda\in(0,1]$ on which conditions will be  
imposed later, let us consider on $\cB_{\alpha,\beta,\gamma}$
the norm \\[0.15cm]
\indent $\ \ \ \ \ \ \ \ \ \ \ \ \ \ \ \ \ \ \ \ \ \ \ \ \ \ \ \ \
\|f\|^{(\lambda)}_{\alpha,\beta,\gamma,\gamma'} =  
m^{(\lambda)}_{\alpha,\beta,\gamma}(f) + |f|^{(\lambda)}_{\alpha,\gamma'}$  
\\[0.15cm]
with $m^{(\lambda)}_{\alpha,\beta,\gamma}(f)$ defined in Section B.0, and 
$|f|^{(\lambda)}_{\alpha,\gamma'} := \sup_{x\in E}\ \frac{|f(x)|}{p_\lambda(x)^{\alpha(\gamma'+1)}}$.
It can be easily shown that the norms $\|\cdot\|_{\alpha,\beta,\gamma}$ and
$\|\cdot\|^{(\lambda)}_{\alpha,\beta,\gamma,\gamma'}$ are equivalent on $\cB_{\alpha,\beta,\gamma}$ (see \cite{aap} Prop.~5.2), and that the norms $|\cdot|_{\alpha,\gamma'}$ and $|\cdot|^{(\lambda)}_{\alpha,\gamma'}$ 
are equivalent on $\widehat{\cB}$. We have to establish  that, if $\lambda\in(0,1]$ is suitably chosen, then for any compact set $K_0$ in $\R$, there exist $\kappa<1$ and $C> 0$ such that: \\[0.15cm]
{\scriptsize $\bullet$} $\forall n\geq1,\ \forall f\in\cB_{\alpha,\beta,\gamma},\ \forall t\in K_0,
\ \ \|Q(t)^nf\|^{(\lambda)}_{\alpha,\beta,\gamma,\gamma'} \leq C\, \kappa^n\,  
\|f\|^{(\lambda)}_{\alpha,\beta,\gamma,\gamma'} +
C\, |f|^{(\lambda)}_{\alpha,\gamma'}$ \\[0.15cm]
\noindent {\scriptsize $\bullet$}  $\ \forall t\in K_0$, $r_{ess}(Q(t))  
\leq \kappa$.\\[0.2cm]
We have $Q(t) = K$ with $q(x) = e^{it\xi(x)}$ satisfying
Conditions (i)-(ii) of Lemma B.1 with $A=1$, $a=0$, $B = D_\lambda\,  
|t|^\alpha\ $ ($D_\lambda>0$) and $b=\alpha s$.
Let $f\in\cB_{\alpha,\beta,\gamma}$.  Because of the presence of $\gamma'$  
in the above norm,
Lemma B.1 cannot be directly applied here. However one can follow the proof of
lemma B.1 and see  
that \\[0.15cm]
\indent $\ \ \ \ \ \ \ \ \ \ \ \ \ \ \ \ \ \ \ \ \ \ \ \ \ \ \ \ \
|Q(t)f|^{(\lambda)}_{\alpha,\gamma'}
\leq \E[\cM^{\alpha(\gamma'+1)}]\, |f|^{(\lambda)}_{\alpha,\gamma'}  
\leq   I\, |f|^{(\lambda)}_{\alpha,\gamma'}$ \\[0.15cm]
  and that
for $x,y\in E$ such that $d(y,x_0) \leq d(x,x_0)$, we have by using in particular the fact that $\gamma'$ has been chosen such that $s+1 +\gamma'-\beta \leq \gamma$: 
\begin{eqnarray*}
|Q(t)f(x)-Q(t)f(y)| &\leq&  m^{(\lambda)}_{\alpha,\beta,\gamma}(f)\, d(x,y)^\alpha\,
\Delta^{(\lambda)}_{\alpha,\beta,\gamma}(x,y)\,  
\E[\cC^\alpha\, \cC_\lambda^{\alpha(\gamma+\beta)}] \\[0.15cm]
&\ & + \ D'_\lambda\, |t|^\alpha\, |f|^{(\lambda)}_{\alpha,\gamma'}\,
d(x,y)^\alpha\, p_\lambda(x)^{\alpha s}\,  p_\lambda(y)^{\alpha(\gamma'+1)}\,
\E[\cC^\alpha\, \cM^{\alpha(\gamma'+s+1)}] \\
&\leq&  m^{(\lambda)}_{\alpha,\beta,\gamma}(f)\, d(x,y)^\alpha\,
\Delta^{(\lambda)}_{\alpha,\beta,\gamma}(x,y)\, \E[\cC^\alpha\, \cC_\lambda^{\alpha(\gamma+\beta)}]\\[0.15cm]
&\ & + \ I\, D'_\lambda\, |t|^\alpha\, |f|^{(\lambda)}_{\alpha,\gamma'}\,
d(x,y)^\alpha\, p_\lambda(x)^{\alpha s}\,  p_\lambda(y)^{\alpha(\gamma'+1-\beta)}\, p_\lambda(y)^{\alpha\beta}\, , 
\end{eqnarray*}
with $p_\lambda(x)^{\alpha s}\, p_\lambda(y)^{\alpha(\gamma'+1-\beta)} \leq 
p_\lambda(x)^{\alpha(s+1 +\gamma'-\beta)} \leq p_\lambda(x)^{\alpha\gamma}$. Thus 
$$\frac{|Q(t)f(x)-Q(t)f(y)|}{d(x,y)^\alpha\, \Delta^{(\lambda)}_{\alpha,\beta,\gamma}(x,y)} \leq m^{(\lambda)}_{\alpha,\beta,\gamma}(f)\, \E[\cC^\alpha\, \cC_\lambda^{\alpha(\gamma+\beta)}]\,
+ \, I\, D'_\lambda\, |t|^\alpha\, |f|^{(\lambda)}_{\alpha,\gamma'}.$$
Besides, by Lebesgue's theorem, we have  $\kappa:=\E[\cC^\alpha\,  
\cC_\lambda^{\alpha(\gamma+\beta)}] < 1$ for sufficiently small $\lambda$.
The previous estimate then easily gives the desired Doeblin-Fortet  
inequalities. \\ 
Since the canonical embedding from
$\cB_{\alpha,\beta,\gamma}$ into $\widehat{\cB}$ is compact
(this easily follows from Ascoli's theorem, see \cite{aap} Lemma 5.4),  
the property $r_{ess}(Q(t)) \leq \kappa$ is then a consequence of
\cite{hen2}.  \fdem


\newpage 


\begin{thebibliography}{99}

\bibitem{als}
{\sc Alsmeyer G.}
{\em On the Harris recurrence of iterated random Lipschitz functions and related convergence rate results.}  
J. of Th.~Probab., Vol.~16, N.~1 (2003). 

\bibitem{aar2}
{\sc Aaronson J.~and Denker M..}
{\em Local limit theorem for partial sums of stationary sequences generated by Gibbs-Markov maps.}   
Stochastic Dynamics 1(2), pp.~193-237, 2001. 
 
\bibitem{aar3}
{\sc Aaronson J.~and Denker M..}
{\em A Local limit theorem for stationary processes in the domain of attraction of a normal distribution.}   
In N. Balakrishnan, I.A. Ibragimov, and V.B. Nevzorov, editors, Asymptotics methods in probabability and statistics 
with applications. Papers from the international conference, St. Petersburg, Russia, 1998, pp.~215-224. Birkh\"auser, 2001. 

\bibitem{bapei}
{\sc Babillot M.~and  Peign\'e M.}
{\em Asymptotic laws for geodesic homology on hyperbolic manifolds with cusps.} 
Bull.~Soc.~Math.~France, 134, Fasc.~1 (2006), 119-163. 

\bibitem{bal}
{\sc Baladi V.}
{\em Positive transfer operators and decay of correlations.}  
Advanced Series in Nonlinear Dynamics 16, World Scientific (2000). 

\bibitem{gouezel}
{\sc B\'alint P., Gou\"ezel S.}
{\em Limit theorems in the stadium billiard.} 
Communications in Mathematical Physics {\bf 263}, 451-512 (2006). 

\bibitem{benda}
{\sc Benda M.}
{\em A central limit theorem for contractive stochastic dynamical systems.}  
J.~App.~Prob. {\bf 35} (1998) 200-205.

\bibitem{bertail}
{\sc Bertail P., Cl\'emen\c con S.}
{\em Edgeworth expansions of suitably normalized sample mean statistics for atomic Markov chains.}  
Prob.~Theory Relat.~Fields 130, 388-414 (2004). 

\bibitem{Billingsley}
{\sc Billingsley P.}
{\em Convergence of probability measures.}
Wiley \& Sons, New-York (1968).

\bibitem{BoattoGolse}
{\sc Boatto S., Golse F.}  
{\em Diffusion approximation of a Knudsen gaz model : dependence of the diffusion constant upon a boundary condition.}
Asymptotic.~Anal., Vol.~31, no 2, 93-111 (2002). 

\bibitem{bolt}
{\sc Bolthausen, E.}  
{\em The Berry-Esseen theorem for strongly mixing Harris recurrent Markov chains.}
Z. Wahrscheinlichkeitstheorie verw. Gebiet {\bf 60} 283-289 (1982).

\bibitem{bre}
{\sc Breiman L.}
{\em Probability}
Classic in Applied Mathematics, SIAM, 1993. 

\bibitem{broi}
{\sc Broise A., Dal'bo F., Peign\'e M.}
{\em \'Etudes spectrales d'op\'erateurs de transfert et applications.} Ast\'erisque 238 (1996).

\bibitem{chazottes-gouezel}
{\sc  Chazottes J.-R., Gou\"ezel S.} 
{\em On almost-sure versions of classical limit theorems for dynamical systems.} 
Probability Theory and Related Fields 138:195-234, 2007. 

\bibitem{chen}
{\sc Chen X.}
{\em Limit theorems for functionals of ergodic Markov chains with general state space. Memoirs of the American Mathematical Society, 
139. (1999)} 

\bibitem{datta}
{\sc Datta S. McCormick W.P.}
{\em On the first-order Edgeworth expansion for a Markov chain. } 
Journal of Multivariate Analysis {\bf 44}, 345-359 (1993).

\bibitem{ded-mer-rio}
{\sc Dedecker J., Merlevède F., Rio E.} 
{\em Rates of convergence for minimal metrics in the
central limit theorem under projective criteria.} 
To appear in Electronic Journal of Probability. 

\bibitem{ded-rio}  
{\sc Dedecker J., Rio E.}
{\em On Esseen's mean central limit theorem for dependent sequences.} 
Ann. Inst. H. Poincare, probab. et statist. 2008 Vol. 44, No 4, 693-726.

\bibitem{diaco}
{\sc Diaconis P.~and  Freedman D.}
{\em Iterated random functions. } 
SIAM Rev. {\bf 41}, no.~1, 45-76 (1999).

\bibitem{Dudley}
{\sc Dudley R. M.}
{\em Real analysis and probability. } 
Wadsworth \& Brooks Cole Math. series,
Pacific Grove CA (1989).

\bibitem{duf}
{\sc Duflo M.}
{\em Random Iterative Models.} 
Applications of Mathematics, Springer-Verlag Berlin Heidelberg (1997). 

\bibitem{ds}
{\sc Dunford N. and Schwartz J.T.}
{\em Linear Operators.}
Part. I : General Theory. New York : Wiley 1958. 

\bibitem{dur}
{\sc Durrett R.}
{\em Probability: theory and examples.}
Wadsworth Brooks (1991).

\bibitem{fel}
{\sc Feller W.}
{\em An introduction to probability theory and its applications, Vol.~II.}
John Wiley and Sons, New York (1971).

\bibitem{fuhlai}
{\sc Fuh C.-D, Lai T.L.}
{\em Asymptotic expansions in multidimensional Markov renewal theory and first passage times for Markov random walks. }
Adv.~in Appl.~Probab. {\bf 33} 652-673 (2001).

\bibitem{gharib}
{\sc Gharib M.}
{\em A uniform estimate for the rate of convergence in the multidimensional central limit theorem for homogeneous Markov chains.}
Internat. J.~Math.~Math.~Sci., Vol.~19, No 3 (1996) 441-450.

\bibitem{gong-wu}
{\sc Gong F.Z., Wu L.M.}
{\em Spectral gap of positive operators and applications.} 
J. Math. Pures Appl. 85 (2006) pp.~151-191. 

\bibitem{gor}
{\sc Gordin M.I., Lifsic B.A.}
{\em On the central limit theorem for stationary Markov processes.} 
Soviet Math.~Dokl.~19 No 2 (1978) 392-394. 

\bibitem{gotze}
{\sc G\"otze F., Hipp C.}
{\em Asymptotic expansions for sums of weakly dependent random vectors.} 
Z.~Wahrscheinlichkeitstheorie verw.~Gebiete, {\bf 64}, 211-239 (1983). 

\bibitem{gou}
{\sc Gou\"ezel S.}
{\em Central limit theorem and stable laws for intermittent maps.} 
Probability Theory and Related Fields {\bf 128}, 82-122 (2004) 


\bibitem{seb-08}
{\sc Gou\"ezel S.} 
{\em Necessary and sufficient conditions for limit theorems in Gibbs-Markov maps. }
Preprint (2008). 

\bibitem{gouliv}
{\sc Gou\"ezel S., Liverani C.}
{\em Banach spaces adapted to Anosov systems. } 
Ergodic Theory Dyn.~Syst. {\bf 26}, 189-217 (2006). 

\bibitem{guibourg}
{\sc Guibourg D.}
{\em Th\'eor\`eme de renouvellement pour  cha\^{\i}nes de Markov g\'eom\'triquement ergodiques. 
Applications aux mod\`eles it\'eratifs Lipschitziens.} 
C.~R.~Acad.~Sci.~Paris, Ser.~I 346 (2008) 435-438.  

\bibitem{guiher}
{\sc Guibourg D., Herv\'e L.}
{\em  A renewal theorem for strongly ergodic Markov chains in dimension $d\geq3$ and in the centered case. } 
Preprint (2009). To appear in Potential Analysis. 

\bibitem{gui}
{\sc Guivarc'h Y.}
{\em Application d'un th\'eor\`eme limite local \`a la transcience et \`a la r\'ecurrence de marches al\'eatoires. } 
Lecture Notes in Math.~Springer, 301-332 (1984) 

\bibitem{guihar}
{\sc Guivarc'h Y., Hardy J.}
{\em Th\'eor\`emes limites pour une classe de cha\^{\i}nes de Markov et applications aux diff\'eomorphismes d'Anosov.} 
Ann.~Inst.~Henri Poincar\'e, Vol.~24, No 1, p.~73-98 (1988). 

\bibitem{gui-lejan}
{\sc Guivarc'h Y., Le Jan Y.}
{\em Asymptotic winding of the geodesic flow on modular surfaces and continuous fractions.}
 Ann.~scient.~Ec.~Norm.~Sup, 4, 26, (1993), 23-50. 

\bibitem{gui-lepage}
{\sc Guivarc'h Y., Le Page E.}
{\em On spectral properties of a family of transfer operators and convergence to stable laws for affine random walks.}
Ergodic Theory Dynam. Systems, 28 (2008) no.~2, pp.~423-446. 

\bibitem{hen1}
{\sc Hennion H.}  
{\em  D\'erivabilit\'e du plus grand exposant caract\'eristique des produits de matrices al\'eatoires 
ind\'ependantes \`a coefficients positifs. } 
Ann.~Inst.~Henri Poincar\'e, Vol.~27, No 1 (1991) p.~27-59. 

\bibitem{hen2}
{\sc Hennion H.}
{\em Sur un th\'eor\`eme spectral et son application aux noyaux lipchitziens.} 
Proceeding of the A.M.S vol. 118 No 2 (1993) 627-634.

\bibitem{hen4}
{\sc Hennion H.}
{\em Quasi-compactness and absolutely continuous kernels.} 
Probab.Theory Related Fields 139 (2007) pp.~451-471. 
\bibitem{hulo}
{\sc Hennion H., Herv\'e L.}
{\em Limit theorems for Markov chains and stochastic properties of dynamical systems by quasi-compactness.} 
Lecture Notes in Mathematics No 1766, Springer (2001). 

\bibitem{aap}
{\sc Hennion H., Herv\'e L.}
{\em Central limit theorems for iterated random lipschitz mappings.} 
Annals of Proba.~Vol.~32 No.~3A (2004) pp.~1934-1984.  

\bibitem{hh-lh}
{\sc Hennion H., Herv\'e L.}
{\em Stable laws and products of positive random Matrices.}
J.~Theor.~Probab., 21, no.~4 (2008) pp.~966-981. 

\bibitem{ihp1}
{\sc Herv\'e L.}
{\em Th\'eor\`eme local pour cha\^{\i}nes de Markov de probabilit\'e de transition quasi-compacte. Applications aux 
cha\^{\i}nes $V$-g\'eom\'etriquement ergodiques et aux mod\`eles it\'eratifs.} 
Ann.~I.~H.~Poincar\'e - PR 41 (2005) 179-196. 

\bibitem{ihp2}
{\sc Herv\'e L.}
{\em Vitesse de convergence dans le th\'eor\`eme limite central pour cha\^{\i}nes de Markov fortement ergodiques. } 
Ann.~Inst.~H.~Poincar\'e Probab. Statist.~{\bf 44} No.~2 (2008) 280-292. 

\bibitem{jlv}  
{\sc Hervé L., Ledoux J., Patilea V.} 
{\em A Berry-Esseen theorem on $M$-estimators for geometrically ergodic Markov chains.} 
Preprint (2009). 

\bibitem{itm}
{\sc C.T. Ionescu-Tulcea, G. Marinescu.}
{\em Th\'eor\`eme ergodique pour des classes d'op\'erations non compl\`etement 
continues.} 
Ann.~of Maths.~52 1 (1950) 140-147. 

\bibitem{jan}
{\sc Jan C.}
{\em Vitesse de convergence dans le TCL pour des processus associ\'es \`a des syst\`emes dynamiques et aux produits de 
matrices al\'eatoires. } 
Th\`ese de doctorat (2001) I.R.M.A.R, Universit\'e de Rennes I. 

\bibitem{jensen}
{\sc Jensen J.L.}
{\em Asymptotic Expansions for strongly mixing Harris recurrent Markov chains. } 
Scand.~J.~Statist.~{\bf 16} 47-63 (1989). 

\bibitem{jones}  
{\sc Jones G.L.}
{\em On the Markov chain central limit theorem.} 
Probability surveys, Vol.~1 (2004) 299-320.  

\bibitem{keli}
{\sc Keller G., Liverani C.}
{\em Stability of the Spectrum for Transfer Operators.} 
Ann.~Scuola Norm.~Sup.~Pisa.~CI.~Sci.~(4) Vol.~XXVIII (1999) 141-152. 

\bibitem{kontmey}
{\sc Kontoyiannis, I., Meyn, S.P.} 
{\em Spectral theory and limit theorems for geometrically ergodic Markov processses.}  
Annals of Applied Probability {\bf 13} (2003) 304-362.

\bibitem{lep82}
{\sc Le Page E.} 
{\em Th\'eor\`emes limites pour les produits de matrices al\'eatoires.} 
Springer Lecture Notes, 928 (1982) 258-303. 

\bibitem{lep83}
{\sc Le Page E.} 
{\em Th\'eor\`emes de renouvellement pour les produits de matrices al\'eatoires. Equations aux diff\'erences al\'eatoires.} 
S\'eminaires de Rennes (1983). 

\bibitem{lep89}   
{\sc Le Page E.}
{\em R\'egularit\'e du plus grand exposant caract\'eristique des produits de matrices al\'eatoires ind\'ependantes et applications. } 
Ann.~Inst.~Henri Poincar\'e, Vol.~25, No 2 (1989) pp.~109-142.

\bibitem{liverani}
{\sc Liverani C.} 
{\em Invariant measure and their properties. 
A functional analytic point of view}
Dynamical Systems. Part II: Topological Geometrical 
and Ergodic Properties of Dynamics. 
Pubblicazioni della Classe di Scienze, 
Scuola Normale Superiore, Pisa. 
Centro di Ricerca Matematica "Ennio De Giorgi" : 
Proceedings (2004).
\bibitem{malin}  
{\sc Malinovskii V.K..}
{\em Limit theorems for Harris Markov chains I. } 
Theory Prob.~Appl.~{\bf 31}, 269-285 (1987).  

\bibitem{mey}  
{\sc S.P.~Meyn, R.L.~Tweedie.}
{\em Markov chains and stochastic stability.} 
Springer Verlag, New York, Heidelberg, Berlin (1993). 

\bibitem{mira}  
{\sc Milhaud X., Raugi A.}
{\em Etude de l'estimateur du maximum de vraisemblance dans le cas d'un
processus auto-r\'egressif: convergence, normalit\'e asymptotique, vitesse de convergence.} 
Ann.~Inst.~H.~Poincar\'e Vol.~25 No 4 (1989) 383-428.

\bibitem{nag1}
{\sc Nagaev S.V.}
{\em Some limit theorems for stationary Markov chains.} 
Theory of probability and its applications 11 4 (1957) 378-406.  

\bibitem{nag2}  
{\sc Nagaev S.V.}
{\em More exact statements of limit theorems for homogeneous Markov chains.} 
Theory of probability and its applications 6 1 (1961) 62-81.

\bibitem{pei}    
{\sc Peign\'e M.} 
{\em Iterated function schemes and spectral decomposition of the associated Markov operator. }
S\'eminaires de Probabilit\'e de Rennes (1993). 

\bibitem{FPAAP}
{\sc P\`ene F.}
{\em Rate of convergence in the multidimensional CLT for
stationary processes. Application to the Knudsen gas 
and to the Sinai billiard}
Ann.~Appl.~Probability, 15 (4), pp.~2331-2392 (2005).

\bibitem{rab-wolf}
{\sc Räbiger F., Wolff M.P.H.}
{\em On the approximation of positive operators and the behaviour of the spectra of the appoximants. } 
Integr. equ. oper. theory. 28 (1997) pp.~72-86. 

\bibitem{rosen}
{\sc Rosenblatt M. }
{\em Markov processes. Structure and asymptotic behavior. } 
Springer-Verlag.~New York (1971). 

\bibitem{Rotar70}
{\sc Rotar V. I. }
{\em A non-uniform estimate for the convergence speed in the multidimensional
central  theorem. } 
Theory Prob.~Applications, 15, pp.~630-648 (1970). 

\bibitem{rou}
{\sc Rousseau-Egele J.}
{\em Un th\'eor\`eme de la limite locale pour une classe de transformations dilatantes. } 
Annals of Proba., 11, 3, pp.~772-788 (1983). 


\bibitem{sev}
{\sc Seva M.}
{\em On the local limit theorem for non-uniformly ergodic Markov chains.} 
J.~Appl.~Prob.~{\bf 32}, 52-62 (1995). 

\bibitem{wu}
{\sc Wu L.M.}
{\em Essential spectral radius for Markov semigroups (I): discrete time case.} 
Probab.Theory Related Fields, {\bf 128}, pp.~255-321 (2004). 

\bibitem{Yurinskii}
{\sc Yurinskii V. V.}
{\em A smoothing inequality for estimates of
the Levy-Prokhorov distance,}
Theory Probab.~Appl. {\bf 20}, pp.~1-10 (1975);
translation from Teor. Veroyatn. Primen.
{\bf 20}, pp.~3-12 (1975).


\end{thebibliography}
\end{document}